\documentclass[twoside,a4paper,12pt]{article}
\usepackage[OT1]{fontenc}
\usepackage[english]{babel}
\usepackage{amsthm,amsmath,amssymb}
\usepackage{authblk}
\usepackage[left=3cm,right=2cm,top=3.5cm,bottom=3cm,asymmetric]{geometry}
\usepackage{geometry}
\usepackage{graphicx}
\usepackage{url}
\usepackage[hidelinks,colorlinks=false,linkcolor=blue,citecolor=blue]{hyperref}
\graphicspath{ {./images/} }
\usepackage[all]{xy}
\usepackage{float}
\usepackage[section]{placeins}
\usepackage[table,xcdraw]{xcolor}
\usepackage[switch,columnwise]{lineno}
\usepackage{slantsc}
\usepackage{enumitem}
\usepackage{lastpage}
\usepackage{longtable}
\usepackage[utf8]{inputenc}
\usepackage[T1]{fontenc}
\usepackage{lmodern} 
\usepackage{geometry}
\usepackage{amsmath}
\usepackage{fouriernc}
\usepackage{etoolbox}
\usepackage{cases}
\usepackage{makecell}
\usepackage[bottom]{footmisc}
\usepackage{fancyhdr}
\usepackage[
singlelinecheck=false
]{caption}
\title{Integer Factorization By Sieving The Delta}
\author{Vishal Mudgal \footnote {17 Simei Street 4, Double Bay Residences, Singapore-529881. Email: vishal.mudgal@gmail.com}}
\date{}
\usepackage{everypage}
\usepackage{caption}
\newlength\newtop
\makeatletter
\let\runauthor\@author
\begin{document}
\normalfont
\maketitle
\section{Abstract}\
Let \textit {$n = p\times q$} ($p < q$) and $\Delta = | p - q |$, where $p,q$ are odd integers, then, it is hypothesized that factorizing this composite $n$ will take O(1) time once the steady state value is reached for any $\Delta$ in $zone_0$ of some observation deck ($od$) with specific dial settings. We also introduce a new factorization approach by looking for $\Delta$ in different $\Delta$ sieve zones. Once $\Delta$ is found and $n$ is already given, one can easily find the factors of this composite $n$ from any one of the following quadratic equations: $p^2 + p\Delta - n = 0$ or $q^2 - q\Delta - n = 0$. The new factorization approach does not rely on congruence of squares or any special properties of $n$, $p$ or $q$ and is only based on sieving the $\Delta$. In addition, some other new factorization approaches are also discussed. Finally, a new trapdoor function is presented which is leveraged to encrypt and decrypt a message with different keys.
\tableofcontents
\section{Introduction}\
If $n = p \times q$ where $p,q$ are some large primes, then no non-quantum algorithm exists today that can either find $p$ or $q$ in polynomial time when only $n$ is given. This problem of factorizing the product of two large primes has kept great minds busy for quite some time, especially since 1977, when Rivest-Shamir-Adleman \cite{rsa} leveraged this problem in building a public-key cryptosystem, famously known as RSA, that is widely used for secure data transmission. However, in 1994, Peter Shor \cite{shor} published an algorithm that can solve the integer factorization problem in polynomial time by using a universal quantum computer. When in future, a large ideal universal quantum computer with enough qubits is able to operate efficiently, then RSA scheme will no longer remain secure. Two key ideas have remained centerstage when attempting to solve the factorization problem. Variety of different methods and approaches have been developed around these two central ideas.
\begin{enumerate}
\item \textbf{Congruence of squares method}
\\
In the 17th century, Pierre de Fermat \cite{fermat1600s} applied the difference-of-squares method to factorization. The idea is to represent an odd integer as the difference of two squares, i.e. $n = a^2 - b^2$. If neither factor equals one, then $n$ can be factored as $(a+b)(a-b)$. This method works well when the two factors are close to each other but in worst case, this method can be far worse than trial division. However, this idea was further expanded by Maurice B. Kraitchik in the 1920s, who reasoned that instead of finding two integers $a$ and $b$, such that, $n = a^2 - b^2$, it is sufficient to find two integers $a$ and $b$ with $a^2 - b^2$ equal to multiple of n, i.e. $a^2 \equiv b^2 \pmod n$ and $a \not\equiv \pm b \pmod n$. The idea of factorization by congruence of squares was born. Most general-purpose factorization algorithms are based on this idea.
\\ \\
 Continued fraction factorization (CFRAC) by Derrick H. Lehmer and Ralph E. Powers  \cite{cfrac} and implemented as a computer algorithm by Michael A. Morrison and John Brillhart \cite{cfrac1}, Square forms factorization by Daniel Shanks \cite{squfof}. Linear sieve by Richard C. Schroeppel \cite{ls}. Random squares method by John D. Dixon \cite{dixon}. Quadratic sieve by Carl B. Pomerance \cite{qs}. The multiple polynomial quadratic sieve by Robert D. Silverman \cite{mpqs}. The special number field sieve by John M. Pollard, and even though this is a special purpose factorization method, it laid the foundations for the general number field sieve (GNFS), making GNFS the fastest non-quantum algorithm. The details of how Pollard's original idea, which was suited for numbers of special form was extended to be applied to any number in general, is described in the book titled "The development of the number field sieve" \cite{nfs}. Finally, some Lattice based factorization methods have also been explored by many researchers to solve the factorization problem.
\\ \\
The aforementioned factorization methods may use many other different ideas and concepts from analytical number theory, analytical algebraic theory or other areas of mathematics, however, achieving congruence of squares remains a key underlying goal to find the nontrivial factors of the semiprime $n$.
\item \textbf{Special properties of semiprime or any one of its unknown factors}
\\
Numbers come in different size and form and sometimes this size and form have special properties, which is leveraged for factorization.  Pollard's rho algorithm by John M. Pollard \cite{pollard1975} is useful when number to be factored has small factors and Pollard's $p - 1$ algorithm, again by John M. Pollard \cite{pollard_1974} is used when the number preceding one of the factors, p-1, is powersmooth. Williams's $p + 1$ algorithm by Hugh C. Williams \cite{williams1982} is used in scenarios when the number to be factored contains one or more prime factors $p$ such that $p + 1$ is smooth, i.e. $p + 1$ only contains small factors and Elliptic-curve factorization method (ECM) by Hendrik W. Lenstra Jr. \cite{lenstra1987} is most suitable for finding small factors not exceeding 50 to 60 digits.
\end{enumerate}
There are many other different and creative ideas that have been applied to factorization, but the above two approaches have remained the front-runners for many years. A lot of research has already been carried out around these two key ideas and while a lot more is still left to be discovered, we thought to experiment with a different underlying approach. 
\\
After a lot of failed experiments and toying with new structures, the inspiration came from nature and especially from the principle of flow of energy. For energy to flow between any two points, a difference of some sorts is required between these points - potential difference, gravitational difference, pressure difference etc. Applying this newly found inspiration to factorization, where, $n = p \times q$ and $p,q$ are some large primes, the only  reasonable potential/gravitational/pressure difference equivalent was observed in:
\begin{equation}
\Delta = |p-q|
\end{equation}
And we embarked on a long but beautiful journey to comprehend this esoteric "$\Delta$". With only $n$ given as an input, it is quite fascinating to observe this property of numbers where one can make the $\Delta$ appear by only performing fundamental arithmetic operations. Along the way, this also led us to launch a systemic study of:
\begin{equation}
\sum = p + q
\end{equation}
Leading us to find the beautiful equilibrium between the two and laid the foundations for a new trapdoor function that can be used in encryption and decryption of data and also opened up another factorization avenue at the same time.
\begin{equation}
\Delta \rightleftharpoons \sum
\end{equation}
\\
We could not find any literature in the public domain where the ideas and observations presented in this paper have been applied to factorization or for encryption/decryption of data, as we also found a new trapdoor function. We spent a lot of time in finding any prior work having similarities with our research and findings, but to the best of our efforts, we couldn't find anything. 
\\ \\
We present our findings focussed on this $\Delta_{|p-q|}$,  $\sum_{p+q}$ and $\Delta_{|p-q|} \rightleftharpoons \sum_{p+q}$, and some other insights that we think we have been blessed with during the course of this journey.
\section{Notation and setup}\
$\Delta=| p - q |$: difference between any two integers and it's resultant absolute value
\\
$p = 2k+1$ or $p=2k$ for $k \in \mathbb{N}$									
\\
$q = p + \Delta$							
\\
$n = p\times q$														
\\
$\lfloor \sqrt{n} \rfloor$: the floor function of $n$ 
\\
$\Delta_{|p-q|}$: delta series (all combinations of $p$ and $q$ with respective $|p-q|$ value)									
\\
$\sum_{p+q}$: sum series (all combinations of $p$ and $q$ with respective $p+q$ value)									
\\
$\Delta_{|p-q|} \rightleftharpoons \sum_{p+q}$: equilibrium state between delta and sum series
\\
$dial_1$:\ relative distance from  $\lfloor \sqrt{n} \rfloor$:	
\begin{subnumcases}{dial_1=}
   a_1 \in \mathbb{Z}_0 & \text{for $\lfloor \sqrt{n} \rfloor = 2k,\ k \in  \mathbb{N}$} \label{eqn:nas1a} \
  \\ [6pt]
   a_2 \in \mathbb{Z}_0 & \text{for $\lfloor \sqrt{n} \rfloor = 2k+1,\ k \in  \mathbb{N}$} \label{eqn:nas1b} \
\end{subnumcases}
$dial_2$: relative distance from $\lfloor \sqrt{n} \rfloor + dial_1$ (or $d_1$):
\begin{subnumcases}{dial_2=}
  v_1 \in \mathbb{Z} & \text{for $\lfloor \sqrt{n} \rfloor + dial_1= 2k,\ k \in  \mathbb{N}$} \label{eqn:nas2a} \
  \\ [6pt]
  v_2  \in \mathbb{Z} & \text{for $\lfloor \sqrt{n} \rfloor + dial_1 = 2k+1,\ k \in  \mathbb{N}$} \label{eqn:nas2b} \
\end{subnumcases}
\\
$dial_P$: dial pair and represented as \{$dial_1,dial_2$\} 
\clearpage\
\\
$dial\ settings$: can be represented in any of the below formats:
\begin{enumerate}[noitemsep]
\item As a dial pair. i.e. $dial_P  = \{dial_1,dial_2\} = \{a_1, a_2, v_1, v_2\}$
\item As a matrix:
\begin{table}[!htbp]
\captionsetup{labelformat=empty}
\phantom{800000000em} 
\begin{tabular}{|c|c|}
\hline 
$a_1$ & $a_2$ \\ 
\Xhline{0.5\arrayrulewidth}
$v_1$ & $v_2$ \\										 						
\Xhline{0.5\arrayrulewidth}
\end{tabular}
\end{table}
\end{enumerate}
$d_1$ =  $\lfloor \sqrt{n} \rfloor + dial_1$
\\
$d_2$ =  $d_1 + dial_2$
\\
$od$: observation deck. Different observation decks will be referred with a numerical subscript value
\\ \\
$od_1$ =  $d_1^2 - n$
\\ \\
$od_2$ =  $d_2^2 - n$
\\ \\
$od_3$ =  $od_2 - od_1$
\\
$od_4$ =  $od_1 + od_2$
\\
$od_5$ =  $od_1 + od_2 + od_3 + od_4$
\\
$od_6$ = $\sqrt{(n \times v_1^2) + (od_1\times od_2)}$   for  $v_1 = v_2 = 2$
\\ \\
$df_x$: difference between consecutive $od$ values. Difference between different observation decks will be referred with a numerical subscript value. For e.g. $df_6$ represents difference between consecutive $od_6$ values
\\ \\
$id$: sequence number starting from 1
\\ \\
$Steady\ State\ Value$: the value which can be expressed as a function of $\Delta$ and some constant $k$ and doesn't change for a particular observation deck with a specific dial setting, for any $\Delta$
\\ \\
$\Delta\ sieve\ zone$: this is the zone where $\Delta$ can be sieved. It is hypothesized that there are infinite $\Delta$ sieve zones
\\ \\
$\Delta\ sieve\ coverage$: sum total of numbers that can be sieved in a particular $\Delta$ sieve zone with specific dial settings
\\ \\
$Switchover\ zone$: the first row of a zone when the steady state value is first observed is defined as the switchover zone ($soz$) for that observation deck
\\ \\
$Reflection\ over\ \{X,Y\}\ or\ ro\{X,Y\}$: when value of an observation deck repeats over an interval from some convergence points/values, it is referred as $ro\{X,Y\}$
\\ \\
\section{Steady State Value}\
For some fixed dials and known $\Delta$, we aim to find nontrivial observation decks, such that:
\begin{equation}
\begin{gathered}
\sum\limits_{i=x}^{k} {a}_i {df}_i = 0$ for $x, k \in \mathbb{N}$, $a \in \mathbb{Z}_0 \text {; and } \\
\sum\limits_{i=x}^{k} {a}_i {od}_i \neq 0$ for $x, k \in \mathbb{N}$, $a \in \mathbb{Z}_0  \label{eqn:dfodsum} \
\end{gathered}
\end{equation}
When this is achieved, it is hypothesized, that a steady state value will be yielded in a specific $\Delta$ sieve zone of an observation deck and this steady state will remain constant for \textbf{any} $\Delta$ in the same observation deck.
There can be many different steady states in different (or same) observation decks and the steady state value as a function of $\Delta$ and constants is expressed as:
\begin{equation}
\begin{gathered}
c_1 \Delta^2 + c_2 \Delta +  k \label{eqn:ssvgen_main_1a}
\end{gathered}
\end{equation}
Note: $c_1, c_2\in \mathbb{Q}$, $k \in \mathbb{Z}$ and $c_1,c_2,k$ are constants for a given observation deck with specific dial settings. \\ \\
Now, finding $\sum\limits_{i=x}^{k} {a}_i {df}_i = 0$ is a subset sum problem (SSP), where target sum T = 0 is required. However, in some edge cases, equation (\ref{eqn:dfodsum}) may yield a trivial solution. The verification process is simple, as any steady state value that is yielded for known $\Delta$ in any observation deck can be verified with "known $\Delta + 4$". Now, if steady state value obtained for known $\Delta$ doesn't hold for "known $\Delta + 4$", then the solution of equation (\ref{eqn:dfodsum}) was trivial, otherwise, this steady state will remain constant for any $\Delta$ in the same observation deck as found for known $\Delta$.
\\
\\
This steady state property is quite interesting, it's like tuning a $\Delta$ sieve zone for a particular observation deck with known $\Delta$ and sieving any $\Delta$ through it (as long as the observation deck yields the corresponding $\Delta$ sieve zone). The zone where steady state is observed is referred as a $\Delta$ sieve zone for that observation deck. There can be many $\Delta$ sieve zones ($zone_0, zone_1,..., zone_x$), with different steady states in the same observation deck.  
\\
\\
Below equations will help in finding the steady state value with specific dial settings in $zone_0$ of $od_4$ and $od_2$. When a steady state value is yielded in $zone_0$ of an observation deck, then it is hypothesized that this steady state will continue until $\infty$. Let's see this for $od_4$ first, when:
\\ \\
$dial_P = \{0,-1,2,2\}$; $\Delta = 4k, p = 2j+1, \ k, j \in \mathbb{N}$; or \\
$dial_P = \{0,-1,2,2\}$; $\Delta = 4k+2, p = 2j, \ k, j \in \mathbb{N}$; or \\
$dial_P = \{-1,0,2,2\}$; $\Delta = 4k+2, p = 2j+1, \ k, j \in \mathbb{N}$; or \\
$dial_P = \{-1,0,2,2\}$; $\Delta = 4k, p = 2j, \ k, j \in \mathbb{N}$
\begin{equation}
SteadyStateValue_{_{od_{_4}}} =  \dfrac{\Delta^2}{2} + 2 \label{eqn:ssv1a}
\end{equation}
\clearpage\
\\
Likewise, for $od_2$, when: \\ \\
$dial_P = \{0,-1,2,2\}$; $\Delta = 4k+2, p = 2j+1, \ k, j \in \mathbb{N}$; or \\
$dial_P = \{0,-1,2,2\}$; $\Delta = 4k, p = 2j, \ k, j \in \mathbb{N}$; or \\
$dial_P = \{-1,0,2,2\}$; $\Delta = 4k, p = 2j+1, \ k, j \in \mathbb{N}$; or \\
$dial_P = \{-1,0,2,2\}$; $\Delta = 4k+2, p = 2j, \ k, j \in \mathbb{N}$ \\
\begin{equation}
SteadyStateValue_{_{od_{_2}}} =  \left(\dfrac{\Delta}{2}\right)^2 \label{eqn:ssv1b}
\end{equation}
Also, when \textit {p < q}, then the value of \textit {p} on first occurrence of steady state value in $zone_0$ (either $od_2$ or $od_4$) and when $dial_P = \{0,-1,2,2\}$ is given as:
\begin{subnumcases}{p=}
  \left(2 \times \left(\dfrac{\Delta - 4}{4}\right) \left(\dfrac{\Delta - 4}{4} + 1\right)\right)+ 1 & \text{for $\Delta =  4k, \ k \in \mathbb{N};\ p = 2j+1, \ j \in \mathbb{N}$} \label{eqn:ssv2a} \
  \\ [8pt]
  \left(\dfrac{\Delta - 6}{4} + 1\right)^2 - \left(\dfrac{\Delta - 6}{4}\right) & \text{for $\Delta =  4k+2, \ k \in \mathbb{N};\ p = 2j+1, \ j \in \mathbb{N}$} \label{eqn:ssv2b} \
  \\ [8pt]
  \left(\dfrac{\Delta - 4}{4}\right)^2 & \text{for $\Delta =  4k, k = 2m+1, \ k,m \in \mathbb{N}; p = 2j, \ j \in \mathbb{N}$} \label{eqn:ssv2c} \
  \\ [8pt]
  \left(\dfrac{\Delta - 4}{4}\right) \left(\dfrac{\Delta - 4}{4} + 1\right) - \left(\dfrac{\Delta - 4}{4} - 1\right)  & \text{for $\Delta =  4k, k = 2m, \ k,m \in \mathbb{N}; p = 2j, \ j \in \mathbb{N}$} \label{eqn:ssv2d} \
  \\ [8pt]
  \left(2 \times \left(\left(\dfrac{\Delta - 6}{4} + 1\right)^2 - 1\right)\right) + 2  & \text{for $\Delta =  4k+2,\ k \in \mathbb{N};\ p = 2j, \ j \in \mathbb{N}$} \label{eqn:ssv2e} \
\end{subnumcases}
\\
Likewise, when $dial_P = \{-1,0,2,2\}$ and when \textit {p < q}, then the value of \textit {p} on first occurrence of steady state value in $zone_0$ (either $od_2$ or $od_4$) is given as:
\begin{subnumcases}{p=}
 \left(\dfrac{\Delta - 4}{4}\right)^2 & \text{for $\Delta =  4k, k = 2m, \ k,m \in \mathbb{N}; p = 2j+1, \ j \in \mathbb{N}$} \label{eqn:ssv4a} \
  \\ [8pt]
  \left(\dfrac{\Delta - 4}{4}\right) \left(\dfrac{\Delta - 4}{4} + 1\right) - \left(\dfrac{\Delta - 4}{4} - 1\right)  & \text{for $\Delta =  4k, k = 2m+1, \ k,m \in \mathbb{N}; p = 2j+1, \ j \in \mathbb{N}$} \label{eqn:ssv4b} \
  \\ [8pt]
 \left(2 \times \left(\left(\dfrac{\Delta - 6}{4} + 1\right)^2 - 1\right)\right) + 3 & \text{for $\Delta =  4k+2, \ k \in \mathbb{N};\ p = 2j+1, \ j \in \mathbb{N}$} \label{eqn:ssv4c} \
  \\ [8pt]
   \left(2 \times \left(\dfrac{\Delta - 4}{4}\right) \left(\dfrac{\Delta - 4}{4} + 1\right)\right)+ 2 & \text{for $\Delta =  4k, \ k \in \mathbb{N};\ p = 2j, \ j \in \mathbb{N}$} \label{eqn:ssv4d} \
  \\ [8pt]
  \left(\dfrac{\Delta - 6}{4}\right) \left(\dfrac{\Delta - 6}{4} + 1\right) + 2& \text{for $\Delta =  4k+2,\ k \in \mathbb{N};\ p = 2j, \ j \in \mathbb{N}$} \label{eqn:ssv4e} \
\end{subnumcases}
A new system and method is presented to factorize any composite $n$, which is a product of two large primes $p$ and $q$, by searching for $\Delta$, which is the difference between the two primes, i.e. $\Delta = |p - q|$. 
\\
\\
Furthermore, we have been able to generalize the steady state value in different observation decks. In this paper we will cover these generalizations for $od_4$ and $od_5$ in any $zone$ and with respective dial settings. 
\\ \\
Steady state from \textbf{\begin{math}od_4\end{math}} (for any $\Delta$) when $dial_1 = \{0,-1\}$ or $dial_1 = \{-1,0\}$ and $v_1=v_2$
\\
\begin{equation}
SteadyStateValue_{_{od_{_4}}} =  \dfrac{(\Delta)^2}{2} + \dfrac{(v_1)^2}{2} \text{ for } v_1=2k \text{ where } k=2j+1 \text{ for } j \in \mathbb{N}_0 \label{eqn:ssvgenod4}
\end{equation}
\\
Steady state from \textbf{\begin{math}od_5\end{math}} (for any $\Delta$) when $dial_1 = \{-2,-1\}$ or $dial_1 = \{-1,-2\}$ and  $v_1=v_2$
\\
\begin{equation}
SteadyStateValue_{_{od_{_5}}} =  \Delta^2 + \dfrac{3}{4}(v_1)^2 \text{ for } v_1=4k \text{ where } k=2j+1 \text{ for } j \in \mathbb{N}_0 \label{eqn:ssvgenod5}
\end{equation}
It is also hypothesized that with the same dial settings, the steady state value will appear in the same $\Delta$ sieve zone of a specific observation deck for any $\Delta$. However, the number of rows in these $\Delta$ sieve zones will be different for different $\Delta$ and will be referred as "$\Delta$ sieve coverage".
\\
\\
We believe this is a new territory and from our research we couldn't find any prior art where $\Delta$ has been the primary focus for factorization purposes. As a result, we had to introduce some new notations and terminologies. There are 47 Tables and 10 Figures in this paper which will explain these ideas and concepts with clear examples. We now present the examples below making references to the respective tables and figures.
\\
\\
Kindly note, the values of $\Delta$ are selected in a manner to accommodate $\Delta = 4k$, $\Delta = 4k + 2$ forms when $\Delta$ is even and $\Delta = 4k+1$, $\Delta = 4k + 3$ forms when $\Delta$ is odd. The concepts will remain valid for any large value of $\Delta$.
\\
\\
``Place Table 1 here.'' Below observations can then be made:
\begin{itemize}[noitemsep]
\item Steady State Value = 74 as per equation (\ref{eqn:ssv1a}) or as per equation (\ref{eqn:ssvgenod4})
\item $dial_2\ (v_1 =v_2) = 2$
\item $p=13$ as per equation (\ref{eqn:ssv2a})
\item \begin{math}zone_0\end{math} for $od_4 $ is defined from \begin{math}id=7 \rightarrow \infty\end{math} 
\end{itemize}
``Place Table 2 here.'' Below observations can then be made:
\begin{itemize}[noitemsep]
\item Steady State Value = 121 as per equation (\ref{eqn:ssv1b})
\item $p=21$ as per equation (\ref{eqn:ssv2b})
\item \begin{math}zone_0\end{math} for $od_2$ is defined from \begin{math}id=11 \rightarrow \infty\end{math} 
\end{itemize}
``Place Table 3 here.'' Below observations can then be made:
\begin{itemize}[noitemsep]
\item Steady State Value = 36 as per equation (\ref{eqn:ssv1b})
\item $p=4$ as per equation (\ref{eqn:ssv2c})
\item \begin{math}zone_0\end{math} for $od_2$ is defined from \begin{math}id=2 \rightarrow \infty\end{math} 
\end{itemize}
``Place Table 4 here.'' Below observations can then be made:
\begin{itemize}[noitemsep]
\item Steady State Value = 64 as per equation (\ref{eqn:ssv1b})
\item $p=10$ as per equation (\ref{eqn:ssv2d})
\item \begin{math}zone_0\end{math} for $od_2$ is defined from \begin{math}id=5 \rightarrow \infty\end{math} 
\end{itemize}
``Place Table 5 here.'' Below observations can then be made:
\begin{itemize}[noitemsep]
\item Steady State Value = 244 as per equation (\ref{eqn:ssv1a}) or as per equation (\ref{eqn:ssvgenod4})
\item $dial_2\ (v_1 =v_2) = 2$
\item $p=50$ as per equation (\ref{eqn:ssv2e})
\item \begin{math}zone_0\end{math} for $od_4 $ is defined from \begin{math}id=25 \rightarrow \infty\end{math} 
\end{itemize}
``Place Table 6 here.'' Below observations can then be made:
\begin{itemize}[noitemsep]
\item Steady State Value = 144 as per equation (\ref{eqn:ssv1b})
\item $p=25$ as per equation (\ref{eqn:ssv4a})
\item \begin{math}zone_0\end{math} for $od_2$ is defined from \begin{math}id=13 \rightarrow \infty\end{math} 
\end{itemize}
``Place Table 7 here.'' Below observations can then be made:
\begin{itemize}[noitemsep]
\item Steady State Value = 36 as per equation (\ref{eqn:ssv1b})
\item $p=5$ as per equation (\ref{eqn:ssv4b})
\item \begin{math}zone_0\end{math} for $od_2$ is defined from \begin{math}id=3 \rightarrow \infty\end{math} 
\end{itemize} 
``Place Table 8 here.'' Below observations can then be made:
\begin{itemize}[noitemsep]
\item Steady State Value = 244 as per equation (\ref{eqn:ssv1a}) or as per equation (\ref{eqn:ssvgenod4})
\item $dial_2\ (v_1=v_2) = 2$
\item $p=51$ as per equation (\ref{eqn:ssv4c})
\item \begin{math}zone_0\end{math} for $od_4 $ is defined from \begin{math}id=26 \rightarrow \infty\end{math} 
\end{itemize} 
``Place Table 9 here.'' Below observations can then be made:
\begin{itemize}[noitemsep]
\item Steady State Value = 74 as per equation (\ref{eqn:ssv1a}) or as per equation (\ref{eqn:ssvgenod4})
\item $dial_2\ (v_1 =v_2) = 2$
\item $p=14$ as per equation (\ref{eqn:ssv4d})
\item \begin{math}zone_0\end{math} for $od_4 $ is defined from \begin{math}id=7 \rightarrow \infty\end{math} 
\end{itemize}
``Place Table 10 here.'' Below observations can then be made:
\begin{itemize}[noitemsep]
\item Steady State Value = 121 as per equation (\ref{eqn:ssv1b})
\item $p=22$ as per equation (\ref{eqn:ssv4e})
\item \begin{math}zone_0\end{math} for $od_2 $ is defined from \begin{math}id=11 \rightarrow \infty\end{math} 
\end{itemize}\
\\
In addition, we will relook at Table \ref{table:A1}, Table \ref{table:A5}, Table \ref{table:B3} and Table \ref{table:B4}, with focus turned to the difference between the consecutive observation deck values, represented as $df_x$, where $x$ is the reference number that identifies the respective observation deck. The hypothesis that a $\Delta$ sieve zone will be yielded when the sum total of these $df_x$ is 0 and mathematically represented as $\sum\limits_{i=x}^{k} {a}_i {df}_i = 0$ for $x, k \in \mathbb{N}$, $a \in \mathbb{Z}$ will also become clear. When this is achieved, the steady state value will be yielded in the respective $\Delta$ sieve zone and its value will be expressed as $\sum\limits_{i=x}^{k} {a}_i {od}_i$ for $x, k \in \mathbb{N}$, $a \in \mathbb{Z}$. One exception to this rule occurs at $switchover\ zones$, where $\Delta$ sieve zone is yielded but $\sum\limits_{i=x}^{k} {a}_i {df}_i \neq 0$
\\ \\
``Place Table 11 here.'' Below observations can then be made:
\begin{itemize}[noitemsep]
\item $\Delta$ sieve zone is yielded in $od_4$ when $df_1 + df_2 = 0$
\item There is an exception to the above rule, at the switchover zones (id=7), \\ where $\Delta$ sieve zone is yielded but $df_1 + df_2 \neq 0$
\end{itemize} 
``Place Table 12 here.'' Below observations can then be made:
\begin{itemize}[noitemsep]
\item $\Delta$ sieve zone is yielded in $od_4$ when $df_1 + df_2 = 0$
\item There is an exception to the above rule, at the switchover zones (id=25), \\ where $\Delta$ sieve zone is yielded but $df_1 + df_2 \neq 0$
\end{itemize} 
``Place Table 13 here.'' Below observations can then be made:
\begin{itemize}[noitemsep]
\item $\Delta$ sieve zone is yielded in $od_4$ when $df_1 + df_2 = 0$
\item There is an exception to the above rule, at the switchover zones (id=26), \\ where $\Delta$ sieve zone is yielded but $df_1 + df_2 \neq 0$
\end{itemize}
``Place Table 14 here.'' Below observations can then be made:
\begin{itemize}[noitemsep]
\item $\Delta$ sieve zone is yielded in $od_4$ when $df_1 + df_2 = 0$
\item There is an exception to the above rule, at the switchover zones (id=7), \\ where $\Delta$ sieve zone is yielded but $df_1 + df_2 \neq 0$
\end{itemize}\
We will continue to expand on this idea to yield $\Delta$ sieve zones in different observation decks. We will see this with few more examples for $od_5$ 
``Place Table 15 here.'' Below observations can then be made:
\begin{itemize}[noitemsep]
\item For $od_5$, $\Delta$ sieve zone will be yielded when $df_1 + df_2 + df_3 + df_4 = 0$ 
\item There is an exception to the above rule, at the switchover zones (id=56), where $\Delta$ sieve zone is yielded but $df_1 + df_2 + df_3 + df_4 \neq 0$
\item $dial_2\ (v_1=v_2) = 4$
\item $zone_0$ is defined from id = 56 to $\infty$
\item $SteadyStateValue_{_{od_{_5}}}=2128$ as per equation (\ref{eqn:ssvgenod5})
\end{itemize}
``Place Table 16 here.'' Below observations can then be made:
\begin{itemize}[noitemsep]
\item For $od_5$, $\Delta$ sieve zone will be yielded when $df_1 + df_2 + df_3 + df_4 = 0$ 
\item There is an exception to the above rule, at the switchover zones (id=61), where $\Delta$ sieve zone is yielded but $df_1 + df_2 + df_3 + df_4 \neq 0$
\item $dial_2\ (v_1=v_2) = 4$
\item $zone_0$ is defined from id = 61 to $\infty$
\item $SteadyStateValue_{_{od_{_5}}}=2316$ as per equation (\ref{eqn:ssvgenod5})
\end{itemize}
It is hypothesized that once the steady state is achieved in $zone_0$ of $od_2$, $od_4$ or $od_5$ for any $\Delta$, this steady state will continue until $\infty$
\\
\\
A quick note with regards to the primary purpose of covering the natural numbers of the form $p,q=2j,$ for $ j \in \mathbb{N}$ (the even numbers) as given in Table \ref{table:A3}, Table \ref{table:A4}, Table \ref{table:A5}, Table \ref{table:B4} and Table \ref{table:B5} is to highlight the applicability of relationship with $\Delta$ in different observation decks (od), which is especially useful when we try to factorize the composite $n$ by taking its multiple, like $4kn$, for $k \in \mathbb{N}$, as this multiple ($k$) will change the $\Delta$
\\
\\
In fact, relationship with $\Delta$ exists for any $kn$, for $k \in \mathbb{N}$. We will cover the "$odd\ \Delta$" use case, i.e. $\Delta = 4k+1 \ and\ \Delta = 4k+3$ for $k \in \mathbb{N}_0$ towards the end of this paper.
\subsection{$dial_P$}\
A dial pair, represented as "$dial_P = \{dial_1, dial_2\} = \{a_1, a_2, v_1, v_2\}$" (or like a matrix as given under the "Notation and setup" section ) plays a critical role in understanding how to sieve the $\Delta$. We have briefly explained how the dial pairs are chosen as per equations: (\ref{eqn:nas1a}), (\ref{eqn:nas1b}), (\ref{eqn:nas2a}) and (\ref{eqn:nas2b}), but we will take a quick moment to explain it a bit more in this subsection.
\\
\\
$dial_1$: Relative distance from $\lfloor \sqrt{n} \rfloor$
\begin{enumerate}
\item We start with the composite $n$ we want to factorize and it's $\sqrt{n}$.  We take the floor value of $\sqrt{n}$ ($\lfloor \sqrt{n} \rfloor$). So, we have $n$ and $\lfloor \sqrt{n} \rfloor$
\item Now, $\lfloor \sqrt{n} \rfloor$ can either be even ($\lfloor \sqrt{n} \rfloor = 2j$ for $j \in \mathbb{N}$) or odd ($\lfloor \sqrt{n} \rfloor = 2j + 1$ for $j \in \mathbb{N}$)
\item If $\lfloor \sqrt{n} \rfloor = 2j$ for $j \in \mathbb{N}$, we select a value $a_1$, likewise, if $\lfloor \sqrt{n} \rfloor = 2j + 1$ for $j \in \mathbb{N}$, we select a value $a_2$. Both $a_1, a_2 \in \mathbb{Z}_0$ and together $a_1, a_2$ are referred as $dial_1$ and represented as $dial_1 = \{a_1, a_2\}$
\item Let's use $dial_1 = \{-1, 1\}$ in below examples:
\begin{enumerate}[noitemsep]
\item If $n_1 = 137$,  $\lfloor \sqrt{n} \rfloor = 11$, then $a_2 = 1$ and $dial_1 = 1$
\item Likewise, if $n_2 = 147$, $\lfloor \sqrt{n} \rfloor = 12$, then $a_1 = -1$ and $dial_1 = -1$
\end{enumerate}
\end{enumerate}
Before we go to $dial_2$, we should see how $dial_1$ is put to use when calculating $d_1$:
\begin{enumerate}
\item $d_1 = \lfloor \sqrt{n} \rfloor + dial_1$
\item From example in point 4(a) above, $d_1 = 11 + 1 = 12$
\item From example in point 4(b) above, $d_1 = 12 + (-1) = 11$
\end{enumerate}
$dial_2$: Relative distance from $d_1$
\begin{enumerate}
\item If $d_1 = 2j$ for $j \in \mathbb{N}$, we select a value $v_1$, likewise, if $d_1 = 2j + 1$ for $j \in \mathbb{N}$, we select a value $v_2$. Both $v_1, v_2 \in \mathbb{Z}$ and together $v_1, v_2$ are referred as $dial_2$ and represented as $dial_2 = \{v_1, v_2\}$
\item Let's use $dial_2 = \{2, 4\}$ in below examples:
\begin{enumerate}[noitemsep]
\item If $n_1 = 137$,  $\lfloor \sqrt{n} \rfloor + dial_1 = 12$, then $v_1 = 2$ and $dial_2 = 2$
\item Likewise, if $n_2 = 147$, $\lfloor \sqrt{n} \rfloor + dial_1 = 11$, then $v_2 = 4$ and $dial_2 = 4$
\end{enumerate}
\end{enumerate}
$d_2$: 
\begin{enumerate}
\item $d_2 = d_1 + dial_2$
\item From above examples, for $n_1 = 137$, $d_1 = 12$, $dial_2 = 2$, then, $d_2 = 12 + 2 = 14$
\item Likewise, for $n_2 = 147$, $d_1 = 11$, $dial_2 = 4$, then, $d_2 = 11 + 4 = 15$
\end{enumerate}
There can be infinite dial pairs and changing their settings can change the $\Delta$ sieve zone within an observation deck. In some cases, $\{v_1, v_2\}$ of $dial_2$ plays a critical role (together with $\Delta$) in calculating the steady state value.
\\
\\
Also, when there are two or more dial pairs, like, $dial_{P1} = \{0, -1, 6, 6\}$ and  $dial_{P2} = \{-2, 1, 16, 16\}$, they may be connected with some relationship(s) with each other, such that, systematic increment to $dial_{P1}$ will require respective changes to $dial_{P2}$ in order to maximize the $\Delta$ sieve coverage in different or same observation decks. 
\\ 
\\
For e.g.:
\begin{itemize}
\item  $v_1{_{_{{dial_{_{P2}}}}}} = (2 \times v_1{_{_{{dial_{_{P1}}}}}}) + 4$; or
 \item $v_1{_{_{{dial_{_{P2}}}}}} = v_1{_{_{{dial_{_{P1}}}}}} + 10$
\end{itemize}
Assuming $v_1 = v_2$ in both dial pairs, if we now increment $dial_{P1} = \{0, -1, 10, 10\}$, there can be two corresponding updates that can be made to $v_1, v_2$ of $dial_{P2}$ as per the above relationships:
\begin{itemize}
\item  $dial_{P2} = \{-2, 1, 24, 24\}$; or
 \item $dial_{P2} = \{-2, 1, 20, 20\}$
\end{itemize}
In addition, intra dial relationship(s) between $a_1$, $a_2$, $v_1$ and $v_2$ will also exist and the increment applied to $a_1$, $a_2$, $v_1$ and $v_2$ to change the $\Delta$ sieve zone will be dependent on $\Delta$, in terms of whether $\Delta$ is of $4k, 4k+1, 4k+2$ or $4k + 3$ form, for $k \in \mathbb{N}$
\subsection{$\Delta$ sieve zones ($zone_0, zone_1,..., zone_x$)}\
This is the zone within an observation deck where $\Delta$ can be sieved. For fixed dial settings, the steady state value remains constant in these zones. As the dial settings change, the $\Delta$ sieve zone will shift and the steady state value will also change. However, this steady state value will continue to remain a function of $\Delta$, dial settings and some constant $k$, as mentioned in equation (\ref{eqn:ssvgen_main_1a})
\\
\\
The property for $\Delta$ sieve zone to shift (from $zone_0 --> zone_1 --> ... --> zone_x$) as dial settings change is also quite interesting. This is like tuning the $zone_0$ of $\Delta$ sieve zone for a particular observation deck first and then sieving the $\Delta$ in any other zones for any $n$ in that $\Delta$ series by appropriately changing the involved dials.
\\
\\
Example data for equations (\ref{eqn:ssvgenod4}) and (\ref{eqn:ssvgenod5}) is given below in Table \ref{table:4k2_od4_zone1} and Table \ref{table:od5_changing_dials_1} respectively.
\\
\\
\textbf{Note:} Since the screen real estate can't accommodate all the Table columns, we will only keep the most relevant columns to clarify the concepts.
\\ \\
``Place Table 17 here.'' Following observations can then be made. The $\Delta$ sieve zone in $od_4$ shifts up and moves to $zone_1$ as we change $dial_2$ to $6$
\begin{itemize}[noitemsep]
\item $\Delta$ sieve zone is yielded when $df_1 + df_2 = 0$ with an exception at $id=6$, referred as $switchover\ zone$ (details in subsection 4.4 below)
\item $zone_1$ for $od_4$ is defined from $6 <= id <26$ 
\item $dial_2\ (v_1=v_2) = 6$
\item $SteadyStateValue_{_{od_{_4}}} = 260$ as per equation (\ref{eqn:ssvgenod4})
\end{itemize}
``Place Table 18 here.'' Following observations can then be made. The $\Delta$ sieve zone in $od_5$ shifts up and moves to $zone_1$ as we change $dial_2$ to $12$
\begin{itemize}[noitemsep]
\item $\Delta$ sieve zone is yielded when $df_1 + df_2 + df_3 + df_4 = 0$ with an exception at $id=8$, referred as $switchover\ zone$
\item $zone_1$ for $od_5$ is defined from $8 <= id < 12$ 
\item $dial_2\ (v_1=v_2) = 12$
\item $SteadyStateValue_{_{od_{_5}}} = 2224$ as per equation (\ref{eqn:ssvgenod5})
\end{itemize}
Likewise, one can make the $\Delta$ sieve zones shift into different zones by changing the dials for different observation decks.
\subsection{$\Delta$ sieve coverage}\
$\Delta$ sieve coverage determines the sum total of numbers which are under the $\Delta$ sieve zone for a given set of dial pairs. As the dial settings change, the $\Delta$ sieve coverage  will also change. In the below examples, one can sieve for $\Delta$ in $od_1, od_2, od_4, od_5$ for given dial pair.
\begin{itemize}
\item ``Place Table 19 here.'' With $\Delta=160$ and $\{a_1, a_2, v_1, v_2\}$ = \{-2, 2, 12, 12\}, we can sieve the $\Delta$ in respective observation decks
\item ``Place Table 20 here.'' With $\Delta=480$ and $\{a_1, a_2, v_1, v_2\}$ = \{-2, 2, 12, 12\}, we can sieve the $\Delta$ in respective observation decks
\end{itemize}
Note the increase in $\Delta$ sieve coverage with increasing $\Delta$ in respective observation decks.
\\ \\ 
Numbers outside of $\Delta$ sieve coverage are referred as "$zoneless$". However, the same numbers may come under respective observation decks $\Delta$ sieve zones as dial settings change and will become "$zoners$".
\\
\\
With one dial pair, we are able to create six observation decks, but given $od_3, od_4, od_5, od_6$ are derived from $od_1$ and $od_2$, we can potentially have infinite such observatories from just one dial pair. It's also important to have non-overlapping $\Delta$ sieve zones distributed amongst these observation decks to increase the overall $\Delta$ sieve coverage.
\\
\\
Let's introduce another dial pair and visualize the increase in non-overlapping $\Delta$ sieve coverage with increasing $\Delta$ sieve zones.
\\ \\
$\Delta = 94$ \\ \\
$dial_{P1} = \{dial_1,dial_2\} = \{ 0,-1,6,6\}$ \\ \\
$dial_{P2} = \{dial_3,dial_4\} = \{ -2,1,16,16\}$ \\ \\
$d_3 =  \lfloor \sqrt{n} \rfloor + dial_3$ \\ \\
$d_4 =  d_3 + dial_4$ \\ \\
$od_7 =  d_3^2 - n$ \\ \\
$od_8 =  d_4^2 - n$ \\ \\
$od_9 = od_2 + od_8$ \\ \\
$od_{10} = od_4 + od_8$ \\ \\
$od_{11} = od_2 + (2 \times od_8) + od_4$ \\ \\
With $dial_{P1} = \{ 0,-1,6,6\}$ and $dial_{P2} = \{ -2,1,16,16\}$, the below steady state values can be observed from respective observation decks for \textbf{any} $\Delta$
\begin{flalign}
&SteadyStateValue_{_{od_{_{9a}}}} =\dfrac{\Delta ^2}{2}  + 72 & \label{eqn:ssvgenod9a}
\end{flalign}
\begin{flalign}
&SteadyStateValue_{_{od_{_{9b}}}}  = \dfrac{\Delta ^2}{2}  + 32 & \label{eqn:ssvgenod9b}
\end{flalign}
\begin{flalign}
&SteadyStateValue_{_{od_{_{10}}}}  = \dfrac{\Delta ^2}{2} + \left(\dfrac{\Delta}{2}\right)^2  + 168 & \label{eqn:ssvgenod10}
\end{flalign}
\begin{flalign}
&SteadyStateValue_{_{od_{_{11}}}}  = \Delta ^2 + \left(\dfrac{\Delta}{2}\right)^2  + 144 & \label{eqn:ssvgenod11}
\end{flalign}
``Place Table 21 here.'' One can see how $\Delta$ is now being sieved in different observation decks. These $\Delta$ sieve zones will remain constant for any $\Delta$ in respective observation decks. We will now be able to see how the $\Delta$ sieve coverage increases with increase in $\Delta$ through Figure 1. ``Place Figure 1 here.'' Below observations can then be made.
\\
\\
Here, we have excluded $od_8$ and $od_7$ as they are the perfect square forms and factorization can be easily achieved with difference of squares method. While the relationship with $\Delta$ still exists, we wanted to get the statistics for non perfect square forms. The \textbf{only} purpose of this graph in Figure 1 is to  confirm the hypothesis stated previously around tuning some $\Delta$ sieve zones in shape of fixed observation decks for any \textbf{one} $\Delta$ and sieving \textbf{any} other $\Delta$ through them. It is this property that makes it quite lucrative to approach the factorization problem via this method. In addition, this also confirmed that the $\Delta$ sieve coverage will also grow with increase in $\Delta$. 
\\
\\
However, one important point to note is, $p$ will also grow at the same rate for its corresponding $\Delta$, bringing $p$ \& $q$ close to each other with reference to their respective $\Delta$, making factorization simple through Fermat's method as well.
\\
\\
``Place Figure 2 here.''
\begin{enumerate}
\item $p$ on x-axis represents the growth of $p$ with increasing $\Delta$ sieve coverage (filled in red) for fixed $\Delta = 1002$
\item Area highlighted in red will come under $\Delta$ sieve coverage from $od_4$ and $od_5$ with base dial starting at $\{0, -1, 8, 8\}$ and $v_1=v_2$ having an incremental rate of +8 with each iteration
\item Light red signifies the non-overlapping coverage zones between $od_4$ and $od_5$, while dark red means overlapping zones between these two observation decks
\end{enumerate}
\subsection{Switchover zone ($soz$)}\
The first row of a zone when the steady state value is first observed is defined as the switchover zone ($soz$) for that observation deck. Based on $dial_1$, $dial_2$, ... , $dial_x$, the steady state values can be different at these switchover zones. Also, while steady state value is yielded at the switchover zones, however, $\sum\limits_{i=x}^{k} {df}_i \neq 0$ here.
\\
\\
``Place Table 22 here.'' Below observations can then be made:
\begin{itemize}[noitemsep]
\item The row in bold is $soz_1$ for $od_4$ with dial settings mentioned in Table 22
\end{itemize}\
``Place Table 23 here.'' Below observations can then be made:
\begin{itemize}[noitemsep]
\item The row in bold is $soz_0$ for $od_5$ with dial settings mentioned in Table 23
\end{itemize} 
\subsection{Switchover point ($sop$)}\
Switchover points are observed from $od_3$ and are quite dense at the beginning for a given $\Delta$ starting with $p=1$, but their frequency decreases as one progresses with the $\Delta$ series with increasing  $p$. ``Place Table 24 here.'' Below observations can then be made:
\begin{itemize}[noitemsep]
\item The rows in bold are switchover points for $od_3$ with dial settings mentioned in the Table 24.
\begin{enumerate}
\item $od_3{_{_{sop{_{_2}}}}} = 120\ (id=2, df_3 = -48)$
\item $od_3{_{_{sop{_{_3}}}}} = 168\ (id=3, df_3 = -48)$
\item $od_3{_{_{sop{_{_6}}}}} = 264\ (id=6, df_3 = -48)$
\end{enumerate}
\end{itemize}
\section{Intra $\Delta$ relationships (a$\Delta$)}\
The purpose is to find the "previous composite n" or "next composite n" with the same $\Delta$ as that of the composite $n$ we are trying to factorize. "previous composite n" will be represented as odd subscripts ($n_1, n_3, n_5, etc$) with higher odd subscripts away from the $n$ we are trying to factorize. Likewise, "next composite n" will be represented as even subscripts ($n_2, n_4, n_6, etc$) with higher even subscripts away from the $n$ we are trying to factorize.
\subsection{od connect} \
Different observation decks (od) can be connected to each other through a nontrivial relationship that can help to yield the previous or next composite n. We will see this with following example. ``Place Table 25 here.'' We need to factorize $n=219781$
\\
\\
As a first step, observation deck values will be calculated for this $n$, i.e. $od_1, od_2, od_3, od_4$ (highlighted in bold in Table 25). Based on this available information, some relationship needs to be established with other observation decks so we can determine $n{_{id=137}}$ (value of $n$ for $id=137$). In this example we can see:
\begin{equation}
(od{_3{_{_{id=136}}}}) - (od{_3{_{_{id=137}}}}) = -8 \label{eqn:intraeg1}\
\end{equation}
\begin{equation}
od{_4{_{_{id=136}}}} - df{_4{_{_{id=137}}}} = od{_4{_{_{id=137}}}} \label{eqn:intraeg2}\
\end{equation}
\begin{equation}
(\dfrac{od{_3{_{_{id=136}}}}}{4} - df{_4{_{_{id=137}}}}) \times 2 = od{_4{_{_{id=137}}}} \label{eqn:intraeg3}
\end{equation}
In above equations (\ref{eqn:intraeg2}) and (\ref{eqn:intraeg3}), "$df{_4{_{_{id=137}}}}$", represents the value of $df_4$ for $id=137$. 
\\
\\
The relationship in equation (\ref{eqn:intraeg3}) is key in finding the next $n$ . However, It's not known how many such nontrivial relationships are out there to make this method effective.
\\
\\
Let's find $n_2$ for $id=137$. From (\ref{eqn:intraeg1}), (\ref{eqn:intraeg2}) and (\ref{eqn:intraeg3}):
\begin{equation}
od_3 = 1884, od_4 = -214, df_4 = 576 \label{eqn:intraeg4}\
\end{equation}
\begin{equation} \label{eqn:intraeg5}\
\begin{split}
od_3 = od_2 - od_1 \\ 
od_4 = od_2 + od_1 
\end{split}
\end{equation}
From (\ref{eqn:intraeg4}) and (\ref{eqn:intraeg5}) and assuming we are not going to hit the $switchover\ point\ (sop)$, $d{_1{_{_{id=137}}}}$ = $d{_1{_{_{id=136}}}} + 2$, so $d{_1{_{_{id=137}}}}$ = 470
\begin{equation} \label{eqn:intrasolveforn}\
\begin{split}
od_1 & = -1049, od_2 = 835 \\
od_1 & = d_1^2 - n \\
-1049 & = 470^2 - n \\
n_2 & = 221949
\\
\end{split}
\end{equation}
Once $n_2$ is known, determining its factors should be simple, assuming $n_2$ is not a difficult semiprime. Once the factors of $n_2$ are determined, it will be easy to find out the nontrivial factors of the composite $n$ we are trying to factorize.
\subsection{Range of previous/next $n$} \
The aim is to first find the range within which the previous or next composite $n$ values will fall, such that they have the same $\Delta$ as that of the composite $n$ we are trying to factorize. Once the range is known, the second goal is to find some nontrivial relationship between min/max value of this range and $od_x, df_x, \Delta$, so the previous or next composite $n$ can be found efficiently. Kindly note, we are not covering this second goal in this paper. ``Place Table 26 here.''
\\
\\
If $n = 1643$, $od_1 = -122$, then:
\\
\phantom{x}\hspace{3ex}$n_3 = 1323$, $d_{1_{n_3}} = 35$, $n_1 =  1479$, $d_{1_{n_1}} = 37$
\\
\phantom{x}\hspace{3ex}$n_2 = 1815$, $d_{1_{n_2}}  = 41$, $n_4 = 1995$, $d_{1_{n_4}}  = 43$
\\
\\
Assuming the previous and next $n$ values are not going to cross the switchover points, which means $d_1$ will be linear and we can find the respective ranges for $n_1$,  $n_2$, ... ,$n_x$
\\
\\
$n_1range=$ $(d_{1_{n_3}}^2 - od_1)$ to $(d_{1_{n_1}}^2 - od_1)$ $= (35^2 + 122)$ to $(37^2 + 122)$ $= 1347$ to $1491$
\\
\phantom{x}\hspace{1ex} $n_1range_{min} = 1347$ and $n_1range_{max} = 1491$
\\
\\
$n_2range =$ $(d_{1_{n_2}}^2 - od_1)$ to $(d_{1_{n_4}}^2 - od_1)$ $= (41^2 + 122)$ to $(43^2 + 122)$ $= 1803$ to $1971$
\\
\phantom{x}\hspace{1ex} $n_2range_{min} = 1803$ and $n_2range_{max}  = 1971$
\\
\\
Summary of the ranges we obtained above:
\begin{itemize}[noitemsep]
\item $n_1$ range: 1347 - 1491
\item $n_2$ range: 1803 - 1971
\end{itemize}
\subsection{Reflection over \{X,Y\}}\
As we are looking for $\Delta$ indirectly in this observation deck by finding the previous or next composite integer having the same $\Delta$ as that of the composite $n$ we are trying to factorize, this goal is achieved by looking for reflected $od_6$ values here.
\\
\\
``Place Table 27 here.'' $od_6 = \sqrt{4n + (od_1\times od_2)}$ for  $\{v_1, v_2\}$ = $\{2, 2\}$. Reflection over \color{red}\textbf{4,4} \color{black} is observed. Since $|X - Y| = 0$, values will be same over the reflection points for some interval. Here \{X,Y\} = \{4,4\} and is represented as: $ro\{4,4\}_{0}$
\\
\\
``Place Table 28 here.'' Reflection over \color{red}\textbf{3,5} \color{black} is observed. Since $|X - Y| = 2$, reflected values will differ by 2 over the reflection points. Here \{X,Y\} = \{3,5\} and is represented as: $ro\{3,5\}_{2}$
\\
\\
\textbf{Note}: 
\begin{enumerate}[noitemsep]
\item Reflection can be some distance $d$ apart from the point(s) of convergence (X, Y)
\item When searching for $od_6$, it's not known at the start whether to search up or down and hence one may need to go in both directions simultaneously
\end{enumerate}
$od_6$ search is carried out based on certain assumptions around how the reflection will occur so as to find a nontrivial $n_1$ or $n_2$, such that the $\Delta$ between the factors of $n_1$ or $n_2$ is same as that of the composite $n$ we are trying to factorize. Below equation is used for this search:
\begin{equation} \label{eqn:od6_searchzone}\
\begin{split}
(od_6)^2 & = d_1^4 + 4(d_1)^3 + n_2^2  + 4(d_1)^2  - 2n_2(d_1)^2 - 4n_2d_1  \\
OR
\\
(od_6)^2 & = d_1^4 + 4(d_1)^3 + n_1^2  + 4(d_1)^2  - 2n_1(d_1)^2 - 4n_1d_1  \\
\end{split}
\end{equation}
\\ Key points with regards to equation (\ref{eqn:od6_searchzone}):
\begin{enumerate}[noitemsep]
\item $n_2 > n > n_1$
\item $d_1$ is dependent on the corresponding $n\ (n_1, n_2, etc)$ and on fixed dial settings
\item There can be many trivial $n_1$ or $n_2$ satisfying the above equation
\item There is also a possibility for some of the $\Delta$s of these trivial $n_1$ or $n_2$ to be in close vicinity of the nontrivial $\Delta$ of composite $n$ we are trying to factorize
\item Reflection is not a mandatory condition for carrying out $od_6$ search
\item The below graphs visualizes all $od_6$ values for different $\Delta$s until steady state is achieved
\end{enumerate}
``Place Figure 3 here.'' and ``Place Figure 4 here.'' Key points with regards to Figure 3 and Figure 4 below:
\begin{enumerate}[noitemsep]
\item The graphs are dependent on $dial_2$, i.e. $v_1 = v_2 = 2$, i.e. graphs will change if $v_1$ or $v_2$ changes
\\
\item The graphs looks like letter "A", increasing (or decreasing) in its size as it attains the steady state. For this reason, this will be referred as "$A_+$" or "$A_-$" graph.
\\
\item We observed the same $A_+$ pattern for large $\Delta$s  (for both $\Delta=4k$ and $\Delta=4k+2$ form, for $k \in \mathbb{N}$) and hence we hypothesized that reflections in $od_6$ are imminent
\\
\item Note: $n$ is on log scale on x-axis
\end{enumerate}
``Place Figure 5 here.'' and ``Place Figure 6 here.'' A quick look at some other graphs with different $v_1 = v_2$ values.
\begin{enumerate}[noitemsep]
\item Figure 5: $v_1 = v_2 = 14$
\item Figure 6: $v_1 = v_2 = 37$
\end{enumerate}
\section{Inter $\Delta$ relationships (r$\Delta$)}\
\begin{equation} \label{eqn:interdeltagen}
known\ \Delta = delta\_val_{_{id=id\_val}} \{dial_{P1},...dial_{Pn}\} \simeq unknown\ \{\Delta_1, \Delta_2, \Delta_3\} \simeq \{n1, n2, n3\}
\end{equation} 
This section is about finding relationships with ``known $\Delta$'' and ``unknown $\Delta$s''. Any $od_x$, $df_x$, $zone_x$ where we are searching for $\Delta$ can have nontrivial relationships with any other $od_x$, $df_x$, $zone_x$ of any one or more known $\Delta$s
\\
\\
Also, with a given known $\Delta$, known dial pairs (one or more) and some fixed $id$ (one or more), many other unknown $\Delta$s can be selected ($\Delta_1, \Delta_2, \Delta_3$) such that for these fixed dial pairs, known $\Delta$ yields a $\Delta$ sieve zone at selected $id$ value(s) in any one of the observation decks, and the same number of dial pairs will yield $\Delta$ sieve zones for the respective unknown $\Delta$s. Conversely, for many unknown $\Delta$s, a known $\Delta$ can be found with some fixed id (one or more), such that the same number of dial pairs will yield $\Delta$ sieve zones for both known $\Delta$ and unknown $\Delta$s.  
\\
\\
Let's see equation (\ref{eqn:interdeltagen} ) with below example:
\\
\begin{equation} \label{eqn:interdeltaeg} \
known\ \Delta =122_{_{id=93}} \{0, -1, 8, 8\} \simeq \{162, 178, 202\} \simeq \{785539, 936863, 441383\}
\end{equation} 
\\
Details of notation used in equations (\ref{eqn:interdeltagen}) and (\ref{eqn:interdeltaeg}):
\begin{itemize}
\item "known $\Delta=122_{_{id=93}} \{0, -1, 8, 8\}$": A $\Delta$ sieve zone is observed (in $od_2$) at $id=93$ for $\Delta=122$ with $dial_P=\{0, -1, 8, 8\}$
\item "$\simeq$": Similar/equal symbol is used for now. Having a $\Delta$ sieve congruence symbol will be useful
\item "\{162, 178, 202\}": Set of unknown $\Delta$s
\item "$\{785539, 936863, 441383\}$": A $\Delta$ sieve zone will be observed for these $n$ for $dial_P=\{0, -1, 8, 8\}$ in one of the observation decks.
\end{itemize}
``Place Table 29 here.''
$n1=785539$, $n2=936863$, $n3=441383$
\\
Here, known $\Delta\{122\}$  is linked with unknown $\Delta\{162, 178, 202\}$ with a fixed reference of $dial_P=\{0,-1,8,8\}$
\\
\\
\textbf{Note:} $dial_P$ of unknown $\Delta$s doesn't necessarily have to be fixed and can change as long as there is some relationship between $dial_P$ of unknown $\Delta$s with $dial_P$ of known $\Delta$ for factorization purposes.
\section{The sum series ($\sum = p + q$)}\
Very similar to $\Delta$ series, instead of studying all $p,q$ linked with the same $\Delta$, we study all $p,q$ linked with the same $\sum$, such that, $\Delta_{|p-q|} = \sum_{p+q}$
\\ \\
``Place Table 30 here.'' and ``Place Table 31 here.''  Below observations can then be made.
\\ \\
The $N$ on $\sum$ series, as described in Table \ref{table:sum_series_1} and Table \ref{table:sum_series_2} is connected with $od_6$ on the $\Delta$ series. It is this relationship between $N$ and $od_6$ that not only lays the foundation for a new factorization method but also a new trapdoor function supporting two different keys for encryption and decryption respectively. This is explained clearly with more examples in the next sections.
\section{The equilibrium of $\Delta_{|p-q|}$ and $\sum_{p+q}$}\
The equilibrium between $N$ on $\sum_{p+q}$ and $od_6$ on $\Delta_{|p-q|}$ series is reached when:
\begin{subnumcases}{N_{_{\Delta \rightleftharpoons \sum}}=}
  od_6{_{_{_{\Delta_{|p-q|}}}}} & \text{for $\Delta = 4k, p = 2j+1, \ k, j \in \mathbb{N}$; $\Delta = 4k+2, p = 2j, \ k, j \in \mathbb{N}$} \label{eqn:ssvod6_1a} \
  \\ [6pt]
  od_6{_{_{_{\Delta_{|p-q|}}}}} + 1  & \text{for $\Delta = 4k+2, p = 2j+1, \ k, j \in \mathbb{N}$; $\Delta = 4k, p = 2j, \ k, j \in \mathbb{N}$} \label{eqn:ssvod6_1b} \
\end{subnumcases}
``Place Table 32 here.'' and ``Place Table 33 here.'' The equilibrium is reached when $N_{_{\Delta \rightleftharpoons \sum}} = od_6{_{_{_{\Delta_{|p-q|}}}}} = 99$ as per equation (\ref{eqn:ssvod6_1a})
\\
\\
``Place Table 34 here.'' In this Table:
\begin{itemize}[noitemsep]
\item $N_{_{\sum_{p+q}}}$ is from Table \ref{table:sum_equilibrim_1}, starting from $N=99$, for $od_6=0$ ($id=6$) and moving down
\\
\item $od_6{_{_{_{\Delta_{|p-q|}}}}}$ is from Table \ref{table:delta_equilibrim_1}, starting from $od_6=99$ ($id=21$) and moving up
\\
\item $N_{_{\sum_{p+q}}}$ - $od_6{_{_{_{\Delta_{|p-q|}}}}}$ = $\Delta \rightleftharpoons \sum$\
\\
\item $df_{_{\Delta\sum}}$: Difference between consecutive $\Delta \rightleftharpoons \sum$\ values
\end{itemize}
``Place Table 35 here.'' From Table 34 and Table 35, we define below points:
\begin{enumerate}[noitemsep]
\item $od_6$ on $\Delta$ series is said to be in "full equilibrium" with $N$ on the $\sum$ series, iff :
\begin{itemize}
\item "$od_6\ +$ (any) $\Delta \rightleftharpoons \sum$" can get $od_6$ to jump to (any) $N$ on it's $\sum$ series, such that, sum of the factors of $N$ equals difference between the factors of $n$ on the $\Delta$ series
\end{itemize}
\item $od_6$ on $\Delta$ series is said to be in "close equilibrium" with $n$ on the $\sum$ series, iff :
\begin{itemize}
\item "$od_6\ +$ (any) $\Delta \rightleftharpoons \sum\ +$ $k,\ k \in \mathbb{N}$" can get $od_6$ to jump to (any) $n$ on it's $\sum$ series, such that, sum of the factors of $N$ equals difference between the factors of $n$ on the $\Delta$ series
\end{itemize}
\end{enumerate}
When we can make the composite $n$ on $\Delta$ series jump to another composite $N$ on $\sum$ series, and assuming the corresponding $N$ on $\sum$ series is not a difficult semiprime, and once the factors of $N$ are found, one can add these factors to get the $\Delta$. Once the $\Delta$ is known, the factors of the composite $n$ we are trying to factorize can be easily determined.
\\
\\
$\Delta \rightleftharpoons \sum$ in Table 34 and Table 35 contain the constants, which $od_6$ on any $\Delta$ series can use to jump to $N$ on the $\sum$ series. These constants will be referred as the "global equilibrium constants $(gec)$". When studying  $gec$ with consecutive $\Delta$s ($4k$ and $4k+2$ forms separately), we can observe $gec$ grows with increasing $\Delta$, however, the non-constant elements ($nce$), connecting $\Delta$ and $\sum$ series also grow at the same rate. Below two equations can summarize the relationship between $gec$ and $nce$.
\begin{equation}
Total\ number\ of\ elements\ on\ \Delta \ series = Number\ of\ gec + Number\ of\ nce
\end{equation}
\begin{equation}
Total\ number\ of\ elements\ on\ \Delta \ series = (2 \times Number\ of\ gec) + residue
\end{equation}
Number of $gec$ and $nce$ are calculated based on continuous comparison between consecutive $\Delta$ and $\Delta + 4$ (for $\Delta=4k, k \in \mathbb{N}$). ``Place Figure 7 here.'' In reference to this figure:
\begin{enumerate}[noitemsep]
\item Start with $\Delta=20$ and get all the elements on $\Delta$ series (from p=1 until steady state)
\item Get all the elements on $\sum=20$ series, such that number of elements between $\Delta$ and $\sum$ series is same
\item Calculate $\Delta \rightleftharpoons \sum$ and increment $\Delta$ by 4
\item Repeat steps 1 to 3 for this incremented $\Delta$ and $\sum$
\item Check how many elements are common ($gec$) and not common ($nce$) between \\ $\Delta \rightleftharpoons \sum = 20$ and $\Delta \rightleftharpoons \sum = 24$
\item Continuously increment $\Delta$ by 4 and keep comparing the $gec$ and $nce$ with the previous $\Delta$ until a preset threshold is reached
\end{enumerate}
``Place Table 36 here.'' and ``Place Table 37 here.'' Below observations can then be made.
\\
\\
From Table \ref{table:delta_4k2_1} and Table \ref{table:sum_4k2_1}, the equilibrium is reached when $N_{_{\Delta \rightleftharpoons \sum}} = od_6{_{_{_{\Delta_{|p-q|}}}}} + 1 = 121$ as per equation (\ref{eqn:ssvod6_1b})
\\
\\
``Place Table 38 here.'' In this Table:
\begin{itemize}[noitemsep]
\item $N_{_{\sum_{p+q}}}$ is from Table \ref{table:sum_4k2_1}, starting from $N=121$, for $od_1=0$ ($id=6$) and moving down
\\
\item $od_6{_{_{_{\Delta_{|p-q|}}}}}$ is from Table \ref{table:delta_4k2_1}, starting from $od_6=120$ ($id=26$) and moving up
\\
\item $N_{_{\sum_{p+q}}}$ - $od_6{_{_{_{\Delta_{|p-q|}}}}}$ = $\Delta \rightleftharpoons \sum$\
\\
\item $df_{_{\Delta\sum}}$: Difference between consecutive $\Delta \rightleftharpoons \sum$\ values
\end{itemize}
 ``Place Table 39 here.'' From Table 38 and Table 39, we define below points:
\begin{enumerate}[noitemsep]
\item $\Delta \rightleftharpoons \sum$ in Table 38 and Table 39 contain the global equilibrium constants ($gec$), which $od_6$ on any $\Delta$ series can use to jump to $N$ on the $\sum$ series. The points mentioned earlier for $\Delta = 4k$ form will remain same for $\Delta = 4k+2$ form as well
\end{enumerate}
``Place Figure 8 here.'' Both Figure 7 and Figure 8 are almost similar in terms of the growth of $gec$ and $nce$. However, there is a difference in growth of $residue$ between $\Delta = 4k$ and $\Delta = 4k  + 2$ for $k \in \mathbb{N}$ form. ``Place Figure 9 here.'' and ``Place Figure 10 here.'' These two figures highlight the difference between the residue graphs.
\\
\\
For completeness, below two examples will cover the use cases when $p = 2k$ for $k \in \mathbb{N}$ form, i.e. $p$ is an even number.  ``Place Table 40 here.'' and ``Place Table 41 here.''  We can observe From Table \ref{table:delta_equilibrim} and Table \ref{table:sum_equilibrim}, the equilibrium is reached when $N_{_{\Delta \rightleftharpoons \sum}} = od_6{_{_{_{\Delta_{|p-q|}}}}} + 1 = 100$ as per equation (\ref{eqn:ssvod6_1b})
\\
\\
``Place Table 42 here.'' and ``Place Table 43 here.'' We can observe From Table \ref{table:delta_4k2} and Table \ref{table:sum_4k2}, the equilibrium is reached when $N_{_{\Delta \rightleftharpoons \sum}} = od_6{_{_{_{\Delta_{|p-q|}}}}} = 120$ as per equation (\ref{eqn:ssvod6_1a})
\section{The Trapdoor}\
Computing $od_6$ from a given $n$ is easy, however, it's not easy to go back to $n$ from $od_6$ even though $\Delta$ is known. Also, we have seen above how $\Delta \rightleftharpoons \sum$ equilibrium connects $od_6$ on $\Delta$ series with $N$ on $\sum$ series and this allows us to find $n$ from $od_6$ easily when this relationship and $\Delta$ are given.
\\
\\
When $\Delta$ is known publicly
\begin{enumerate}
\item Algorithm
\begin{enumerate}[noitemsep]
\item Message encryption and generating the private key:
\begin{enumerate}[noitemsep]
\item Select some large $\Delta$ (such that integer $p < p_{ssv}$)
\item Select below dial pairs:
\begin{enumerate}
\item For $\Delta = 4k$ for $k \in \mathbb{N}$, $\Delta_{dial_{P1}} = \{0,-1,2,2\}$, $\sum_{dial_{P2}}=\{-1,0,2,2\}$
\item For $\Delta = 4k+2$ for $k \in \mathbb{N}$, $\Delta_{dial_{P1}} = \{-1,0,2,2\}$, $\sum_{dial_{P2}}=\{-1,0,2,2\}$
\end{enumerate}
\item Compute steady state value ($ssv$) from  equation (\ref{eqn:ssvgenod4}) 
\item Compute $p_{_{ssv}}$ at this steady state value from equation (\ref{eqn:ssv2a}). As dials are fixed, equation (\ref{eqn:ssv2a}) can be used for both $\Delta = 4k$ and $\Delta = 4k+2$ for $k \in \mathbb{N}$
\item Compute $od_{6_{_{ssv}}}$ at this steady state value
\item Convert message M into an integer p
\item Compute $p_{_{dist}}=p_{_{ssv}} - p$
\item Compute $q = p + \Delta$
\item Compute $n = p \times q$
\item Compute $od_1$ and $od_2$
\item Compute $od_6 = \sqrt{(4n) + (od_1\times od_2)}$ \\
\item $od_6$ is the ciphertext \\
\item Compute $N_{_{\sum}}$ :
\begin {enumerate}
\item From equation (\ref{eqn:ssvod6_1a}), $N_{_{\Delta \rightleftharpoons \sum}} = od_{6_{_{ssv}}}$
\item From equation (\ref{eqn:ssvod6_1b}), $N_{_{\Delta \rightleftharpoons \sum}} = od_{6_{_{ssv}}} + 1$
\end{enumerate}
\item From $N_{_{\Delta \rightleftharpoons \sum}}$ and $\Delta$, compute $p_{_{\Delta \rightleftharpoons \sum}},q_{_{\Delta \rightleftharpoons \sum}}$ on the sum series, such that \\ $p_{_{\Delta \rightleftharpoons \sum}} \times q_{_{\Delta \rightleftharpoons \sum}} = N_{_{\Delta \rightleftharpoons \sum}}$ and $p_{_{\Delta \rightleftharpoons \sum}} < q_{_{\Delta \rightleftharpoons \sum}}$ (Note: $\Delta = (p - q)_{_{\Delta series}} = (p + q)_{_{\sum series}}$) \\
\item Compute $p_{_{sum\_series}} = p_{_{\Delta \rightleftharpoons \sum}} - p_{_{dist}}$
\item Compute $q_{_{sum\_series}} = \Delta - p_{_{sum\_series}}$
\item Compute $N_{_{sum\_series}} = p_{_{sum\_series}} \times q_{_{sum\_series}}$ \\
\item Compute $(\Delta \rightleftharpoons \sum)_{_{private}} = N_{_{sum\_series}} - od_6$ \\
\item $(\Delta \rightleftharpoons \sum)_{_{private}}$  is the private key \\
\end{enumerate}
\item Decryption when cipher text: $od_6$, private key: $(\Delta \rightleftharpoons \sum)_{_{private}}$ and $\Delta$ are given:
\begin{enumerate}[noitemsep]
\item Compute $N_{_{sum\_series}} = od_6 + (\Delta \rightleftharpoons \sum)_{_{private}}$ \\
\item Compute 2 factors of $N_{_{sum\_series}}$: $p_{_{sum\_series}}$ and  $q_{_{sum\_series}}$, such that \\ $\Delta = p_{_{sum\_series}} + q_{_{sum\_series}}$ and $p_{_{sum\_series}}$ < $q_{_{sum\_series}}$ \\
\item Compute steady state value ($ssv$) from  equation (\ref{eqn:ssvgenod4}) \\
\item Compute $od_{6_{_{ssv}}}$ at this steady state value \\
\item Compute $p_{_{ssv}}$ at this steady state value from equation (\ref{eqn:ssv2a}). As dials are fixed, equation (\ref{eqn:ssv2a}) can be used for both $\Delta = 4k$ and $\Delta = 4k+2$ for $k \in \mathbb{N}$ \\
\item Compute $N_{_{\Delta \rightleftharpoons \sum}}$ :
\begin {enumerate}
\item From equation (\ref{eqn:ssvod6_1a}), $N_{_{\Delta \rightleftharpoons \sum}} = od_{6_{_{ssv}}}$
\item From equation (\ref{eqn:ssvod6_1b}), $N_{_{\Delta \rightleftharpoons \sum}} = od_{6_{_{ssv}}} + 1$ \\
\end{enumerate}
\item From $N_{_{\Delta \rightleftharpoons \sum}}$ and $\Delta$, compute $p_{_{\Delta \rightleftharpoons \sum}},q_{_{\Delta \rightleftharpoons \sum}}$ on the sum series, such that \\ $p_{_{\Delta \rightleftharpoons \sum}} \times q_{_{\Delta \rightleftharpoons \sum}} = N_{_{\Delta \rightleftharpoons \sum}}$ and $p_{_{\Delta \rightleftharpoons \sum}} < q_{_{\Delta \rightleftharpoons \sum}}$\\
\item Compute $q_{_{dist}} = q_{_{sum\_series}} - q_{_{\Delta \rightleftharpoons \sum}}$ \\
\item $p = p_{_{ssv}} - q_{_{dist}}$ \\
\item This $p$ is the integer we encrypted earlier and can be converted into original message M
\end{enumerate}
\end{enumerate}
\end{enumerate}
In the below example, we will see the steps involved to encrypt the message (without using any padding scheme) by selecting a $\Delta$, generate the private key and decrypt the encrypted message when $\Delta$ and private key are provided.
\begin{enumerate}
\item Example-1
\begin{enumerate}[noitemsep]
\item Message encryption and generating the private key:
\begin{enumerate}[noitemsep]
\item $\Delta = 137136$
\item $\Delta_{dial_{P1}} = \{0,-1,2,2\}$, $\sum_{dial_{P2}}=\{-1,0,2,2\}$
\item Steady State Value = 9403141250
\item $p_{_{ssv}} = 2350716745$
\item $od_{6_{_{ssv}}}$ = 4701570623
\item Convert message "AUM" (using ASCII code) into an integer p = \textbf{658577}
\item $p_{_{dist}}$ = 2350058168
\item $q$ = 795713
\item $n$ = 524038280401
\item $od_1$ = -1279185 and $od_2$ = 1616435
\item $od_6$ = \textbf{168623}
\item $od_6$ is the ciphertext
\item $N_{_{\Delta \rightleftharpoons \sum}}$ = 4701570623
\item $p_{_{\Delta \rightleftharpoons \sum}}$ = 68567 and $q_{_{\Delta \rightleftharpoons \sum}}$ = 68569
\item $p_{_{sum\_series}}$ = -2349989601
\item $q_{_{sum\_series}}$ = 2350126737
\item $N_{_{sum\_series}}$ = -5522773392982061937
\item $(\Delta \rightleftharpoons \sum)_{_{private}}$ = \textbf{-5522773392982230560}
\item $(\Delta \rightleftharpoons \sum)_{_{private}}$  is the private key \\
\end{enumerate}
\item Decryption when cipher text: $od_6$, private key: $(\Delta \rightleftharpoons \sum)_{_{private}}$ and $\Delta$ are given:
\begin{enumerate}[noitemsep]
\item $N_{_{sum\_series}}$ = -5522773392982061937
\item $p_{_{sum\_series}}$ = -2349989601 and  $q_{_{sum\_series}}$ =  2350126737
\item Steady state value  = 9403141250
\item $od_{6_{_{ssv}}}$ = 4701570623
\item $p_{_{ssv}}$ = 2350716745
\item $N_{_{\Delta \rightleftharpoons \sum}}$ = 4701570623
\item $p_{_{\Delta \rightleftharpoons \sum}}$ =  68567 and $q_{_{\Delta \rightleftharpoons \sum}}$ = 68569
\item $q_{_{dist}}$ = 2350058168
\item $p$ = \textbf{658577}
\item This $p$ is the integer we encrypted earlier. With available ASCII mapping, this can be converted back to the original message "AUM" 
\end{enumerate}
\end{enumerate}
\end{enumerate}
\section{Odd $\Delta$}\
\begin{subnumcases}{SteadyStateValue_{{od}_5}=}
  \Delta^2 + 3 & \text{for \begin{math}\Delta =  4k+3, p = 2j+1, \ k, j \in \mathbb{N}; dial_{P1} = \{-1,0,2,2\}\end{math}}  \label{eqn:odd_delta_1} \
  \\ [6pt]
  \Delta^2 + 3  & \text{for \begin{math}\Delta =  4k+1, p = 2j+1, \ k, j \in \mathbb{N}; dial_{P1} = \{0,-1,2,2\}\end{math}} \label{eqn:odd_delta_2} \  
 \\ [6pt]
  \Delta^2 + 3 & \text{for \begin{math}\Delta =  4k+1, p = 2j, \ k, j \in \mathbb{N}; dial_{P1} = \{-1,0,2,2\}\end{math}}  \label{eqn:odd_delta_3} \
  \\ [6pt]
  \Delta^2 + 3  & \text{for \begin{math}\Delta =  4k+3, p = 2j, \ k, j \in \mathbb{N}; dial_{P1} = \{0,-1,2,2\}\end{math}} \label{eqn:odd_delta_4} \  
\end{subnumcases}\
\\
The principle of tuning the $\Delta$ sieve zone remains the same for odd $\Delta$ as well, i.e. we are looking for $\sum {df}_x = 0$, where $x$ corresponds to those observation decks whose respective $df$ when summed gives 0. One observation deck can be used multiple times to achieve this zero sum state.
\\
\\
As we are observing the steady state from $od_5$, we are looking for zone(s) where $df_1 + df_2 + df_3 + df_4 = 0$, i.e. $\sum\limits_{x=1}^{4} {df}_x = 0$. Steady State is yielded at the switchover zone as well, however, $\sum\limits_{x=1}^{4} {df}_x \neq 0$ here.
\\
\\
``Place Table 44 here.'' Below observations can then be made:
\begin{itemize}[noitemsep]
\item $SteadyStateValue_{_{od_{_5}}}=532$ as per equation (\ref{eqn:odd_delta_1})
\item $df_1 + df_2 + df_3 + df_4 = 0$ is achieved at $id = 19$
\end{itemize}
``Place Table 45 here.'' Below observations can then be made:
\begin{itemize}[noitemsep]
\item $SteadyStateValue_{_{od_{_5}}}=628$ as per equation (\ref{eqn:odd_delta_2})
\item $df_1 + df_2 + df_3 + df_4 = 0$ is achieved at $id = 22$
\end{itemize}
``Place Table 46 here.'' Below observations can then be made:
\begin{itemize}[noitemsep]
\item $SteadyStateValue_{_{od_{_5}}}=628$ as per equation (\ref{eqn:odd_delta_3})
\item $df_1 + df_2 + df_3 + df_4 = 0$ is achieved at $id = 22$
\end{itemize}
``Place Table 47 here.'' Below observations can then be made:
\begin{itemize}[noitemsep]
\item $SteadyStateValue_{_{od_{_5}}}=532$ as per equation (\ref{eqn:odd_delta_4})
\item $df_1 + df_2 + df_3 + df_4 = 0$ is achieved at $id = 18$
\end{itemize}
\section{Final Remarks}\
This is our first number theory paper and even though we have consulted numerous online and offline resources and have reviewed this paper and re-reviewed the review countless times, we still might have had missed some important guidelines and we seek your forgiveness if a spelling or grammar or structure or some obvious rules of the paper writing terrain were not followed correctly and worst, if we still left that typo dangling out there. Our primary focus was to clearly bring out the dimension of $\Delta_{|p-q|}$, $\sum_{p+q}$ and the relationship between the two, i.e. $\Delta_{|p-q|} \rightleftharpoons \sum_{p+q}$. We believe, we have covered a lot of ground here but we are also cognizant that this might just be the beginning of a more beautiful journey some of you will be willing to embark upon and find more esoteric connections and hopefully one day, with some hard work and a lot of grace, we may find that factorization was always in O(1), Uff!
\\
The property of $\Delta$ sieve zones being yielded when $\sum\limits_{i=x}^{k} {a}_i {df}_i = 0$ for $x, k \in \mathbb{N}$, $a \in \mathbb{Z}$ and \\ \\ the respective $\Delta$ sieve values getting expressed as: $c_1 \Delta^2 + c_2 \Delta +  k$, where $c_1, c_2, k$ are constants for $c_1 , c_2 \in \mathbb{Q}$, $k\in \mathbb{Z}$ for \textbf{any} $\Delta$, changed our world. It meant, we could reuse the effort spent in factorizing a composite $n$ to factorize other composite numbers. In addition, the property of these $\Delta$ sieve zones to switch into different zones with changing $dial_2$ gave us another powerful instrument in the toolchain.  We now have two tools with wonderful ability and flexibility to be executed in parallel and using only fundamental arithmetic operations to sieve the $\Delta$
\\
\\
In our almost decade long search, the first ray of hope came in the form of observing one $a\Delta$ relationship, but it was $ro\{X,Y\}$ which gave us the required motivation to continue with the pursuit, as we took the pulsating "As" as some form of sign from the nature.
\\
\\
There is tremendous life philosophy in $\Delta$ sieve zones, if one looked at them with that lens. We could see how nature is trying to communicate with us through numbers and saying, once we know and acknowledge the difference(s), either with others or with ourselves, the chaos ends and we attain the steady state. 
\\
\\
Finally, it was the equilibrium of $\Delta_{|p-q|} \rightleftharpoons \sum_{p+q}$ which took us back in time, some 13.8 billion years ago, right at that "moment", and gave us some insights and ideas on what might have had happened and helped us in providing an explanation of where all the antimatter disappeared. This research started by taking inspiration from Physics and it was the payback time. We had to make an assumption, a mighty one, but not an absurd one, an assumption which some of us may be able to relate to. What iff, there was a multiplicative force and under its influence, two antimatter particles ended up creating the matter, like when we multiply two negative numbers and we get a positive number? So, while we are searching for antimatter all around us and not able to find much, Sherlock says, we are not able to find it because the matter we see today was made from antimatter. The actual matter is somewhere else, the real one, the one which doesn't obey the laws of observable "antimatter" universe ... now this is getting absurd!
\\
\subsection{Applications}
\begin{enumerate}
  \item \textbf{Integer Factorization}
  \begin{enumerate}[noitemsep]
  	\item With appropriate selection of dials and observation decks leading to sufficient number of $\Delta$ sieve zones, one can factorize the composite $n$, as described in section 4 - "Steady State Value" and is a general purpose factorization method \\
	\item ro\{X,Y\} is also a general purpose factorization method and will lead to nontrivial factorization of the composite $n$, as described in subsection 5.3 - "Reflection over \{X,Y\}" \\
	\item $a\Delta$ relationships are special purpose and in certain instances will lead to nontrivial factorization of the composite $n$, as described in subsection 5.1 - "od connect" and subsection 5.2 - "Range of previous/next $n$" \\
	\item $r\Delta$ relationships are also special purpose and in certain instances will lead to nontrivial factorization of the composite $n$, as described in section 6 - "Inter $\Delta$ relationships (r$\Delta$)" \\
	\item When some $\Delta \rightleftharpoons \sum$ constant is added to $od_6$ on the $\Delta$ series and if it yields $N$ on the corresponding $\sum$ series, then this will lead to nontrivial factorization of the composite $n$, as described in section 8 - "The equilibrium of $\Delta_{|p-q|}$ and $\sum_{p+q}$"
  \end{enumerate}
\item \textbf{RSA Strengthening} \\
	When the composite $n$ is obtained from two random primes, we may need to make sure this random $n$ is RSA-safe by passing it through the factorization methods discussed in this paper to eliminate the possibility of nontrivial factorization, thereby leading to strengthening of RSA
\item \textbf{RSA - Master Key} ($known\ \Delta \simeq unknown \{\Delta_1, \Delta_2, \Delta_3, ..., \Delta_x\}$)\\
$r\Delta$ (with known $\Delta$ and unknown $\Delta$s) allows for creation of an RSA-master key in form of known $\Delta$, as described in section 6 - "Inter $\Delta$ relationships (r$\Delta$)". Information about dials and observation deck sieve formulae will be part of the private key.
   \begin{enumerate}[noitemsep]
   	\item The RSA  master key is likely to have a smaller private key footprint when a message is cyclically encrypted many times over with different public keys \\
	\item Different public keys can be used to encrypt the message which can be decrypted with the same private key as long as they are linked with the same known $\Delta$
   \end{enumerate}
\item \textbf{The New Trapdoor}\\
Section 9 - "The Trapdoor", introduces a new trapdoor function. The $\Delta \rightleftharpoons \sum$ equilibrium provides the foundational framework to build message encryption/decryption schemes, digital signatures and authentication schemes.
\end{enumerate}
\subsection{Further Research}
There are some questions at this point that requires further research.
\begin{enumerate}
\item By having as many fixed dials and observation decks as required, which can yield the required number of $\Delta$ sieve zones and \underline{without} changing $dial_2$, will we able to sieve any $\Delta$ in polynomial time?
\item By having as many fixed dials and observation decks as required, which can yield the required number of $\Delta$ sieve zones and \underline{with} changing $dial_2$, will we able to sieve any $\Delta$ in polynomial time?
\item Will the composite $n$, calculated from two large random primes $p$ and $q$, used in RSA be required to calculate a safety index to confirm whether it is safe from factorization by different factorization methods presented in this paper?
\item When different semiprimes are calculated by multiplying two large random primes $p$ and $q$ in RSA, what $\Delta$ sieve zones do these semiprimes are found for their respective $\Delta$s? Is there a pattern to randomness, whereby the probability for certain kind of $\Delta$ sieve zones to appear is higher than the others?
\item Can there be an equivalent quantum algorithm for methods described in this paper? If yes, will it require less or more qubits in comparison to Shor's algorithm \cite{shor} for integer factorization?
\item Will there always be some dial settings that will yield $zone_0$ in any observation deck? Note: $zone_0$ is the zone where $\Delta$ can be sieved and will continue until $\infty$
\item What is the angle of reflection on rising or falling edge of intra-A or inter-A to observe a nontrivial reflection? Does it fall within any range? and whether knowing this angle of reflection provide any advantage in the nontrivial factorization of the composite $n$?
\item What are the different counting functions looking like?
\begin{enumerate}
\item Number of $zones$ for any $\Delta$?
\item Number of $zoners$ and $zoneless$ for any $\Delta$ and for all possible $dial$ combinations?
\item How many different $dial_1$ combinations exists such that the $\Delta$ sieve coverage zones doesn't change with changing $dial_2$ for any $\Delta$?
\end{enumerate}
\item Can we find similar patterns emerge when we apply the ideas in this paper to discrete logarithms?
\item Can lattices be used in anyway to further speed-up the search of $\Delta$ sieve zones?
\end{enumerate}
\clearpage\
\appendix
\section {Tables}
\begin{table}[!htbp] 
 \caption{$\Delta_{|p-q|}=12;\ p = 2j+1,\ j \in \mathbb{N}_0;\ q = p + \Delta$; $dial_P=\{0,-1,2,2\}$}
  \label{table:A1}
\begin{tabular}{|c|c|} 
\hline 
0 & -1 \\ 
\Xhline{0.5\arrayrulewidth}
2 & 2 										 						\\
\Xhline{0.5\arrayrulewidth}
\end{tabular}
   $dial\ settings$
\\
\begin{tabular}{|c|c|c|c|c|c|c|c|c|c|c|c|c|c|c|}
\hline\hline
$id$ & $p$ & $q$ & $n=pq$ & $\lfloor \sqrt{n} \rfloor$ & $d_1$ & $d_2$ & $d_1^2$ & $d_2^2$ & $od_1$ & $od_2$ & $od_3$ & $od_4$ & $od_6$ \\[0.5ex]
\hline
1 & 1 & 13 & 13 & 3 & 2 & 4 & 4 & 16 & -9 & 3 & 12 & -6 & 5 \\
\Xhline{0.5\arrayrulewidth}
2 & 3 & 15 & 45 & 6 & 6 & 8 & 36 & 64 & -9 & 19 & 28 & 10 & 3 \\
\Xhline{0.5\arrayrulewidth}
3 & 5 & 17 & 85 & 9 & 8 & 10 & 64 & 100 & -21 & 15 & 36 & -6 & 5 \\
\Xhline{0.5\arrayrulewidth}
4 & 7 & 19 & 133 & 11 & 10 & 12 & 100 & 144 & -33 & 11 & 44 & -22 & 13 \\
\Xhline{0.5\arrayrulewidth}
5 & 9 & 21 & 189 & 13 & 12 & 14 & 144 & 196 & -45 & 7 & 52 & -38 & 21 \\
\Xhline{0.5\arrayrulewidth}
6 & 11 & 23 & 253 & 15 & 14 & 16 & 196 & 256 & -57 & 3 & 60 & -54 & 29 \\
\Xhline{0.5\arrayrulewidth}
\textbf{7} & \textbf{13} & 25 & 325 & 18 & 18 & 20 & 324 & 400 & -1 & 75 & 76 & \textbf{74} & \textbf{35} \\
\Xhline{0.5\arrayrulewidth}
8 & 15 & 27 & 405 & 20 & 20 & 22 & 400 & 484 & -5 & 79 & 84 & \textbf{74} & \textbf{35} \\
\Xhline{0.5\arrayrulewidth}
 ... & ... & ... & ... & ... & ... & ... & ... & ...  & ... & ... & ... & \textbf{74}  & \textbf{35} \\
\Xhline{0.5\arrayrulewidth}
\begin{math}\infty\end{math} & \begin{math}\infty\end{math}  &  \begin{math}\infty\end{math} &  \begin{math}\infty\end{math}  & ... & ... & ... & ... & ...  & ... & ... & ... & \textbf{74}  & \textbf{35} \\ [1ex] 
\Xhline{0.5\arrayrulewidth}
\end{tabular}
\end{table}
\begin{table}[!htbp] 
  \caption{$\Delta_{|p-q|}=22;\ p = 2j+1,\ j \in \mathbb{N}_0;\ q = p + \Delta$; $dial_P=\{0,-1,2,2\}$}
  \label{table:A2}
\begin{tabular}{|c|c|}
\hline
0 & -1 \\
\Xhline{0.5\arrayrulewidth}
2 & 2 										 						\\
\Xhline{0.5\arrayrulewidth}
\end{tabular}
   $dial\ settings$
\\
\\
\begin{tabular}{|c|c|c|c|c|c|c|c|c|c|c|c|c|c|c|}
\hline\hline
$id$ & $p$ & $q$ & $n=pq$ & $\lfloor \sqrt{n} \rfloor$ & $d_1$ & $d_2$ & $d_1^2$ & $d_2^2$ & $od_1$ & $od_2$ & $od_3$ & $od_4$ & $od_6$ \\[0.5ex]
\hline
1 & 1 & 23 & 23 & 4 & 4 & 6 & 16 & 36 & -7 & 13 & 20 & 6 & 1 \\
\Xhline{0.5\arrayrulewidth}
... & ... & ... & ... & ... & ... & ... & ... & ... & ... & ... & ... & ... & ... \\
\Xhline{0.5\arrayrulewidth}
... & ... & ... & ... & ... & ... & ... & ... & ... & ... & ... & ... & ... & ... \\
\Xhline{0.5\arrayrulewidth}
8 & 15 & 37 & 555 & 23 & 22 & 24 & 484 & 576 & -71 & 21 & 92 & -50 & 27 \\
\Xhline{0.5\arrayrulewidth}
9 & 17 & 39 & 663 & 25 & 24 & 26 & 576 & 676 & -87 & 13 & 100 & -74 & 39 \\
\Xhline{0.5\arrayrulewidth}
10 & 19 & 41 & 779 & 27 & 26 & 28 & 676 & 784 & -103 & 5 & 108 & -98 & 51 \\
\Xhline{0.5\arrayrulewidth}
\textbf{11} & \textbf{21} & 43 & 903 & 30 & 30 & 32 & 900 & 1024 & -3 & \textbf{121} & 124 & 118 & 57 \\
\Xhline{0.5\arrayrulewidth}
12 & 23 & 45 & 1035 & 32 & 32 & 34 & 1024 & 1156 & -11 & \textbf{121} & 132 & 110 & 53 \\
\Xhline{0.5\arrayrulewidth}
 ... & ... & ... & ... & ... & ... & ... & ... & ...  & ... & \textbf{121} & ... & ...  & ...\\
\Xhline{0.5\arrayrulewidth}
\begin{math}\infty\end{math} & \begin{math}\infty\end{math}  &  \begin{math}\infty\end{math} &  \begin{math}\infty\end{math}  & ... & ... & ... & ... & ...  & ... & \textbf{121} & ... & ... & ... \\ [1ex] 
\Xhline{0.5\arrayrulewidth}
\end{tabular}
\end{table}
\begin{table}[!htbp] 
  \caption{$\Delta_{|p-q|}=12;\ p = 2j,\ j \in \mathbb{N};\ q = p + \Delta$; $dial_P=\{0,-1,2,2\}$}
  \label{table:A3}
\begin{tabular}{|c|c|}
\hline
0 & -1 \\
\Xhline{0.5\arrayrulewidth}
2 & 2 										 						\\
\Xhline{0.5\arrayrulewidth}
\end{tabular}
   $dial\ settings$
\\
\\
\begin{tabular}{|c|c|c|c|c|c|c|c|c|c|c|c|c|c|c|}
\hline\hline
$id$ & $p$ & $q$ & $n=pq$ & $\lfloor \sqrt{n} \rfloor$ & $d_1$ & $d_2$ & $d_1^2$ & $d_2^2$ & $od_1$ & $od_2$ & $od_3$ & $od_4$ & $od_6$ \\[0.5ex]
\hline
1 & 2 & 14 & 28 & 5 & 4 & 6 & 16 & 36 & -12 & 8 & 20 & -4 & 4 \\
\Xhline{0.5\arrayrulewidth}
\textbf{2} & \textbf{4} & 16 & 64 & 8 & 8 & 10 & 64 & 100 & 0 & \textbf{36} & 36 & 36 & 16 \\
\Xhline{0.5\arrayrulewidth}
3 & 6 & 18 & 108 & 10 & 10 & 12 & 100 & 144 & -8 & \textbf{36} & 44 & 28 & 12 \\
\Xhline{0.5\arrayrulewidth}
 ... & ... & ... & ... & ... & ... & ... & ... & ...  & ... & \textbf{36} & ... & ...  & ... \\
\Xhline{0.5\arrayrulewidth}
\begin{math}\infty\end{math} & \begin{math}\infty\end{math}  &  \begin{math}\infty\end{math} &  \begin{math}\infty\end{math}  & ... & ... & ... & ... & ...  & ... & \textbf{36} & ... & ...  & ... \\ [1ex]
\Xhline{0.5\arrayrulewidth}
\end{tabular}
\end{table}
\begin{table}[!htbp] 
 \caption{$\Delta_{|p-q|}=16;\ p = 2j,\ j \in  \mathbb{N};\ q = p + \Delta$; $dial_P=\{0,-1,2,2\}$}
  \label{table:A4}
\begin{tabular}{|c|c|}
\hline
0 & -1 \\
\Xhline{0.5\arrayrulewidth}
2 & 2 										 						\\
\Xhline{0.5\arrayrulewidth}
\end{tabular}
   $dial\ settings$
\\
\\
\begin{tabular}{|c|c|c|c|c|c|c|c|c|c|c|c|c|c|c|}
\hline\hline
$id$ & $p$ & $q$ & $n=pq$ & $\lfloor \sqrt{n} \rfloor$ & $d_1$ & $d_2$ & $d_1^2$ & $d_2^2$ & $od_1$ & $od_2$ & $od_3$ & $od_4$ & $od_6$ \\[0.5ex]
\hline
1 & 2 & 18 & 36 & 6 & 6 & 8 & 36 & 64 & 0 & 28 & 28 & 28 & 12 \\
\Xhline{0.5\arrayrulewidth}
2 & 4 & 20 & 80 & 8 & 8 & 10 & 64 & 100 & -16 & 20 & 36 & 4 & 0 \\
\Xhline{0.5\arrayrulewidth}
3 & 6 & 22 & 132 & 11 & 10 & 12 & 100 & 144 & -32 & 12 & 44 & -20 & 12 \\
\Xhline{0.5\arrayrulewidth}
4 & 8 & 24 & 192 & 13 & 12 & 14 & 144 & 196 & -48 & 4 & 52 & -44 & 24 \\
\Xhline{0.5\arrayrulewidth}
\textbf{5} & \textbf{10} & 26 & 260 & 16 & 16 & 18 & 256 & 324 & -4 & \textbf{64} & 68 & 60 & 28 \\
\Xhline{0.5\arrayrulewidth}
6 & 12 & 28 & 336 & 18 & 18 & 20 & 324 & 400 & -12 & \textbf{64} & 76 & 52 & 24 \\
\Xhline{0.5\arrayrulewidth}
 ... & ... & ... & ... & ... & ... & ... & ... & ...  & ... & \textbf{64} & ... & ...  & ...\\
\Xhline{0.5\arrayrulewidth}
\begin{math}\infty\end{math} & \begin{math}\infty\end{math}  &  \begin{math}\infty\end{math} &  \begin{math}\infty\end{math}  & ... & ... & ... & ... & ...  & ... & \textbf{64} & ... & ...  & ... \\ [1ex]
\Xhline{0.5\arrayrulewidth}
\end{tabular}
\vspace{25em}
\end{table}
\begin{table}[!htbp] 
 \caption{$\Delta_{|p-q|}=22;\ p = 2j,\ j \in \mathbb{N};\ q = p + \Delta$; $dial_P=\{0,-1,2,2\}$}
  \label{table:A5}
\begin{tabular}{|c|c|}
\hline
0 & -1 \\
\Xhline{0.5\arrayrulewidth}
2 & 2 										 						\\
\Xhline{0.5\arrayrulewidth}
\end{tabular}
   $dial\ settings$
\\
\\
\begin{tabular}{|c|c|c|c|c|c|c|c|c|c|c|c|c|c|c|}
\hline\hline
$id$ & $p$ & $q$ & $n=pq$ & $\lfloor \sqrt{n} \rfloor$ & $d_1$ & $d_2$ & $d_1^2$ & $d_2^2$ & $od_1$ & $od_2$ & $od_3$ & $od_4$ & $od_6$ \\[0.5ex]
\hline
1 & 2 & 24 & 48 & 6 & 6 & 8 & 36 & 64 & -12 & 16 & 28 & 4 & 0 \\
\Xhline{0.5\arrayrulewidth}
2 & 4 & 26 & 104 & 10 & 10 & 12 & 100 & 144 & -4 & 40 & 44 & 36 & 16 \\
\Xhline{0.5\arrayrulewidth}
3 & 6 & 28 & 168 & 12 & 12 & 14 & 144 & 196 & -24 & 28 & 52 & 4 & 0 \\
\Xhline{0.5\arrayrulewidth}
4 & 8 & 30 & 240 & 15 & 14 & 16 & 196 & 256 & -44 & 16 & 60 & -28 & 16 \\
\Xhline{0.5\arrayrulewidth}
5 & 10 & 32 & 320 & 17 & 16 & 18 & 256 & 324 & -64 & 4 & 68 & -60 & 32 \\
\Xhline{0.5\arrayrulewidth}
6 & 12 & 34 & 408 & 20 & 20 & 22 & 400 & 484 & -8 & 76 & 84 & 68 & 32 \\
\Xhline{0.5\arrayrulewidth}
7 & 14 & 36 & 504 & 22 & 22 & 24 & 484 & 576 & -20 & 72 & 92 & 52 & 24 \\
\end{tabular}
\vspace{25em}
\end{table}
\begin{table}[!htbp] 
 \caption{$\Delta_{|p-q|}=24;\ p = 2j+1,\ j \in \mathbb{N}_0;\ q = p + \Delta$; $dial_P=\{-1,0,2,2\}$}
  \label{table:B1}
\begin{tabular}{|c|c|}
\hline
-1 & 0 \\
\Xhline{0.5\arrayrulewidth}
2 & 2 										 						\\
\Xhline{0.5\arrayrulewidth}
\end{tabular}
   $dial\ settings$
\\
\\
\begin{tabular}{|c|c|c|c|c|c|c|c|c|c|c|c|c|c|c|}
\hline\hline
$id$ & $p$ & $q$ & $n=pq$ & $\lfloor \sqrt{n} \rfloor$ & $d_1$ & $d_2$ & $d_1^2$ & $d_2^2$ & $od_1$ & $od_2$ & $od_3$ & $od_4$ & $od_6$ \\[0.5ex]
\hline
1 & 1 & 25 & 25 & 5 & 5 & 7 & 25 & 49 & 0 & 24 & 24 & 24 & 10 \\
\Xhline{0.5\arrayrulewidth}
2 & 3 & 27 & 81 & 9 & 9 & 11 & 81 & 121 & 0 & 40 & 40 & 40 & 18 \\
\Xhline{0.5\arrayrulewidth}
3 & 5 & 29 & 145 & 12 & 11 & 13 & 121 & 169 & -24 & 24 & 48 & 0 & 2 \\
\Xhline{0.5\arrayrulewidth}
4 & 7 & 31 & 217 & 14 & 13 & 15 & 169 & 225 & -48 & 8 & 56 & -40 & 22 \\
\Xhline{0.5\arrayrulewidth}
5 & 9 & 33 & 297 & 17 & 17 & 19 & 289 & 361 & -8 & 64 & 72 & 56 & 26 \\
\Xhline{0.5\arrayrulewidth}
6 & 11 & 35 & 385 & 19 & 19 & 21 & 361 & 441 & -24 & 56 & 80 & 32 & 14 \\
\Xhline{0.5\arrayrulewidth}
7 & 13 & 37 & 481 & 21 & 21 & 23 & 441 & 529 & -40 & 48 & 88 & 8 & 2 \\
\end{tabular}
\vspace{10em}
\end{table}
\begin{table}[!htbp] 
 \caption{$\Delta_{|p-q|}=12;\ p = 2j+1,\ j \in \mathbb{N}_0;\ q = p + \Delta$; $dial_P=\{-1,0,2,2\}$}
  \label{table:B2}
\begin{tabular}{|c|c|}
\hline
-1 & 0 \\
\Xhline{0.5\arrayrulewidth}
2 & 2 										 						\\
\Xhline{0.5\arrayrulewidth}
\end{tabular}
   $dial\ settings$
\\
\\
\begin{tabular}{|c|c|c|c|c|c|c|c|c|c|c|c|c|c|c|}
\hline\hline
$id$ & $p$ & $q$ & $n=pq$ & $\lfloor \sqrt{n} \rfloor$ & $d_1$ & $d_2$ & $d_1^2$ & $d_2^2$ & $od_1$ & $od_2$ & $od_3$ & $od_4$ & $od_6$ \\[0.5ex]
\hline
1 & 1 & 13 & 13 & 3 & 3 & 5 & 9 & 25 & -4 & 12 & 16 & 8 & 2 \\
\Xhline{0.5\arrayrulewidth}
2 & 3 & 15 & 45 & 6 & 5 & 7 & 25 & 49 & -20 & 4 & 24 & -16 & 10 \\
\Xhline{0.5\arrayrulewidth}
\textbf{3} & \textbf{5} & 17 & 85 & 9 & 9 & 11 & 81 & 121 & -4 & \textbf{36} & 40 & 32 & 14 \\
\Xhline{0.5\arrayrulewidth}
4 & 7 & 19 & 133 & 11 & 11 & 13 & 121 & 169 & -12 & \textbf{36} & 48 & 24 & 10 \\
\Xhline{0.5\arrayrulewidth}
 ... & ... & ... & ... & ... & ... & ... & ... & ...  & ... & \textbf{36} & ... & ...  & ...\\
\Xhline{0.5\arrayrulewidth}
\begin{math}\infty\end{math} & \begin{math}\infty\end{math}  &  \begin{math}\infty\end{math} &  \begin{math}\infty\end{math}  & ... & ... & ... & ... & ...  & ... & \textbf{36} & ... & ...  & ... \\ [1ex]
\Xhline{0.5\arrayrulewidth}
\end{tabular}
\end{table}
\begin{table}[!htbp] 
 \caption{$\Delta_{|p-q|}=22;\ p = 2j+1,\ j \in \mathbb{N}_0;\ q = p + \Delta$; $dial_P=\{-1,0,2,2\}$}
  \label{table:B3}
\begin{tabular}{|c|c|}
\hline
-1 & 0 \\
\Xhline{0.5\arrayrulewidth}
2 & 2 										 						\\
\Xhline{0.5\arrayrulewidth}
\end{tabular}
   $dial\ settings$
\\
\\
\begin{tabular}{|c|c|c|c|c|c|c|c|c|c|c|c|c|c|c|}
\hline\hline
$id$ & $p$ & $q$ & $n=pq$ & $\lfloor \sqrt{n} \rfloor$ & $d_1$ & $d_2$ & $d_1^2$ & $d_2^2$ & $od_1$ & $od_2$ & $od_3$ & $od_4$ & $od_6$ \\[0.5ex]
\hline
1 & 1 & 23 & 23 & 4 & 3 & 5 & 9 & 25 & -14 & 2 & 16 & -12 & 8 \\
\Xhline{0.5\arrayrulewidth}
2 & 3 & 25 & 75 & 8 & 7 & 9 & 49 & 81 & -26 & 6 & 32 & -20 & 12 \\
\Xhline{0.5\arrayrulewidth}
3 & 5 & 27 & 135 & 11 & 11 & 13 & 121 & 169 & -14 & 34 & 48 & 20 & 8 \\
\Xhline{0.5\arrayrulewidth}
4 & 7 & 29 & 203 & 14 & 13 & 15 & 169 & 225 & -34 & 22 & 56 & -12 & 8 \\
\Xhline{0.5\arrayrulewidth}
5 & 9 & 31 & 279 & 16 & 15 & 17 & 225 & 289 & -54 & 10 & 64 & -44 & 24 \\
\Xhline{0.5\arrayrulewidth}
6 & 11 & 33 & 363 & 19 & 19 & 21 & 361 & 441 & -2 & 78 & 80 & 76 & 36 \\
\Xhline{0.5\arrayrulewidth}
7 & 13 & 35 & 455 & 21 & 21 & 23 & 441 & 529 & -14 & 74 & 88 & 60 & 28 \\
\end{tabular}
\end{table}
\begin{table}[!htbp] 
 \caption{$\Delta_{|p-q|}=12;\ p = 2j,\ j \in \mathbb{N};\ q = p + \Delta$; $dial_P=\{-1,0,2,2\}$}
  \label{table:B4}
\begin{tabular}{|c|c|}
\hline
-1 & 0 \\
\Xhline{0.5\arrayrulewidth}
2 & 2 										 						\\
\Xhline{0.5\arrayrulewidth}
\end{tabular}
   $dial\ settings$
\\
\\
\begin{tabular}{|c|c|c|c|c|c|c|c|c|c|c|c|c|c|c|}
\hline\hline
$id$ & $p$ & $q$ & $n=pq$ & $\lfloor \sqrt{n} \rfloor$ & $d_1$ & $d_2$ & $d_1^2$ & $d_2^2$ & $od_1$ & $od_2$ & $od_3$ & $od_4$ & $od_6$ \\[0.5ex]
\hline
1 & 2 & 14 & 28 & 5 & 5 & 7 & 25 & 49 & -3 & 21 & 24 & 18 & 7 \\
\Xhline{0.5\arrayrulewidth}
2 & 4 & 16 & 64 & 8 & 7 & 9 & 49 & 81 & -15 & 17 & 32 & 2 & 1 \\
\Xhline{0.5\arrayrulewidth}
3 & 6 & 18 & 108 & 10 & 9 & 11 & 81 & 121 & -27 & 13 & 40 & -14 & 9 \\
\Xhline{0.5\arrayrulewidth}
4 & 8 & 20 & 160 & 12 & 11 & 13 & 121 & 169 & -39 & 9 & 48 & -30 & 17 \\
\Xhline{0.5\arrayrulewidth}
5 & 10 & 22 & 220 & 14 & 13 & 15 & 169 & 225 & -51 & 5 & 56 & -46 & 25 \\
\Xhline{0.5\arrayrulewidth}
6 & 12 & 24 & 288 & 16 & 15 & 17 & 225 & 289 & -63 & 1 & 64 & -62 & 33 \\
\Xhline{0.5\arrayrulewidth}
\textbf{7} & \textbf{14} & 26 & 364 & 19 & 19 & 21 & 361 & 441 & -3 & 77 & 80 &  \textbf{74} &  \textbf{35} \\
\Xhline{0.5\arrayrulewidth}
8 & 16 & 28 & 448 & 21 & 21 & 23 & 441 & 529 & -7 & 81 & 88 &  \textbf{74} &  \textbf{35} \\
\Xhline{0.5\arrayrulewidth}
... & ... & ... & ... & ... & ... & ... & ... & ... & ... & ... & ... &  \textbf{74} &  \textbf{35} \\
\Xhline{0.5\arrayrulewidth}
$\infty$ & $\infty$ & $\infty$ & $\infty$ & ... & ... & ... & ... & ... & ... & ... & ... &  \textbf{74} &  \textbf{35} \\
\Xhline{0.5\arrayrulewidth}
\end{tabular}
\end{table}

\begin{table}[!htbp] 
 \caption{$\Delta_{|p-q|}=22;\ p = 2j,\ j \in \mathbb{N};\ q = p + \Delta$; $dial_P=\{-1,0,2,2\}$}
  \label{table:B5}
\begin{tabular}{|c|c|}
\hline
-1 & 0 \\
\Xhline{0.5\arrayrulewidth}
2 & 2 										 						\\
\Xhline{0.5\arrayrulewidth}
\end{tabular}
   $dial\ settings$
\\
\\
\begin{tabular}{|c|c|c|c|c|c|c|c|c|c|c|c|c|c|c|}
\hline\hline
$id$ & $p$ & $q$ & $n=pq$ & $\lfloor \sqrt{n} \rfloor$ & $d_1$ & $d_2$ & $d_1^2$ & $d_2^2$ & $od_1$ & $od_2$ & $od_3$ & $od_4$ & $od_6$ \\[0.5ex]
\hline
1 & 2 & 24 & 48 & 6 & 5 & 7 & 25 & 49 & -23 & 1 & 24 & -22 & 13 \\
\Xhline{0.5\arrayrulewidth}
2 & 4 & 26 & 104 & 10 & 9 & 11 & 81 & 121 & -23 & 17 & 40 & -6 & 5 \\
\Xhline{0.5\arrayrulewidth}
3 & 6 & 28 & 168 & 12 & 11 & 13 & 121 & 169 & -47 & 1 & 48 & -46 & 25 \\
\Xhline{0.5\arrayrulewidth}
4 & 8 & 30 & 240 & 15 & 15 & 17 & 225 & 289 & -15 & 49 & 64 & 34 & 15 \\
\Xhline{0.5\arrayrulewidth}
5 & 10 & 32 & 320 & 17 & 17 & 19 & 289 & 361 & -31 & 41 & 72 & 10 & 3 \\
\Xhline{0.5\arrayrulewidth}
6 & 12 & 34 & 408 & 20 & 19 & 21 & 361 & 441 & -47 & 33 & 80 & -14 & 9 \\
\Xhline{0.5\arrayrulewidth}
7 & 14 & 36 & 504 & 22 & 21 & 23 & 441 & 529 & -63 & 25 & 88 & -38 & 21 \\
\Xhline{0.5\arrayrulewidth}
8 & 16 & 38 & 608 & 24 & 23 & 25 & 529 & 625 & -79 & 17 & 96 & -62 & 33 \\
\Xhline{0.5\arrayrulewidth}
9 & 18 & 40 & 720 & 26 & 25 & 27 & 625 & 729 & -95 & 9 & 104 & -86 & 45 \\
\Xhline{0.5\arrayrulewidth}
10 & 20 & 42 & 840 & 28 & 27 & 29 & 729 & 841 & -111 & 1 & 112 & -110 & 57 \\
\Xhline{0.5\arrayrulewidth}
\textbf{11} & \textbf{22} & 44 & 968 & 31 & 31 & 33 & 961 & 1089 & -7 & \textbf{121} & 128 & 114 & 55 \\
\Xhline{0.5\arrayrulewidth}
12 & 24 & 46 & 1104 & 33 & 33 & 35 & 1089 & 1225 & -15 & \textbf{121} & 136 & 106 & 51 \\
\Xhline{0.5\arrayrulewidth}
... & ... & ... & ... & ... & ... & ... & ... & ... & ... & \textbf{121} & ... & ... & ... \\
\Xhline{0.5\arrayrulewidth}
$\infty$ & $\infty$ & $\infty$ & $\infty$ & ... & ... & ... & ... & ... & ... & \textbf{121} & ... & ... & ... \\
\Xhline{0.5\arrayrulewidth}
\end{tabular}
\end{table}\
\begin{table}[!htbp] 
 \caption{$\Delta_{|p-q|}=12;\ p = 2j+1,\ j \in \mathbb{N}_0;\ q = p + \Delta$; $dial_P=\{0,-1,2,2\}$}
 \label{table:C1}
\begin{tabular}{|c|c|}
\hline 
0 & -1 \\ 
\Xhline{0.5\arrayrulewidth}
2 & 2 										 						\\
\Xhline{0.5\arrayrulewidth}
\end{tabular}
   $dial\ settings$
\\
\\
\begin{tabular}{|c|c|c|c|c|c|c|c|c|c|c|c|c|c|c|c|}
\hline\hline
$id$ & $p$ & $q$ & $n=pq$ & $od_1$ & $df_1$ & $od_2$ & $df_2$ & $od_3$ & $df_3$ & $od_4$ & $df_4$ & $zone$\\[0.5ex]
\hline
\Xhline{0.5\arrayrulewidth}
1 & 1 & 13 & 13 & -9 & 0 & 3 & 0 & 12 & 0 & -6 & 0 & ... \\
\Xhline{0.5\arrayrulewidth}
2 & 3 & 15 & 45 & -9 & 0 & 19 & -16 & 28 & -16 & 10 & -16 & ... \\
\Xhline{0.5\arrayrulewidth}
3 & 5 & 17 & 85 & -21 & 12 & 15 & 4 & 36 & -8 & -6 & 16 & ... \\
\Xhline{0.5\arrayrulewidth}
4 & 7 & 19 & 133 & -33 & 12 & 11 & 4 & 44 & -8 & -22 & 16 & ... \\
\Xhline{0.5\arrayrulewidth}
5 & 9 & 21 & 189 & -45 & 12 & 7 & 4 & 52 & -8 & -38 & 16 & ... \\
\Xhline{0.5\arrayrulewidth}
6 & 11 & 23 & 253 & -57 & 12 & 3 & 4 & 60 & -8 & -54 & 16 & ... \\
\Xhline{0.5\arrayrulewidth}
7 & 13 & 25 & 325 & -1 & \textbf{-56} & 75 & \textbf{-72} & 76 & -16 & \textbf{74} & \textbf{-128} & \textbf{0} \\
\end{tabular}
\vspace{30em}
\end{table}
\begin{table}[!htbp] 
 \caption{$\Delta_{|p-q|}=22;\ p = 2j,\ j \in \mathbb{N};\ q = p + \Delta$; $dial_P=\{0,-1,2,2\}$}
 \label{table:C2}
\begin{tabular}{|c|c|}
\hline 
0 & -1 \\ 
\Xhline{0.5\arrayrulewidth}
2 & 2 										 						\\
\Xhline{0.5\arrayrulewidth}
\end{tabular}
   $dial\ settings$
\\
\begin{tabular}{|c|c|c|c|c|c|c|c|c|c|c|c|c|c|c|c|}
\hline\hline
$id$ & $p$ & $q$ & $n=pq$ & $od_1$ & $df_1$ & $od_2$ & $df_2$ & $od_3$ & $df_3$ & $od_4$ & $df_4$ & $zone$\\[0.5ex]
\hline
\Xhline{0.5\arrayrulewidth}
1 & 2 & 24 & 48 & -12 & 0 & 16 & 0 & 28 & 0 & 4 & 0 & ... \\
\Xhline{0.5\arrayrulewidth}
2 & 4 & 26 & 104 & -4 & -8 & 40 & -24 & 44 & -16 & 36 & -32 & ... \\
\Xhline{0.5\arrayrulewidth}
3 & 6 & 28 & 168 & -24 & 20 & 28 & 12 & 52 & -8 & 4 & 32 & ... \\
\Xhline{0.5\arrayrulewidth}
4 & 8 & 30 & 240 & -44 & 20 & 16 & 12 & 60 & -8 & -28 & 32 & ... \\
\Xhline{0.5\arrayrulewidth}
5 & 10 & 32 & 320 & -64 & 20 & 4 & 12 & 68 & -8 & -60 & 32 & ... \\
\Xhline{0.5\arrayrulewidth}
6 & 12 & 34 & 408 & -8 & -56 & 76 & -72 & 84 & -16 & 68 & -128 & ... \\
\Xhline{0.5\arrayrulewidth}
7 & 14 & 36 & 504 & -20 & 12 & 72 & 4 & 92 & -8 & 52 & 16 & ... \\
\end{tabular}
\vspace{30em}
\end{table}
\begin{table}[!htbp] 
 \caption{$\Delta_{|p-q|}=22;\ p = 2j+1,\ j \in \mathbb{N}_0;\ q = p + \Delta$; $dial_P=\{-1,0,2,2\}$}
 \label{table:C3}
\begin{tabular}{|c|c|}
\hline 
-1 & 0 \\ 
\Xhline{0.5\arrayrulewidth}
2 & 2 										 						\\
\Xhline{0.5\arrayrulewidth}
\end{tabular}
   $dial\ settings$
\\
\begin{tabular}{|c|c|c|c|c|c|c|c|c|c|c|c|c|c|c|c|}
\hline\hline
$id$ & $p$ & $q$ & $n=pq$ & $od_1$ & $df_1$ & $od_2$ & $df_2$ & $od_3$ & $df_3$ & $od_4$ & $df_4$ & $zone$\\[0.5ex]
\hline
\Xhline{0.5\arrayrulewidth}
1 & 1 & 23 & 23 & -14 & 0 & 2 & 0 & 16 & 0 & -12 & 0 & ... \\
\Xhline{0.5\arrayrulewidth}
2 & 3 & 25 & 75 & -26 & 12 & 6 & -4 & 32 & -16 & -20 & 8 & ... \\
\Xhline{0.5\arrayrulewidth}
3 & 5 & 27 & 135 & -14 & -12 & 34 & -28 & 48 & -16 & 20 & -40 & ... \\
\Xhline{0.5\arrayrulewidth}
4 & 7 & 29 & 203 & -34 & 20 & 22 & 12 & 56 & -8 & -12 & 32 & ... \\
\Xhline{0.5\arrayrulewidth}
5 & 9 & 31 & 279 & -54 & 20 & 10 & 12 & 64 & -8 & -44 & 32 & ... \\
\Xhline{0.5\arrayrulewidth}
6 & 11 & 33 & 363 & -2 & -52 & 78 & -68 & 80 & -16 & 76 & -120 & ... \\
\Xhline{0.5\arrayrulewidth}
7 & 13 & 35 & 455 & -14 & 12 & 74 & 4 & 88 & -8 & 60 & 16 & ... \\
\end{tabular}
\vspace{15em}
\end{table}
\begin{table}[!htbp] 
 \caption{$\Delta_{|p-q|}=12;\ p = 2j,\ j \in \mathbb{N};\ q = p + \Delta$; $dial_P=\{-1,0,2,2\}$}
 \label{table:C4}
\begin{tabular}{|c|c|}
\hline 
-1 & 0 \\ 
\Xhline{0.5\arrayrulewidth}
2 & 2 										 						\\
\Xhline{0.5\arrayrulewidth}
\end{tabular}
   $dial\ settings$
\\
\begin{tabular}{|c|c|c|c|c|c|c|c|c|c|c|c|c|c|c|c|}
\hline\hline
$id$ & $p$ & $q$ & $n=pq$ & $od_1$ & $df_1$ & $od_2$ & $df_2$ & $od_3$ & $df_3$ & $od_4$ & $df_4$ & $zone$\\[0.5ex]
\hline
\Xhline{0.5\arrayrulewidth}
1 & 2 & 14 & 28 & -3 & 0 & 21 & 0 & 24 & 0 & 18 & 0 & ... \\
\Xhline{0.5\arrayrulewidth}
2 & 4 & 16 & 64 & -15 & 12 & 17 & 4 & 32 & -8 & 2 & 16 & ... \\
\Xhline{0.5\arrayrulewidth}
3 & 6 & 18 & 108 & -27 & 12 & 13 & 4 & 40 & -8 & -14 & 16 & ... \\
\Xhline{0.5\arrayrulewidth}
4 & 8 & 20 & 160 & -39 & 12 & 9 & 4 & 48 & -8 & -30 & 16 & ... \\
\Xhline{0.5\arrayrulewidth}
5 & 10 & 22 & 220 & -51 & 12 & 5 & 4 & 56 & -8 & -46 & 16 & ... \\
\Xhline{0.5\arrayrulewidth}
6 & 12 & 24 & 288 & -63 & 12 & 1 & 4 & 64 & -8 & -62 & 16 & ... \\
\Xhline{0.5\arrayrulewidth}
7 & 14 & 26 & 364 & -3 & \textbf{-60} & 77 & \textbf{-76} & 80 & -16 & \textbf{74} & \textbf{-136} & \textbf{0} \\
\Xhline{0.5\arrayrulewidth}
8 & 16 & 28 & 448 & -7 & \textbf{4} & 81 & \textbf{-4} & 88 & -8 & \textbf{74} & \textbf{0} & \textbf{0} \\
\Xhline{0.5\arrayrulewidth}
9 & 18 & 30 & 540 & -11 & \textbf{4} & 85 & \textbf{-4} & 96 & -8 & \textbf{74} & \textbf{0} & \textbf{0} \\
\Xhline{0.5\arrayrulewidth}
10 & 20 & 32 & 640 & -15 & \textbf{4} & 89 & \textbf{-4} & 104 & -8 & \textbf{74} & \textbf{0} & \textbf{0} \\
\Xhline{0.5\arrayrulewidth}
11 & 22 & 34 & 748 & -19 & \textbf{4} & 93 & \textbf{-4} & 112 & -8 & \textbf{74} & \textbf{0} & \textbf{0} \\
\Xhline{0.5\arrayrulewidth}
 ... & ... & ... & ... & ... & \textbf{4} & ... & \textbf{-4} & ...  & ... & \textbf{74} & \textbf{0} & \textbf{0} \\
\Xhline{0.5\arrayrulewidth}
 $\infty$ & $\infty$ & $\infty$ & $\infty$ & ... & \textbf{4} & ... & \textbf{-4} & ...  & ... & \textbf{74} & \textbf{0} & \textbf{0} \\
\Xhline{0.5\arrayrulewidth}
\end{tabular}
\end{table}
\begin{table}[!htbp] 
 \caption{$\Delta_{|p-q|}=46;\ p = 2j + 1,\ j \in \mathbb{N}_0;\ q = p + \Delta$; $dial_P=\{-1,-2,4,4\}$}
  \label{table:15}
\begin{tabular}{|c|c|} 
\hline 
-1 & -2 \\ 
\Xhline{0.5\arrayrulewidth}
4 & 4 										 						\\
\Xhline{0.5\arrayrulewidth}
\end{tabular}
   $dial\ settings$
\\
\\
\begin{tabular}{|c|c|c|c|c|c|c|c|c|c|c|c|c|c|c|c|}
\hline\hline
$id$ & $p$ & $q$ & $n=pq$ & $od_1$ & $df_1$ & $od_2$ & $df_2$ & $od_3$ & $df_3$ & $od_4$ & $df_4$ & $od_5$ \\[0.5ex]
\hline 
\Xhline{0.5\arrayrulewidth}
1 & 1 & 47 & 47 & -22 & 0 & 34 & 0 & 56 & 0 & 12 & 0 & 80 \\
\Xhline{0.5\arrayrulewidth}
2 & 3 & 49 & 147 & -26 & 4 & 78 & -44 & 104 & -48 & 52 & -40 & 208 \\
\Xhline{0.5\arrayrulewidth}
3 & 5 & 51 & 255 & -86 & 60 & 34 & 44 & 120 & -16 & -52 & 104 & 16 \\
\Xhline{0.5\arrayrulewidth}
4 & 7 & 53 & 371 & -82 & -4 & 70 & -36 & 152 & -32 & -12 & -40 & 128 \\
\Xhline{0.5\arrayrulewidth}
5 & 9 & 55 & 495 & -54 & -28 & 130 & -60 & 184 & -32 & 76 & -88 & 336 \\
\Xhline{0.5\arrayrulewidth}
6 & 11 & 57 & 627 & -98 & 44 & 102 & 28 & 200 & -16 & 4 & 72 & 208 \\
\Xhline{0.5\arrayrulewidth}
7 & 13 & 59 & 767 & -142 & 44 & 74 & 28 & 216 & -16 & -68 & 72 & 80 \\
\end{tabular}
\end{table}
\begin{table}[!htbp] 
 \caption{$\Delta_{|p-q|}=48;\ p = 2j + 1,\ j \in \mathbb{N}_0;\ q = p + \Delta$; $dial_P=\{-2,-1,4,4\}$}
  \label{table:16}
\begin{tabular}{|c|c|} 
\hline 
-2 & -1 \\ 
\Xhline{0.5\arrayrulewidth}
4 & 4 										 						\\
\Xhline{0.5\arrayrulewidth}
\end{tabular}
   $dial\ settings$
\\
\\
\begin{tabular}{|c|c|c|c|c|c|c|c|c|c|c|c|c|c|c|c|}
\hline\hline
$id$ & $p$ & $q$ & $n=pq$ & $od_1$ & $df_1$ & $od_2$ & $df_2$ & $od_3$ & $df_3$ & $od_4$ & $df_4$ & $od_5$ \\[0.5ex]
\hline 
\Xhline{0.5\arrayrulewidth}
1 & 1 & 49 & 49 & -13 & 0 & 51 & 0 & 64 & 0 & 38 & 0 & 140 \\
\Xhline{0.5\arrayrulewidth}
2 & 3 & 51 & 153 & -53 & 40 & 43 & 8 & 96 & -32 & -10 & 48 & 76 \\
\Xhline{0.5\arrayrulewidth}
3 & 5 & 53 & 265 & -69 & 16 & 59 & -16 & 128 & -32 & -10 & 0 & 108 \\
\Xhline{0.5\arrayrulewidth}
4 & 7 & 55 & 385 & -61 & -8 & 99 & -40 & 160 & -32 & 38 & -48 & 236 \\
\Xhline{0.5\arrayrulewidth}
5 & 9 & 57 & 513 & -113 & 52 & 63 & 36 & 176 & -16 & -50 & 88 & 76 \\
\Xhline{0.5\arrayrulewidth}
6 & 11 & 59 & 649 & -73 & -40 & 135 & -72 & 208 & -32 & 62 & -112 & 332 \\
\Xhline{0.5\arrayrulewidth}
7 & 13 & 61 & 793 & -117 & 44 & 107 & 28 & 224 & -16 & -10 & 72 & 204 \\
\end{tabular}
\vspace{20em}
\end{table}
\begin{table}[!htbp] 
 \caption{$\Delta_{|p-q|}=22;\ p = 2j + 1,\ j \in \mathbb{N}_0;\ q = p + \Delta$; $dial_P=\{-1,0,6,6\}$}
  \label{table:4k2_od4_zone1}
\begin{tabular}{|c|c|}
\hline 
-1 & 0 \\ 
\Xhline{0.5\arrayrulewidth}
6 & 6 										 						\\
\Xhline{0.5\arrayrulewidth}
\end{tabular}
   $dial\ settings$
\\
\\
\begin{tabular}{|c|c|c|c|c|c|c|c|c|c|c|c|c|c|c|c|}
\hline\hline
$id$ & $p$ & $q$ & $n=pq$ & $od_1$ & $df_1$ & $od_2$ & $df_2$ & $od_3$ & $df_3$ & $od_4$ & $df_4$ & $zone$\\[0.5ex]
\hline
\Xhline{0.5\arrayrulewidth}
1 & 1 & 23 & 23 & -14 & 0 & 58 & 0 & 72 & 0 & 44 & 0 & ... \\
\Xhline{0.5\arrayrulewidth}
2 & 3 & 25 & 75 & -26 & 12 & 94 & -36 & 120 & -48 & 68 & -24 & ... \\
\Xhline{0.5\arrayrulewidth}
3 & 5 & 27 & 135 & -14 & -12 & 154 & -60 & 168 & -48 & 140 & -72 & ... \\
\Xhline{0.5\arrayrulewidth}
4 & 7 & 29 & 203 & -34 & 20 & 158 & -4 & 192 & -24 & 124 & 16 & ... \\
\Xhline{0.5\arrayrulewidth}
5 & 9 & 31 & 279 & -54 & 20 & 162 & -4 & 216 & -24 & 108 & 16 & ... \\
\Xhline{0.5\arrayrulewidth}
6 & 11 & 33 & 363 & -2 & -52 & 262 & -100 & 264 & -48 & \textbf{260} & \textbf{-152} & \textbf{1} \\
\Xhline{0.5\arrayrulewidth}
7 & 13 & 35 & 455 & -14 & 12 & 274 & -12 & 288 & -24 & \textbf{260} & \textbf{0} & \textbf{1} \\
\end{tabular}
\vspace{20em}
\end{table}
\begin{table}[!htbp] 
 \caption{$\Delta_{|p-q|}=46;\ p = 2j + 1,\ j \in \mathbb{N}_0;\ q = p + \Delta$; $dial_P=\{-1,-2,12,12\}$}
  \label{table:od5_changing_dials_1}
\begin{tabular}{|c|c|} 
\hline 
-1 & -2 \\ 
\Xhline{0.5\arrayrulewidth}
12 & 12 \\
\Xhline{0.5\arrayrulewidth}
\end{tabular}
   $dial\ settings$
\\
\\
\begin{tabular}{|c|c|c|c|c|c|c|c|c|c|c|c|c|c|c|c|c|}
\hline\hline
$id$ & $p$ & $q$ & $n=pq$ & $od_1$ & $df_1$ & $od_2$ & $df_2$ & $od_3$ & $df_3$ & $od_4$ & $df_4$ & $od_5$ & $zone$ \\[0.5ex]
\hline 
\Xhline{0.5\arrayrulewidth}
1 & 1 & 47 & 47 & -22 & 0 & 242 & 0 & 264 & 0 & 220 & 0 & 704 & ... \\
\Xhline{0.5\arrayrulewidth}
2 & 3 & 49 & 147 & -26 & 4 & 382 & -140 & 408 & -144 & 356 & -136 & 1120 & ... \\
\Xhline{0.5\arrayrulewidth}
3 & 5 & 51 & 255 & -86 & 60 & 370 & 12 & 456 & -48 & 284 & 72 & 1024 & ... \\
\Xhline{0.5\arrayrulewidth}
4 & 7 & 53 & 371 & -82 & -4 & 470 & -100 & 552 & -96 & 388 & -104 & 1328 & ... \\
\Xhline{0.5\arrayrulewidth}
5 & 9 & 55 & 495 & -54 & -28 & 594 & -124 & 648 & -96 & 540 & -152 & 1728 & ... \\
\Xhline{0.5\arrayrulewidth}
6 & 11 & 57 & 627 & -98 & 44 & 598 & -4 & 696 & -48 & 500 & 40 & 1696 & ... \\
\Xhline{0.5\arrayrulewidth}
7 & 13 & 59 & 767 & -142 & 44 & 602 & -4 & 744 & -48 & 460 & 40 & 1664 & ... \\
\end{tabular}
\end{table}
\begin{table}[!htbp] 
 \caption{$\Delta$ sieve coverage e.g. with $\Delta = 160$}
  \label{table:7}
\begin{tabular}{|c|c|c|c|c|c|} 
\hline 
$od_x$ & id start & id end & $\Delta$ sieve coverage & steady state value & equation reference	    	 \\ 
\Xhline{0.5\arrayrulewidth}
$od_2$ & 79 & 86 & 8 & 6400	 & (\ref{eqn:ssv1b})				 						\\
\Xhline{0.5\arrayrulewidth}
$od_4$ & 163 & 190 & 28 & 12872 & 	(\ref{eqn:ssvgenod4})							 	\\
\Xhline{0.5\arrayrulewidth} 
$od_5$ & 191 & 228 & 38 & 25708 & (\ref{eqn:ssvgenod5})							 	\\
\Xhline{0.5\arrayrulewidth} 
$od_1$ & 761 & 1560 & 800 & 6400 & (\ref{eqn:ssv1b})					 				\\
\Xhline{0.5\arrayrulewidth}
\end{tabular}
\end{table}
\begin{table}[!htbp] 
  \caption{$\Delta$ sieve coverage e.g. with $\Delta = 480$}
  \label{table:9}
\begin{tabular}{|c|c|c|c|c|c|} 
\hline 
$od_x$ & id start & id end & $\Delta$ sieve coverage & steady state value & equation reference 			\\
\Xhline{0.5\arrayrulewidth}
$od_2$ & 913 & 991 & 78 & 57600 & (\ref{eqn:ssv1b})				 						\\
\Xhline{0.5\arrayrulewidth}
$od_4$ & 1683 & 1939 & 256 & 115272  & (\ref{eqn:ssvgenod4})				 				\\
\Xhline{0.5\arrayrulewidth} 
$od_5$ & 1940 & 2281 & 341 & 230508 & (\ref{eqn:ssvgenod5})					 			\\
\Xhline{0.5\arrayrulewidth} 
$od_1$ & 7081 & 14280 & 7199 & 57600 & (\ref{eqn:ssv1b})			 						\\
\Xhline{0.5\arrayrulewidth} 
\end{tabular}
\end{table}
\begin{table}[!htbp] 
  \caption{$\Delta$ sieve zones and coverage with two dials}
  \label{table:10}
\begin{tabular}{|c|c|c|c|c|c|}
\hline
$od_x$ & id start & id end & $\Delta$ sieve coverage & steady state value & equation reference					 \\ 
\Xhline{0.5\arrayrulewidth}
$od_8$ & 14 & 15 & 2 & 2209 & (\ref{eqn:ssv1b})				 									\\
\Xhline{0.5\arrayrulewidth}
$od_8$ & 20 & 22 & 3 & 2209 &  (\ref{eqn:ssv1b})							 						\\
\Xhline{0.5\arrayrulewidth} 
$od_{9a}$ & 30 & 34 & 5 & 4490 & (\ref{eqn:ssvgenod9a})					 						\\
\Xhline{0.5\arrayrulewidth} 
$od_{9b}$ & 35 & 40 & 6 & 4450 & (\ref{eqn:ssvgenod9b})					 						\\
\Xhline{0.5\arrayrulewidth} 
$od_{11}$ & 49 & 57 & 9 & 11189 & (\ref{eqn:ssvgenod11})					 						\\
\Xhline{0.5\arrayrulewidth} 
$od_{10}$ & 58 & 70 & 13 & 6795 & (\ref{eqn:ssvgenod10})					 						\\
\Xhline{0.5\arrayrulewidth} 
$od_{7}$ & 530 & $\infty$ & $530-\infty$ & 2209 & (\ref{eqn:ssv1b})					 				\\
\Xhline{0.5\arrayrulewidth} 
\end{tabular}
\end{table}
\begin{table}[!htbp] 
 \caption{$\Delta_{|p-q|}=22;\ p = 2j + 1,\ j \in \mathbb{N}_0;\ q = p + \Delta$; $dial_P=\{-1,0,6,6\}$}
  \label{table:11}
\begin{tabular}{|c|c|}
\hline 
-1 & 0 \\ 
\Xhline{0.5\arrayrulewidth}
6 & 6  \\
\Xhline{0.5\arrayrulewidth}
\end{tabular}
   $dial\ settings$
\\
\\
\begin{tabular}{|c|c|c|c|c|c|c|c|c|c|c|c|c|c|c|c|}
\hline\hline
$id$ & $p$ & $q$ & $n=pq$ & $od_1$ & $df_1$ & $od_2$ & $df_2$ & $od_3$ & $df_3$ & $od_4$ & $df_4$ & $zone$\\[0.5ex]
\hline
\Xhline{0.5\arrayrulewidth}
1 & 1 & 23 & 23 & -14 & 0 & 58 & 0 & 72 & 0 & 44 & 0 & ... \\
\Xhline{0.5\arrayrulewidth}
2 & 3 & 25 & 75 & -26 & 12 & 94 & -36 & 120 & -48 & 68 & -24 & ... \\
\Xhline{0.5\arrayrulewidth}
3 & 5 & 27 & 135 & -14 & -12 & 154 & -60 & 168 & -48 & 140 & -72 & ... \\
\Xhline{0.5\arrayrulewidth}
4 & 7 & 29 & 203 & -34 & 20 & 158 & -4 & 192 & -24 & 124 & 16 & ... \\
\Xhline{0.5\arrayrulewidth}
5 & 9 & 31 & 279 & -54 & 20 & 162 & -4 & 216 & -24 & 108 & 16 & ... \\
\Xhline{0.5\arrayrulewidth}
\textbf{6} & \textbf{11} & \textbf{33} & \textbf{363} & \textbf{-2} & \textbf{-52} & \textbf{262} & \textbf{-100} & \textbf{264} & \textbf{-48} & \textbf{260} & \textbf{-152} & \textbf{1} \\
\Xhline{0.5\arrayrulewidth}
7 & 13 & 35 & 455 & -14 & 12 & 274 & -12 & 288 & -24 & 260 & 0 & 1 \\
\Xhline{0.5\arrayrulewidth}
8 & 15 & 37 & 555 & -26 & 12 & 286 & -12 & 312 & -24 & 260 & 0 & 1 \\
\Xhline{0.5\arrayrulewidth}
 ... & ... & ... & ... & ... & ... & ... & ... & ...  & ... & ... \\
\Xhline{0.5\arrayrulewidth}
\end{tabular}
\end{table}
\begin{table}[!htbp] 
 \caption{$\Delta_{|p-q|}=22;\ p = 2j + 1,\ j \in \mathbb{N}_0;\ q = p + \Delta$; $dial_P=\{-1,0,4,4\}$}
  \label{table:12}
\begin{tabular}{|c|c|} 
\hline 
-1 & 0 \\ 
\Xhline{0.5\arrayrulewidth}
4 & 4 										 						\\
\Xhline{0.5\arrayrulewidth}
\end{tabular}
   $dial\ settings$
\\
\\
\begin{tabular}{|c|c|c|c|c|c|c|c|c|c|c|}
\hline\hline
$id$ & $p$ & $q$ & $n=pq$ & $od_1$ & $od_2$ & $od_3$ & $od_4$ & $od_5$ & $df_5$ & $zone$\\[0.5ex]
\hline 
\Xhline{0.5\arrayrulewidth}
1 & 1 & 23 & 23 & -14 & 26 & 40 & 12 & 64 & 0 & ... \\
\Xhline{0.5\arrayrulewidth}
2 & 3 & 25 & 75 & -26 & 46 & 72 & 20 & 112 & -48 & ... \\
\Xhline{0.5\arrayrulewidth}
3 & 5 & 27 & 135 & -14 & 90 & 104 & 76 & 256 & -144 & ... \\
\Xhline{0.5\arrayrulewidth}
4 & 7 & 29 & 203 & -34 & 86 & 120 & 52 & 224 & 32 & ... \\
\Xhline{0.5\arrayrulewidth}
5 & 9 & 31 & 279 & -54 & 82 & 136 & 28 & 192 & 32 & ... \\
\Xhline{0.5\arrayrulewidth}
\textbf{6} & \textbf{11} & \textbf{33} & \textbf{363} & \textbf{-2} & \textbf{166} & \textbf{168} & \textbf{164} & \textbf{496} & \textbf{-304} & \textbf{0} \\
\Xhline{0.5\arrayrulewidth}
7 & 13 & 35 & 455 & -14 & 170 & 184 & 156 & 496 & 0 & 0 \\
\Xhline{0.5\arrayrulewidth}
 ... & ... & ... & ... & ... & ... & ... & ... & ...  & ... & ... \\
\Xhline{0.5\arrayrulewidth}
\end{tabular}
\end{table}
\begin{table}[!htbp] 
 \caption{$\Delta_{|p-q|}=22;\ p = 2j + 1,\ j \in \mathbb{N}_0;\ q = p + \Delta$; $dial_P=\{-1,0,6,6\}$}
  \label{table:13}
\begin{tabular}{|c|c|}
\hline 
-1 & 0 \\ 
\Xhline{0.5\arrayrulewidth}
6 & 6  \\
\Xhline{0.5\arrayrulewidth}
\end{tabular}
   $dial\ settings$
\\
\\
\begin{tabular}{|c|c|c|c|c|c|c|c|c|c|c|c|c|c|c|c|}
\hline\hline
$id$ & $p$ & $q$ & $n=pq$ & $od_1$ & $df_1$ & $od_2$ & $df_2$ & $od_3$ & $df_3$ & $od_4$ & $df_4$ & $zone$\\[0.5ex]
\hline
\Xhline{0.5\arrayrulewidth}
1 & 1 & 23 & 23 & -14 & 0 & 58 & 0 & 72 & 0 & 44 & 0 & ... \\
\Xhline{0.5\arrayrulewidth}
\textbf{2} & \textbf{3} & \textbf{25} & \textbf{75} & \textbf{-26} & \textbf{12} & \textbf{94} & \textbf{-36} & \textbf{120} & \textbf{-48} & \textbf{68} & \textbf{-24} & \textbf{...} \\
\Xhline{0.5\arrayrulewidth}
\textbf{3} & \textbf{5} & \textbf{27} & \textbf{135} & \textbf{-14} & \textbf{-12} & \textbf{154} & \textbf{-60} & \textbf{168} & \textbf{-48} & \textbf{140} & \textbf{-72} & \textbf{...} \\
\Xhline{0.5\arrayrulewidth}
4 & 7 & 29 & 203 & -34 & 20 & 158 & -4 & 192 & -24 & 124 & 16 & ... \\
\Xhline{0.5\arrayrulewidth}
5 & 9 & 31 & 279 & -54 & 20 & 162 & -4 & 216 & -24 & 108 & 16 & ... \\
\Xhline{0.5\arrayrulewidth}
\textbf{6} & \textbf{11} & \textbf{33} & \textbf{363} & \textbf{-2} & \textbf{-52} & \textbf{262} & \textbf{-100} & \textbf{264} & \textbf{-48} & \textbf{260} & \textbf{-152} & \textbf{1} \\
\Xhline{0.5\arrayrulewidth}
7 & 13 & 35 & 455 & -14 & 12 & 274 & -12 & 288 & -24 & 260 & 0 & 1 \\
\Xhline{0.5\arrayrulewidth}
8 & 15 & 37 & 555 & -26 & 12 & 286 & -12 & 312 & -24 & 260 & 0 & 1 \\
\Xhline{0.5\arrayrulewidth}
 ... & ... & ... & ... & ... & ... & ... & ... & ...  & ... & ... \\
\Xhline{0.5\arrayrulewidth}
\end{tabular}
\end{table}
\begin{table}[!htbp] 
 \caption{\begin{math}\Delta_{|p-q|}=540\end{math}}
  \label{table:14}
\begin{tabular}{|c|c|}
\hline 
0 & -1 \\ 
\Xhline{0.5\arrayrulewidth}
2 & 2 										 						\\
\Xhline{0.5\arrayrulewidth} 
\end{tabular}
   $dial\ settings$
\\
\begin{tabular}{|c|c|c|c|c|c|c|c|c|c|c|c|c|}
\hline\hline
$id$ & $p$ & $q$ & $n=pq$ & $d_1$ & $od_1$ & $df_1$ & $od_2$ & $df_2$ & $od_3$ & $df_3$ & $od_4$ & $df_4$  \\[0.5ex]
\hline
1 & 1 & 541 & 541 & 22 & -57 & 0 & 35 & 0 & 92 & 0 & -22 & 0  \\
\Xhline{0.5\arrayrulewidth}
... & ... & ... & ... & ... & ... & ... & ... & ... & ... & ... & ... & ...  \\
\Xhline{0.5\arrayrulewidth}
136 & 271 & 811 & \textbf{219781} & \textbf{468} & \textbf{-757} & 292 & \textbf{1119} & 284 & \textbf{1876} & -8 & \textbf{362} & 576  \\
\Xhline{0.5\arrayrulewidth}
137 & 273 & 813 & 221949 & 470 & -1049 & 292 & 835 & 284 & 1884 & -8 & -214 & 576  \\
\Xhline{0.5\arrayrulewidth}
... & ... & ... & ... & ... & ... & ... & ... & ... & ... & ... & ... & ...  \\
\Xhline{0.5\arrayrulewidth}
\end{tabular}
\end{table}
\begin{table}[!htbp] 
 \caption{\begin{math}\Delta_{|p-q|}=22\end{math}}
  \label{table:range15}
\begin{tabular}{|c|c|}
\hline 
-1 & 0 \\ 
\Xhline{0.5\arrayrulewidth}
2 & 2 										 						\\
\Xhline{0.5\arrayrulewidth} 
\end{tabular}
   $dial\ settings$
\\
\begin{tabular}{|c|c|c|c|c|c|c|c|c|c|c|c|c|c|c|}
\hline\hline
$id$ & $p$ & $q$ & $n=pq$ & $\lfloor \sqrt{n} \rfloor$ & $d_1$ & $d_2$ & $d_1^2$ & $d_2^2$ & $od_1$ & $od_2$ & $od_3$ & $od_4$ & $od_6$ \\[0.5ex]
\hline
1 & 1 & 23 & 23 & 4 & 3 & 5 & 9 & 25 & -14 & 2 & 16 & -12 & 8 \\
\Xhline{0.5\arrayrulewidth}
... & ... & ... & ... & ... & ... & ... & ... & ... & ... & ... & ... & ... & ... \\
\Xhline{0.5\arrayrulewidth}
13 & 25 & 47 & 1175 & 34 & 33 & 35 & 1089 & 1225 & -86 & 50 & 136 & -36 & 20 \\
\Xhline{0.5\arrayrulewidth}
14 & 27 & 49 & 1323 & 36 & \textbf{35} & 37 & 1225 & 1369 & -98 & 46 & 144 & -52 & 28 \\
\Xhline{0.5\arrayrulewidth}
15 & 29 & 51 & 1479 & 38 & \textbf{37} & 39 & 1369 & 1521 & -110 & 42 & 152 & -68 & 36 \\
\Xhline{0.5\arrayrulewidth}
16 & 31 & 53 & \textbf{1643} & 40 & \textbf{39} & 41 & 1521 & 1681 & \textbf{-122} & 38 & 160 & -84 & 44 \\
\Xhline{0.5\arrayrulewidth}
17 & 33 & 55 & 1815 & 42 & \textbf{41} & 43 & 1681 & 1849 & -134 & 34 & 168 & -100 & 52 \\
\Xhline{0.5\arrayrulewidth}
18 & 35 & 57 & 1995 & 44 & \textbf{43} & 45 & 1849 & 2025 & -146 & 30 & 176 & -116 & 60 \\
\Xhline{0.5\arrayrulewidth}
19 & 37 & 59 & 2183 & 46 & 45 & 47 & 2025 & 2209 & -158 & 26 & 184 & -132 & 68 \\
\Xhline{0.5\arrayrulewidth}
20 & 39 & 61 & 2379 & 48 & 47 & 49 & 2209 & 2401 & -170 & 22 & 192 & -148 & 76 \\
\Xhline{0.5\arrayrulewidth}
\begin{math}\infty\end{math} & \begin{math}\infty\end{math}  &  \begin{math}\infty\end{math} &  \begin{math}\infty\end{math}  & ... & ... & ... & ... & ...  & ... & ... & ... & 244  & 120 \\ [1ex]
\Xhline{0.5\arrayrulewidth}
\end{tabular}
\end{table}
\begin{table}[!htbp] 
 \caption{\begin{math}\Delta_{|p-q|}=22\end{math}}
  \label{table:roxy16}
\begin{tabular}{|c|c|}
\hline 
-1 & 0 \\ 
\Xhline{0.5\arrayrulewidth}
2 & 2 										 						\\
\Xhline{0.5\arrayrulewidth} 
\end{tabular}
   $dial\ settings$
\\
\begin{tabular}{|c|c|c|c|c|c|c|c|c|c|c|c|c|c|c|}
\hline\hline
$id$ & $p$ & $q$ & $n=pq$ & $\lfloor \sqrt{n} \rfloor$ & $d_1$ & $d_2$ & $d_1^2$ & $d_2^2$ & $od_1$ & $od_2$ & $od_3$ & $od_4$ & $od_6$ \\[0.5ex]
\hline
1 & 1 & 23 & 23 & 4 & 3 & 5 & 9 & 25 & -14 & 2 & 16 & -12 & 8 \\
\Xhline{0.5\arrayrulewidth}
2 & 3 & 25 & 75 & 8 & 7 & 9 & 49 & 81 & -26 & 6 & 32 & -20 & 12 \\
\Xhline{0.5\arrayrulewidth}
3 & 5 & 27 & 135 & 11 & 11 & 13 & 121 & 169 & -14 & 34 & 48 & 20 & 8 \\
\Xhline{0.5\arrayrulewidth}
4 & 7 & 29 & 203 & 14 & 13 & 15 & 169 & 225 & -34 & 22 & 56 & -12 & 8 \\
\Xhline{0.5\arrayrulewidth}
5 & 9 & 31 & 279 & 16 & 15 & 17 & 225 & 289 & -54 & 10 & 64 & -44 & 24 \\
\Xhline{0.5\arrayrulewidth}
6 & 11 & 33 & 363 & 19 & 19 & 21 & 361 & 441 & -2 & 78 & 80 & 76 & \color{blue}\textbf{36} \\
\Xhline{0.5\arrayrulewidth}
7 & 13 & 35 & 455 & 21 & 21 & 23 & 441 & 529 & -14 & 74 & 88 & 60 & \textbf{28} \\
\end{tabular}
\vspace{30em}
\end{table}
\begin{table}[!htbp] 
 \caption{\begin{math}\Delta_{|p-q|}=20\end{math}}
  \label{table:17}
\begin{tabular}{|c|c|}
\hline 
0 & -1 \\ 
\Xhline{0.5\arrayrulewidth}
2 & 2 										 						\\
\Xhline{0.5\arrayrulewidth} 
\end{tabular}
   $dial\ settings$
\\
\begin{tabular}{|c|c|c|c|c|c|c|c|c|c|c|c|c|c|c|}
\hline\hline
$id$ & $p$ & $q$ & $n=pq$ & $\lfloor \sqrt{n} \rfloor$ & $d_1$ & $d_2$ & $d_1^2$ & $d_2^2$ & $od_1$ & $od_2$ & $od_3$ & $od_4$ & $od_6$ \\[0.5ex]
\hline
1 & 1 & 21 & 21 & 4 & 4 & 6 & 16 & 36 & -5 & 15 & 20 & 10 & 3 \\
\Xhline{0.5\arrayrulewidth}
2 & 3 & 23 & 69 & 8 & 8 & 10 & 64 & 100 & -5 & 31 & 36 & 26 & 11 \\
\Xhline{0.5\arrayrulewidth}
3 & 5 & 25 & 125 & 11 & 10 & 12 & 100 & 144 & -25 & 19 & 44 & -6 & 5 \\
\Xhline{0.5\arrayrulewidth}
4 & 7 & 27 & 189 & 13 & 12 & 14 & 144 & 196 & -45 & 7 & 52 & -38 & 21 \\
\Xhline{0.5\arrayrulewidth}
5 & 9 & 29 & 261 & 16 & 16 & 18 & 256 & 324 & -5 & 63 & 68 & 58 & \color{black}\textbf{27} \\
\Xhline{0.5\arrayrulewidth}
6 & 11 & 31 & 341 & 18 & 18 & 20 & 324 & 400 & -17 & 59 & 76 & 42 & \color{blue}\textbf{19} \\
\Xhline{0.5\arrayrulewidth}
7 & 13 & 33 & 429 & 20 & 20 & 22 & 400 & 484 & -29 & 55 & 84 & 26 & \color{black}\textbf{11} \\
\end{tabular}
\end{table}
\begin{table}[!htbp] 
 \caption{$r\Delta$ example}
  \label{table:23}
\begin{tabular}{|c|c|c|c|c|}
\hline 
known $\Delta$ & id & $dial_P$ & $\Delta$ sieve value & $od_x$\\ 
\Xhline{0.5\arrayrulewidth}
122 & 93 & \{0, -1, 8, 8\} & 3721 & $od_2$									 						\\
\Xhline{0.5\arrayrulewidth} 
\end{tabular}
\\
\\
\begin{tabular}{|c|c|c|c|c|}
\hline
Unknown $\Delta$ & $n$ & $dial_P$ & $\Delta$ sieve value & $od_x$ \\[0.5ex]
\hline
$162$ & $n1$ &   \{0, -1, 8, 8\} & $13154$ & $od_4$ \\
\Xhline{0.5\arrayrulewidth}
$178$ & $n2$ &   \{0, -1, 8, 8\} & $31732$ & $od_5$  \\
\Xhline{0.5\arrayrulewidth}
$202$ & $n3$ &   \{0, -1, 8, 8\} & $10201$ & $od_2$  \\
\Xhline{0.5\arrayrulewidth}
\end{tabular}
\end{table}
\begin{table}[!htbp] 
 \caption{$\sum_{p+q}=22;\ p = 2j+ 1$ for $j \in \mathbb{N}_0;\ q=\sum-p;\ dial_P = \{-1,0,2,2\}$}
  \label{table:sum_series_1}
\begin{tabular}{|c|c|}
\hline 
-1 & 0 \\ 
\Xhline{0.5\arrayrulewidth}
2 & 2 										 						\\
\Xhline{0.5\arrayrulewidth} 
\end{tabular}
   $dial\ settings$
\\
\begin{tabular}{|c|c|c|c|c|c|c|c|c|c|c|c|c|c|c|}
\hline\hline
$id$ & $p$ & $q$ & $N=pq$ & $\lfloor \sqrt{n} \rfloor$ & $d_1$ & $d_2$ & $d_1^2$ & $d_2^2$ & $od_1$ & $od_2$ & $od_3$ & $od_4$ & $od_6$ \\[0.5ex]
\hline
1 & 1 & 21 & 21 & 4 & 3 & 5 & 9 & 25 & -12 & 4 & 16 & -8 & 6 \\
\Xhline{0.5\arrayrulewidth}
2 & 3 & 19 & 57 & 7 & 7 & 9 & 49 & 81 & -8 & 24 & 32 & 16 & 6 \\
\Xhline{0.5\arrayrulewidth}
3 & 5 & 17 & 85 & 9 & 9 & 11 & 81 & 121 & -4 & 36 & 40 & 32 & 14 \\
\Xhline{0.5\arrayrulewidth}
4 & 7 & 15 & 105 & 10 & 9 & 11 & 81 & 121 & -24 & 16 & 40 & -8 & 6 \\
\Xhline{0.5\arrayrulewidth}
5 & 9 & 13 & 117 & 10 & 9 & 11 & 81 & 121 & -36 & 4 & 40 & -32 & 18 \\
\Xhline{0.5\arrayrulewidth}
6 & 11 & 11 & 121 & 11 & 11 & 13 & 121 & 169 & 0 & 48 & 48 & 48 & 22 \\
\Xhline{0.5\arrayrulewidth}
7 & 13 & 9 & 117 & 10 & 9 & 11 & 81 & 121 & -36 & 4 & 40 & -32 & 18 \\
\Xhline{0.5\arrayrulewidth}
8 & 15 & 7 & 105 & 10 & 9 & 11 & 81 & 121 & -24 & 16 & 40 & -8 & 6 \\
\Xhline{0.5\arrayrulewidth}
9 & 17 & 5 & 85 & 9 & 9 & 11 & 81 & 121 & -4 & 36 & 40 & 32 & 14 \\
\Xhline{0.5\arrayrulewidth}
10 & 19 & 3 & 57 & 7 & 7 & 9 & 49 & 81 & -8 & 24 & 32 & 16 & 6 \\
\Xhline{0.5\arrayrulewidth}
11 & 21 & 1 & 21 & 4 & 3 & 5 & 9 & 25 & -12 & 4 & 16 & -8 & 6 \\
\Xhline{0.5\arrayrulewidth}
\end{tabular}
\end{table}
\begin{table}[!htbp] 
 \caption{$\sum_{p+q}=20;\ p = 2j+ 1$ for $j \in \mathbb{N}_0;\ q=\sum-p;\ dial_P = \{-1,0,2,2\}$}
  \label{table:sum_series_2}
\begin{tabular}{|c|c|}
\hline 
-1 & 0 \\ 
\Xhline{0.5\arrayrulewidth}
2 & 2 										 						\\
\Xhline{0.5\arrayrulewidth} 
\end{tabular}
   $dial\ settings$
\\
\begin{tabular}{|c|c|c|c|c|c|c|c|c|c|c|c|c|c|c|}
\hline\hline
$id$ & $p$ & $q$ & $N=pq$ & $\lfloor \sqrt{n} \rfloor$ & $d_1$ & $d_2$ & $d_1^2$ & $d_2^2$ & $od_1$ & $od_2$ & $od_3$ & $od_4$ & $od_6$ \\[0.5ex]
\hline
1 & 1 & 19 & 19 & 4 & 3 & 5 & 9 & 25 & -10 & 6 & 16 & -4 & 4 \\
\Xhline{0.5\arrayrulewidth}
2 & 3 & 17 & 51 & 7 & 7 & 9 & 49 & 81 & -2 & 30 & 32 & 28 & 12 \\
\Xhline{0.5\arrayrulewidth}
3 & 5 & 15 & 75 & 8 & 7 & 9 & 49 & 81 & -26 & 6 & 32 & -20 & 12 \\
\Xhline{0.5\arrayrulewidth}
4 & 7 & 13 & 91 & 9 & 9 & 11 & 81 & 121 & -10 & 30 & 40 & 20 & 8 \\
\Xhline{0.5\arrayrulewidth}
5 & 9 & 11 & 99 & 9 & 9 & 11 & 81 & 121 & -18 & 22 & 40 & 4 & 0 \\
\Xhline{0.5\arrayrulewidth}
6 & 11 & 9 & 99 & 9 & 9 & 11 & 81 & 121 & -18 & 22 & 40 & 4 & 0 \\
\Xhline{0.5\arrayrulewidth}
7 & 13 & 7 & 91 & 9 & 9 & 11 & 81 & 121 & -10 & 30 & 40 & 20 & 8 \\
\Xhline{0.5\arrayrulewidth}
8 & 15 & 5 & 75 & 8 & 7 & 9 & 49 & 81 & -26 & 6 & 32 & -20 & 12 \\
\Xhline{0.5\arrayrulewidth}
9 & 17 & 3 & 51 & 7 & 7 & 9 & 49 & 81 & -2 & 30 & 32 & 28 & 12 \\
\Xhline{0.5\arrayrulewidth}
10 & 19 & 1 & 19 & 4 & 3 & 5 & 9 & 25 & -10 & 6 & 16 & -4 & 4 \\
\Xhline{0.5\arrayrulewidth}
\end{tabular}
\end{table}
\begin{table}[!htbp] 
 \caption{$\Delta_{|p-q|}=20, p = 2j+1$ for $j \in \mathbb{N}$}
  \label{table:delta_equilibrim_1}
\begin{tabular}{|c|c|}
\hline 
0 & -1 \\ 
\Xhline{0.5\arrayrulewidth}
2 & 2 										 						\\
\Xhline{0.5\arrayrulewidth} 
\end{tabular}
   $dial\ settings$
\\
\begin{tabular}{|c|c|c|c|c|c|c|c|c|c|c|c|c|c|c|}
\hline\hline
$id$ & $p$ & $q$ & $n=pq$ & $\lfloor \sqrt{n} \rfloor$ & $d_1$ & $d_2$ & $d_1^2$ & $d_2^2$ & $od_1$ & $od_2$ & $od_3$ & $od_4$ & $od_6$ \\[0.5ex]
\hline
1 & 1 & 21 & 21 & 4 & 4 & 6 & 16 & 36 & -5 & 15 & 20 & 10 & 3 \\
\Xhline{0.5\arrayrulewidth}
... & ... & ... & ... & ... & ... & ... & ... & ... & ... & ... & ... & ... & ... \\
\Xhline{0.5\arrayrulewidth}
17 & 33 & 53 & 1749 & 41 & 40 & 42 & 1600 & 1764 & -149 & 15 & 164 & -134 & 69 \\
\Xhline{0.5\arrayrulewidth}
18 & 35 & 55 & 1925 & 43 & 42 & 44 & 1764 & 1936 & -161 & 11 & 172 & -150 &77 \\
\Xhline{0.5\arrayrulewidth}
19 & 37 & 57 & 2109 & 45 & 44 & 46 & 1936 & 2116 & -173 & 7 & 180 & -166 & 85 \\
\Xhline{0.5\arrayrulewidth}
20 & 39 & 59 & 2301 & 47 & 46 & 48 & 2116 & 2304 & -185 & 3 & 188 & -182 & 93 \\
\Xhline{0.5\arrayrulewidth}
21 & 41 & 61 & 2501 & 50 & 50 & 52 & 2500 & 2704 & -1 & 203 & 204 & 202 & \color{blue}\textbf{99} \\
\Xhline{0.5\arrayrulewidth}
22 & 43 & 63 & 2709 & 52 & 52 & 54 & 2704 & 2916 & -5 & 207 & 212 & 202 & \color{blue}\textbf{99} \\
\Xhline{0.5\arrayrulewidth}
... & ... & ... & ... & ... & ... & ... & ... & ... & ... & ... & ... & 202 & \color{blue}\textbf{99} \\
\Xhline{0.5\arrayrulewidth}
$\infty$ & $\infty$ & $\infty$ & $\infty$ & $\infty$ & $\infty$ & $\infty$ & $\infty$ & $\infty$ & $\infty$ & $\infty$ & $\infty$ & 202 & \color{blue}\textbf{99} \\
\Xhline{0.5\arrayrulewidth}
\end{tabular}
\end{table}
\begin{table}[!htbp] 
 \caption{$\sum_{p+q}=20, p = 2j+1$ for $j \in \mathbb{N}$}
  \label{table:sum_equilibrim_1}
\begin{tabular}{|c|c|}
\hline 
-1 & 0 \\ 
\Xhline{0.5\arrayrulewidth}
2 & 2 										 						\\
\Xhline{0.5\arrayrulewidth} 
\end{tabular}
   $dial\ settings$
\\
\begin{tabular}{|c|c|c|c|c|c|c|c|c|c|c|c|c|c|c|}
\hline\hline
$id$ & $p$ & $q$ & $N=pq$ & $\lfloor \sqrt{n} \rfloor$ & $d_1$ & $d_2$ & $d_1^2$ & $d_2^2$ & $od_1$ & $od_2$ & $od_3$ & $od_4$ & $od_6$ \\[0.5ex]
\hline
1 & 1 & 19 & 19 & 4 & 3 & 5 & 9 & 25 & -10 & 6 & 16 & -4 & 4 \\
\Xhline{0.5\arrayrulewidth}
2 & 3 & 17 & 51 & 7 & 7 & 9 & 49 & 81 & -2 & 30 & 32 & 28 & 12 \\
\Xhline{0.5\arrayrulewidth}
3 & 5 & 15 & 75 & 8 & 7 & 9 & 49 & 81 & -26 & 6 & 32 & -20 & 12 \\
\Xhline{0.5\arrayrulewidth}
4 & 7 & 13 & 91 & 9 & 9 & 11 & 81 & 121 & -10 & 30 & 40 & 20 & 8 \\
\Xhline{0.5\arrayrulewidth}
5 & 9 & 11 & \color{blue}\textbf{99} & 9 & 9 & 11 & 81 & 121 & -18 & 22 & 40 & 4 & \color{red}\textbf{0} \\
\Xhline{0.5\arrayrulewidth}
6 & 11 & 9 & \color{blue}\textbf{99} & 9 & 9 & 11 & 81 & 121 & -18 & 22 & 40 & 4 & \color{red}\textbf{0} \\
\Xhline{0.5\arrayrulewidth}
7 & 13 & 7 & 91 & 9 & 9 & 11 & 81 & 121 & -10 & 30 & 40 & 20 & 8 \\
\Xhline{0.5\arrayrulewidth}
8 & 15 & 5 & 75 & 8 & 7 & 9 & 49 & 81 & -26 & 6 & 32 & -20 & 12 \\
\Xhline{0.5\arrayrulewidth}
9 & 17 & 3 & 51 & 7 & 7 & 9 & 49 & 81 & -2 & 30 & 32 & 28 & 12 \\
\Xhline{0.5\arrayrulewidth}
10 & 19 & 1 & 19 & 4 & 3 & 5 & 9 & 25 & -10 & 6 & 16 & -4 & 4 \\
\Xhline{0.5\arrayrulewidth}
... & ... & ... & ... & ... & ... & ... & ... & ... & ... & ... & ... & ... & ... \\
\Xhline{0.5\arrayrulewidth}
\end{tabular}
\end{table}
\begin{table}[!htbp] 
\parbox{.45\linewidth}{
\begin{tabular}{|c|c|}
\hline 
$\Delta$ & $\sum$ \\ 
\Xhline{0.5\arrayrulewidth}
\{0,-1,2,2\} & \{-1,0,2,2\} \\
\Xhline{0.5\arrayrulewidth} 
\end{tabular}
   $dial\ settings$ \\ \\
\begin{tabular}{|c|c|c|c|c|c|c|c|c|c|c|c|c|c|c|}
\hline
$N_{_{\sum_{p+q}}}$ & $od_6{_{_{_{\Delta_{|p-q|}}}}}$  & $\Delta \rightleftharpoons \sum$ & $df_{_{\Delta\sum}}$\\[0.5ex]
\hline
99 & 99& 0 &  \\
\Xhline{0.5\arrayrulewidth}
91 & 93 & -2 &  2 \\
\Xhline{0.5\arrayrulewidth}
75 & 85 & -10 &  8 \\
\Xhline{0.5\arrayrulewidth}
51 & 77 & -26 &  16 \\
\Xhline{0.5\arrayrulewidth}
19 & 69 & -50 &  24 \\
\Xhline{0.5\arrayrulewidth}
-21 & 61 & -82 &  32 \\
\Xhline{0.5\arrayrulewidth}
-69 & 53 & -122 &  40 \\
\Xhline{0.5\arrayrulewidth}
-125 & 45 & -170 &  48 \\
\Xhline{0.5\arrayrulewidth}
-189 & 37 & -226 &  56 \\
\Xhline{0.5\arrayrulewidth}
-261 & 29 & -290 &  64 \\
\Xhline{0.5\arrayrulewidth}
\end{tabular}
\captionsetup{justification=centering}
{
\caption{\begin{math}\Delta_{|p-q|} = \sum_{p+q}=20\end{math}}
}
}
 \hfill
\parbox{.45\linewidth}{
\begin{tabular}{|c|c|} 
\hline 
 $\Delta$ & $\sum$ \\ 
\Xhline{0.5\arrayrulewidth}
\{0,-1,2,2\} & \{-1,0,2,2\} \\
\Xhline{0.5\arrayrulewidth} 
\end{tabular}
   $dial\ settings$ \\ \\
\begin{tabular}{|c|c|c|c|c|c|c|c|c|c|c|c|c|c|c|}
\hline
$N_{_{\sum_{p+q}}}$ & $od_6{_{_{_{\Delta_{|p-q|}}}}}$  & $\Delta \rightleftharpoons \sum$ & $df_{_{\Delta\sum}}$\\[0.5ex]
\hline
399 & 399 & 0 &  \\
\Xhline{0.5\arrayrulewidth}
391 & 393 & -2 &  2 \\
\Xhline{0.5\arrayrulewidth}
375 & 385 & -10 &  8 \\
\Xhline{0.5\arrayrulewidth}
351 & 377 & -26 &  16 \\
\Xhline{0.5\arrayrulewidth}
319 & 369 & -50 &  24 \\
\Xhline{0.5\arrayrulewidth}
279 & 361 & -82 &  32 \\
\Xhline{0.5\arrayrulewidth}
231 & 353 & -122 &  40 \\
\Xhline{0.5\arrayrulewidth}
175 & 345 & -170 &  48 \\
\Xhline{0.5\arrayrulewidth}
111 & 337 & -226 &  56 \\
\Xhline{0.5\arrayrulewidth}
39 & 329 & -290 &  64 \\
\Xhline{0.5\arrayrulewidth}
\end{tabular}
\captionsetup{justification=centering}
{
 \caption{\begin{math}\Delta_{|p-q|} = \sum_{p+q}=40\end{math}} \
}
}
\end{table}
\begin{table}[!htbp] 
 \caption{$\Delta_{|p=q|}=22, p = 2j+1$ for $j \in \mathbb{N}$}
  \label{table:delta_4k2_1}
\begin{tabular}{|c|c|}
\hline 
-1 & 0 \\ 
\Xhline{0.5\arrayrulewidth}
2 & 2 										 						\\
\Xhline{0.5\arrayrulewidth} 
\end{tabular}
   $dial\ settings$
\\
\begin{tabular}{|c|c|c|c|c|c|c|c|c|c|c|c|c|c|c|}
\hline\hline
$id$ & $p$ & $q$ & $n=pq$ & $\lfloor \sqrt{n} \rfloor$ & $d_1$ & $d_2$ & $d_1^2$ & $d_2^2$ & $od_1$ & $od_2$ & $od_3$ & $od_4$ & $od_6$ \\[0.5ex]
\hline
1 & 1 & 23 & 23 & 4 & 3 & 5 & 9 & 25 & -14 & 2 & 16 & -12 & 8 \\
\Xhline{0.5\arrayrulewidth}
... & ... & ... & ... & ... & ... & ... & ... & ... & ... & ... & ... & ... & ... \\
\Xhline{0.5\arrayrulewidth}
21 & 41 & 63 & 2583 & 50 & 49 & 51 & 2401 & 2601 & -182 & 18 & 200 & -164 & 84 \\
\Xhline{0.5\arrayrulewidth}
22 & 43 & 65 & 2795 & 52 & 51 & 53 & 2601 & 2809 & -194 & 14 & 208 & -180 & 92 \\
\Xhline{0.5\arrayrulewidth}
23 & 45 & 67 & 3015 & 54 & 53 & 55 & 2809 & 3025 & -206 & 10 & 216 & -196 & 100 \\
\Xhline{0.5\arrayrulewidth}
24 & 47 & 69 & 3243 & 56 & 55 & 57 & 3025 & 3249 & -218 & 6 & 224 & -212 & 108 \\
\Xhline{0.5\arrayrulewidth}
25 & 49 & 71 & 3479 & 58 & 57 & 59 & 3249 & 3481 & -230 & 2 & 232 & -228 & 116 \\
\Xhline{0.5\arrayrulewidth}
26 & 51 & 73 & 3723 & 61 & 61 & 63 & 3721 & 3969 & -2 & 246 & 248 & 244 & \color{blue}\textbf{120} \\
\Xhline{0.5\arrayrulewidth}
27 & 53 & 75 & 3975 & 63 & 63 & 65 & 3969 & 4225 & -6 & 250 & 256 & 244 & \color{blue}\textbf{120} \\
\Xhline{0.5\arrayrulewidth}
 ... & ... & ... & ... & ... & ... & ... & ... & ...  & ... & ... & ... & 244  & \color{blue}\textbf{120}\\
\Xhline{0.5\arrayrulewidth}
\begin{math}\infty\end{math} & \begin{math}\infty\end{math}  &  \begin{math}\infty\end{math} &  \begin{math}\infty\end{math}  & ... & ... & ... & ... & ...  & ... & ... & ... & 244  & \color{blue}\textbf{120} \\ [1ex]
\Xhline{0.5\arrayrulewidth}
\end{tabular}
\\
\end{table}
\begin{table}[!htbp] 
 \caption{$\sum_{p+q}=22, p = 2j+1$ for $j \in \mathbb{N}$}
  \label{table:sum_4k2_1}
\begin{tabular}{|c|c|}
\hline 
-1 & 0 \\ 
\Xhline{0.5\arrayrulewidth}
2 & 2 										 						\\
\Xhline{0.5\arrayrulewidth} 
\end{tabular}
   $dial\ settings$
\\
\begin{tabular}{|c|c|c|c|c|c|c|c|c|c|c|c|c|c|c|}
\hline\hline
$id$ & $p$ & $q$ & $N=pq$ & $\lfloor \sqrt{n} \rfloor$ & $d_1$ & $d_2$ & $d_1^2$ & $d_2^2$ & $od_1$ & $od_2$ & $od_3$ & $od_4$ & $od_6$ \\[0.5ex]
\hline
1 & 1 & 21 & 21 & 4 & 3 & 5 & 9 & 25 & -12 & 4 & 16 & -8 & 6 \\
\Xhline{0.5\arrayrulewidth}
2 & 3 & 19 & 57 & 7 & 7 & 9 & 49 & 81 & -8 & 24 & 32 & 16 & 6 \\
\Xhline{0.5\arrayrulewidth}
3 & 5 & 17 & 85 & 9 & 9 & 11 & 81 & 121 & -4 & 36 & 40 & 32 & 14 \\
\Xhline{0.5\arrayrulewidth}
4 & 7 & 15 & 105 & 10 & 9 & 11 & 81 & 121 & -24 & 16 & 40 & -8 & 6 \\
\Xhline{0.5\arrayrulewidth}
5 & 9 & 13 & 117 & 10 & 9 & 11 & 81 & 121 & -36 & 4 & 40 & -32 & 18 \\
\Xhline{0.5\arrayrulewidth}
6 & 11 & 11 & \color{blue}\textbf{121} & 11 & 11 & 13 & 121 & 169 & \color{red}\textbf{0} & 48 & 48 & 48 & 22 \\
\Xhline{0.5\arrayrulewidth}
7 & 13 & 9 & 117 & 10 & 9 & 11 & 81 & 121 & -36 & 4 & 40 & -32 & 18 \\
\Xhline{0.5\arrayrulewidth}
8 & 15 & 7 & 105 & 10 & 9 & 11 & 81 & 121 & -24 & 16 & 40 & -8 & 6 \\
\Xhline{0.5\arrayrulewidth}
9 & 17 & 5 & 85 & 9 & 9 & 11 & 81 & 121 & -4 & 36 & 40 & 32 & 14 \\
\Xhline{0.5\arrayrulewidth}
10 & 19 & 3 & 57 & 7 & 7 & 9 & 49 & 81 & -8 & 24 & 32 & 16 & 6 \\
\Xhline{0.5\arrayrulewidth}
11 & 21 & 1 & 21 & 4 & 3 & 5 & 9 & 25 & -12 & 4 & 16 & -8 & 6 \\
\Xhline{0.5\arrayrulewidth}
\end{tabular}
\end{table}
\begin{table}[!htbp] 
\parbox{.45\linewidth}{
\begin{tabular}{|c|c|}
\hline 
$\Delta$ & $\sum$ \\ 
\Xhline{0.5\arrayrulewidth}
\{-1,0,2,2\} & \{-1,0,2,2\} \\
\Xhline{0.5\arrayrulewidth} 
\end{tabular}
   $dial\ settings$ \\ \\
\begin{tabular}{|c|c|c|c|c|c|c|c|c|c|c|c|c|c|c|}
\hline
$N_{_{\sum_{p+q}}}$ & $od_6{_{_{_{\Delta_{|p-q|}}}}}$  & $\Delta \rightleftharpoons \sum$ & $df_{_{\Delta\sum}}$\\[0.5ex]
\hline
121 & 120 & 1 &  \\
\Xhline{0.5\arrayrulewidth}
117 & 116 & 1 &  0 \\
\Xhline{0.5\arrayrulewidth}
105 & 108 & -3 &  4 \\
\Xhline{0.5\arrayrulewidth}
85 & 100 & -15 &  12 \\
\Xhline{0.5\arrayrulewidth}
57 & 92 & -35 &  20 \\
\Xhline{0.5\arrayrulewidth}
21 & 84 & -63 &  28 \\
\Xhline{0.5\arrayrulewidth}
-23 & 76 & -99 &  36 \\
\Xhline{0.5\arrayrulewidth}
-75 & 68 & -143 &  44 \\
\Xhline{0.5\arrayrulewidth}
-135 & 60 & -195 &  52 \\
\Xhline{0.5\arrayrulewidth}
\end{tabular}
\captionsetup{justification=centering}
{
\caption{\begin{math}\Delta_{|p-q|} = \sum_{p+q}=22\end{math}}
}
}
 \hfill
\parbox{.45\linewidth}{
\begin{tabular}{|c|c|} 
\hline 
 $\Delta$ & $\sum$ \\ 
\Xhline{0.5\arrayrulewidth}
\{-1,0,2,2\} & \{-1,0,2,2\} \\
\Xhline{0.5\arrayrulewidth} 
\end{tabular}
   $dial\ settings$ \\ \\
\begin{tabular}{|c|c|c|c|c|c|c|c|c|c|c|c|c|c|c|}
\hline
$N_{_{\sum_{p+q}}}$ & $od_6{_{_{_{\Delta_{|p-q|}}}}}$  & $\Delta \rightleftharpoons \sum$ & $df_{_{\Delta\sum}}$\\[0.5ex]
\hline
441 & 440 & 1 &  \\
\Xhline{0.5\arrayrulewidth}
437 & 436 & 1 &  0 \\
\Xhline{0.5\arrayrulewidth}
425 & 428 & -3 &  4 \\
\Xhline{0.5\arrayrulewidth}
405 & 420 & -15 &  12 \\
\Xhline{0.5\arrayrulewidth}
377 & 412 & -35 &  20 \\
\Xhline{0.5\arrayrulewidth}
341 & 404 & -63 &  28 \\
\Xhline{0.5\arrayrulewidth}
297 & 396 & -99 &  36 \\
\Xhline{0.5\arrayrulewidth}
245 & 388 & -143 &  44 \\
\Xhline{0.5\arrayrulewidth}
185 & 380 & -195 &  52 \\
\Xhline{0.5\arrayrulewidth}
\end{tabular}
\captionsetup{justification=centering}
{
 \caption{\begin{math}\Delta_{|p-q|} = \sum_{p+q}=42\end{math}} \
}
}
\end{table}
\begin{table}[!htbp] 
 \caption{$\Delta_{|p-q|}=20, p = 2j$ for $j \in \mathbb{N}$}
  \label{table:delta_equilibrim}
\begin{tabular}{|c|c|}
\hline 
-1 & 0 \\ 
\Xhline{0.5\arrayrulewidth}
2 & 2 										 						\\
\Xhline{0.5\arrayrulewidth} 
\end{tabular}
   $dial\ settings$
\\
\begin{tabular}{|c|c|c|c|c|c|c|c|c|c|c|c|c|c|c|}
\hline\hline
$id$ & $p$ & $q$ & $n=pq$ & $\lfloor \sqrt{n} \rfloor$ & $d_1$ & $d_2$ & $d_1^2$ & $d_2^2$ & $od_1$ & $od_2$ & $od_3$ & $od_4$ & $od_6$ \\[0.5ex]
\hline
1 & 2 & 22 & 44 & 6 & 5 & 7 & 25 & 49 & -19 & 5 & 24 & -14 & 9 \\
\Xhline{0.5\arrayrulewidth}
... & ... & ... & ... & ... & ... & ... & ... & ... & ... & ... & ... & ... & ... \\
\Xhline{0.5\arrayrulewidth}
17 & 34 & 54 & 1836 & 42 & 41 & 43 & 1681 & 1849 & -155 & 13 & 168 & -142 & 73 \\
\Xhline{0.5\arrayrulewidth}
18 & 36 & 56 & 2016 & 44 & 43 & 45 & 1849 & 2025 & -167 & 9 & 176 & -158 & 81 \\
\Xhline{0.5\arrayrulewidth}
19 & 38 & 58 & 2204 & 46 & 45 & 47 & 2025 & 2209 & -179 & 5 & 184 & -174 & 89 \\
\Xhline{0.5\arrayrulewidth}
20 & 40 & 60 & 2400 & 48 & 47 & 49 & 2209 & 2401 & -191 & 1 & 192 & -190 & 97 \\
\Xhline{0.5\arrayrulewidth}
21 & 42 & 62 & 2604 & 51 & 51 & 53 & 2601 & 2809 & -3 & 205 & 208 & 202 & \color{blue}\textbf{99} \\
\Xhline{0.5\arrayrulewidth}
22 & 44 & 64 & 2816 & 53 & 53 & 55 & 2809 & 3025 & -7 & 209 & 216 & 202 & \color{blue}\textbf{99} \\
\Xhline{0.5\arrayrulewidth}
23 & 46 & 66 & 3036 & 55 & 55 & 57 & 3025 & 3249 & -11 & 213 & 224 & 202 & \color{blue}\textbf{99} \\
\Xhline{0.5\arrayrulewidth}
... & ... & ... & ... & ... & ... & ... & ... & ... & ... & ... & ... & 202 & \color{blue}\textbf{99} \\
\Xhline{0.5\arrayrulewidth}
$\infty$ & $\infty$ & $\infty$ & $\infty$ & $\infty$ & $\infty$ & $\infty$ & $\infty$ & $\infty$ & $\infty$ & $\infty$ & $\infty$ & 202 & \color{blue}\textbf{99} \\
\Xhline{0.5\arrayrulewidth}
\end{tabular}
\end{table}
\begin{table}[!htbp] 
 \caption{$\sum_{p+q}=20, p = 2j$ for $j \in \mathbb{N}$}
  \label{table:sum_equilibrim}
\begin{tabular}{|c|c|}
\hline 
-1 & 0 \\ 
\Xhline{0.5\arrayrulewidth}
2 & 2 										 						\\
\Xhline{0.5\arrayrulewidth} 
\end{tabular}%
\end{table}
\begin{table}[!htbp] 
 \caption{$\Delta_{|p=q|}=22, p = 2j$ for $j \in \mathbb{N}$}
  \label{table:delta_4k2}
\begin{tabular}{|c|c|}
\hline 
0 & -1 \\ 
\Xhline{0.5\arrayrulewidth}
2 & 2 										 						\\
\Xhline{0.5\arrayrulewidth} 
\end{tabular}%
\\
\end{table}
\begin{table}[!htbp] 
 \caption{$\sum_{p+q}=22, p = 2j$ for $j \in \mathbb{N}$}
  \label{table:sum_4k2}
\begin{tabular}{|c|c|}
\hline 
-1 & 0 \\ 
\Xhline{0.5\arrayrulewidth}
2 & 2 										 						\\
\Xhline{0.5\arrayrulewidth} 
\end{tabular}%
\end{table}
\begin{table}[!htbp] 
 \caption{$\Delta_{|p-q|}=23, p = 2j+1$ for $j \in \mathbb{N}$}
  \label{table:odd_delta_1}
\begin{tabular}{|c|c|} 
\hline 
-1 & 0 \\ 
\Xhline{0.5\arrayrulewidth}
2 & 2 										 						\\
\Xhline{0.5\arrayrulewidth}
\end{tabular}%
\vspace{30em}
\end{table} 
\begin{table}[!htbp] 
 \caption{$\Delta_{|p-q|}=25, p = 2j+1 for j \in \mathbb{N}$}
  \label{table:odd_delta_2}
\begin{tabular}{|c|c|} 
\hline 
0 & -1 \\ 
\Xhline{0.5\arrayrulewidth}
2 & 2 \\
\Xhline{0.5\arrayrulewidth}
\end{tabular}%
\vspace{20em}
\end{table}
\begin{table}[!htbp] 
 \caption{$\Delta_{|p-q|}=25, p = 2j$ for $j \in \mathbb{N}$}
  \label{table:odd_delta_3}
\begin{tabular}{|c|c|} 
\hline 
-1 & 0 \\ 
\Xhline{0.5\arrayrulewidth}
2 & 2 \\
\Xhline{0.5\arrayrulewidth}
\end{tabular}%
\vspace{20em}
\end{table}
\begin{table}[!htbp] 
 \caption{$\Delta_{|p-q|}=23, p = 2j$ for $j \in \mathbb{N}$}
  \label{table:odd_delta_1}
%
\vspace{15em}
\end{table} 
\section {Figures} 
\begin{figure}[htbp!] 
\captionsetup{justification=centering}
{
  \includegraphics[scale=1]{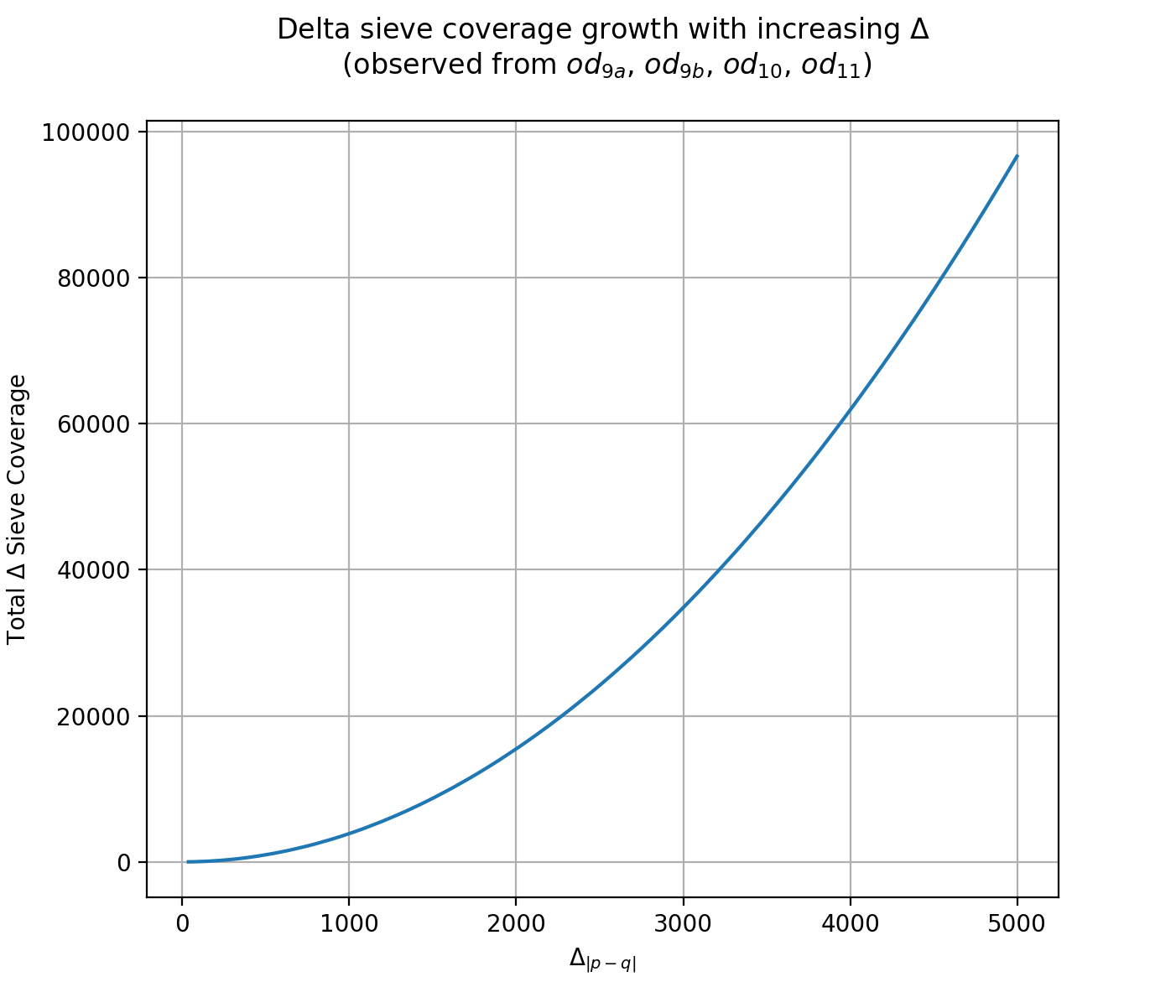}  
  \caption{$\Delta$ Sieve Coverage Growth}
}
  \label{fig:dsc_growth}
\end{figure}
\begin{figure}[htbp!] 
{
\captionsetup{justification=centering}
{
  \includegraphics[scale=1]{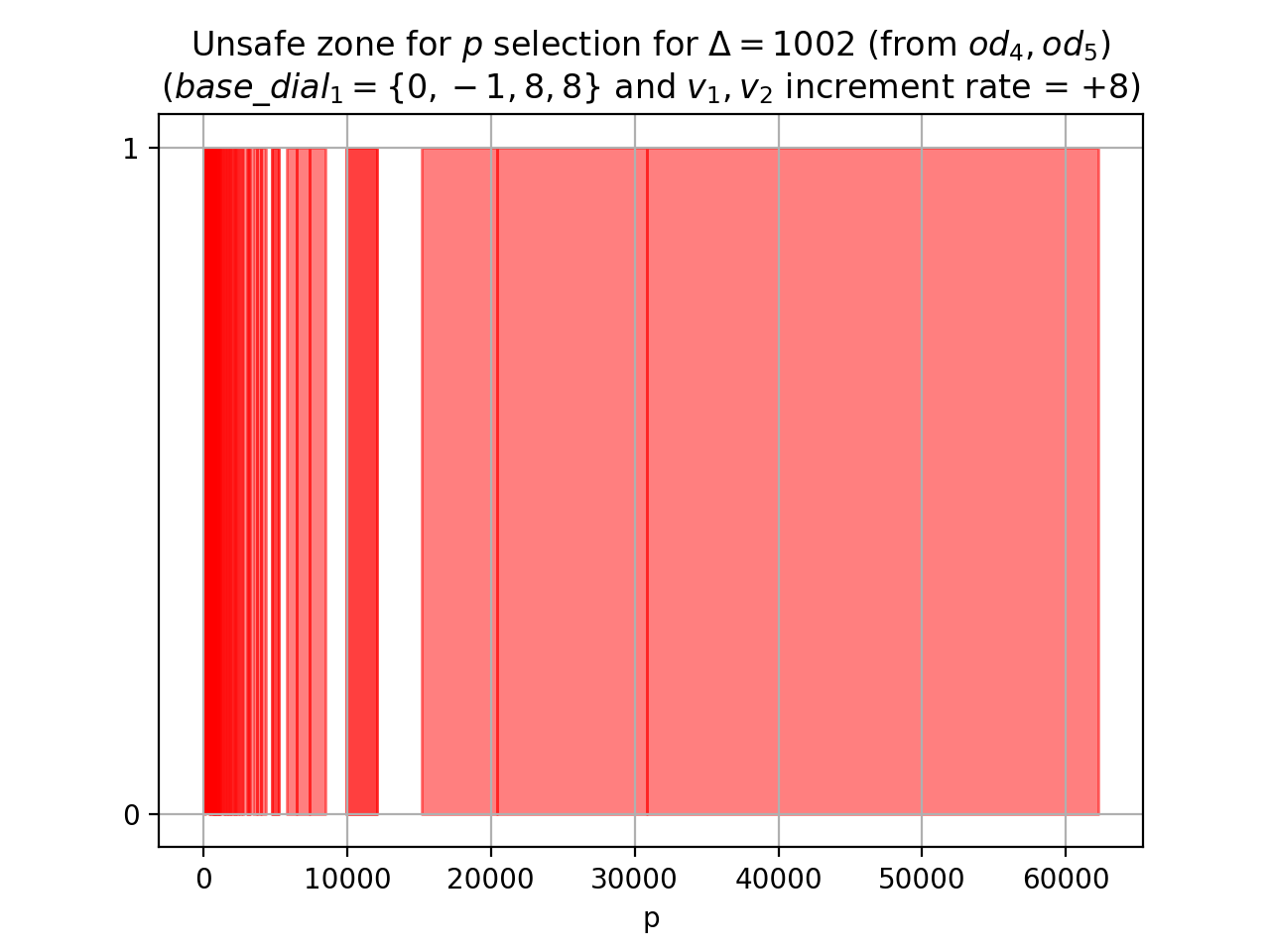}  
  \caption{RSA unsafe zone for p selection from $od_4, od_5$}
}
  \label{fig:rsa_unsafe_zone}}
\end{figure} 
\begin{figure}[htbp!] 
{
\captionsetup{justification=centering}
{
  \includegraphics[scale=1]{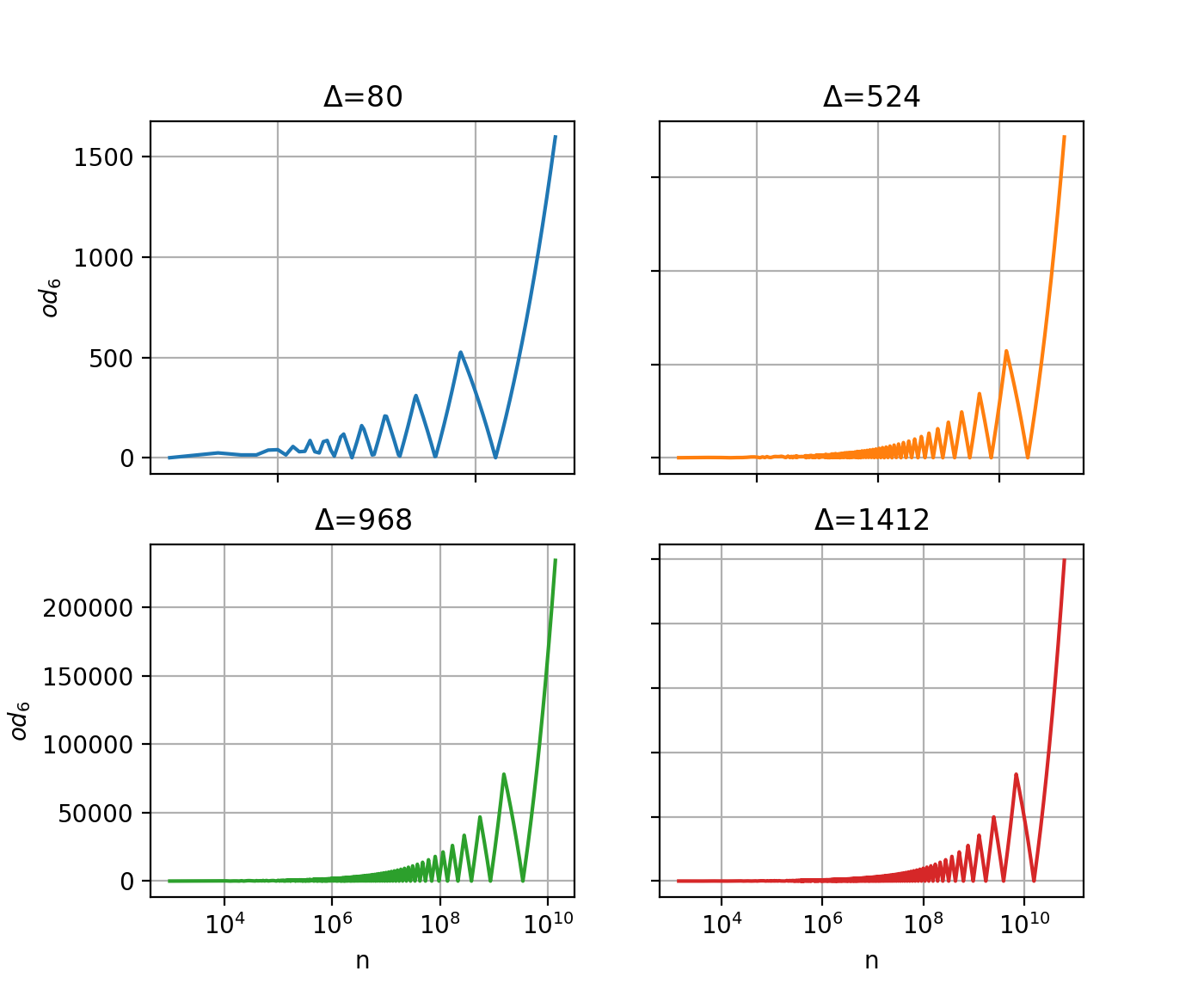}
  \caption{\begin{math}A_+\end{math} graph where \begin{math}v_1 = v_2 = 2\end{math} for \begin{math}\Delta=4k, k \in \mathbb{N}\end{math} form}
}
  \label{fig:a_plus_4k}}
\end{figure}
\begin{figure}[htbp!] 
{
\captionsetup{justification=centering}
{
  \includegraphics[scale=1]{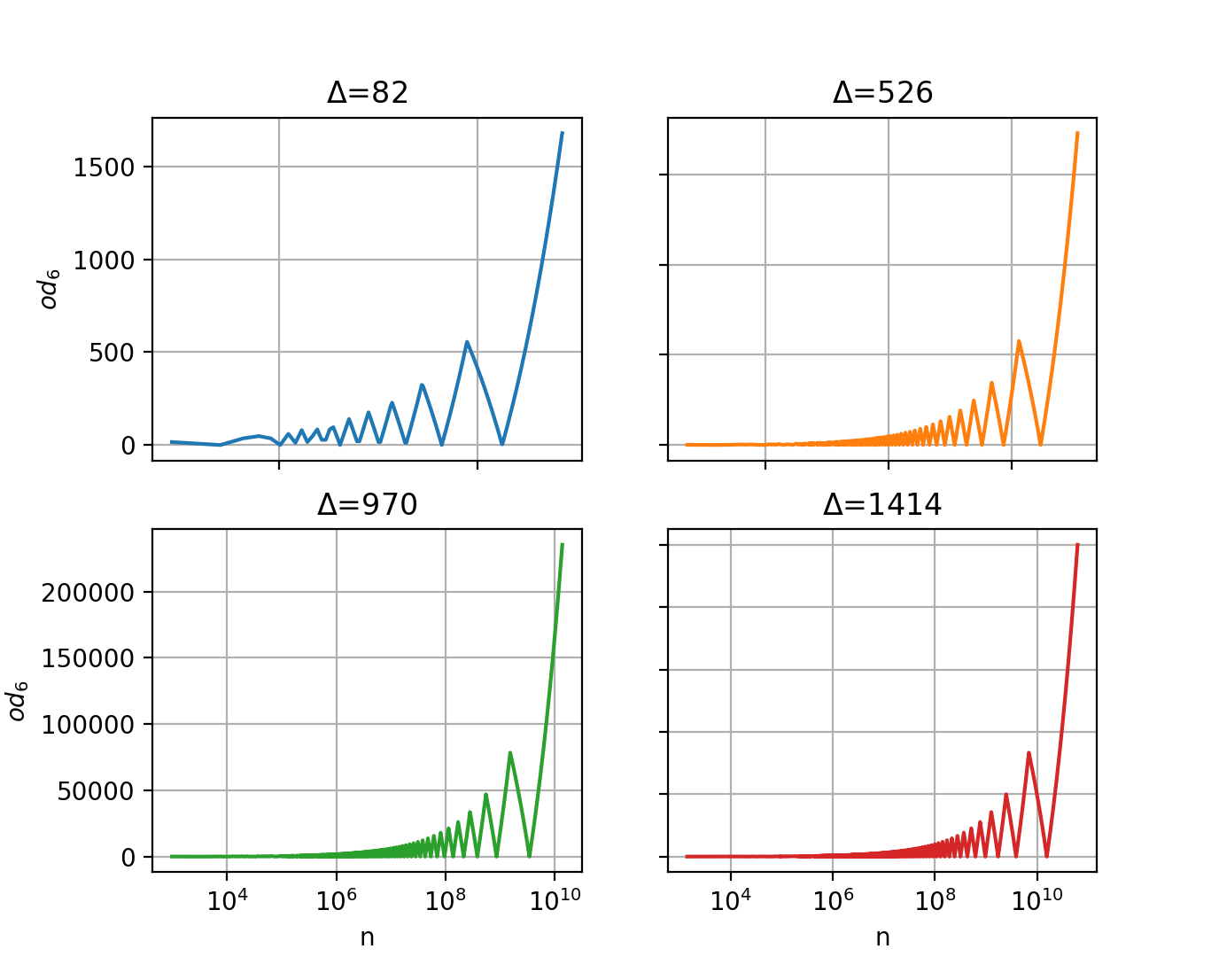}
  \caption{\begin{math}A_+\end{math} graph where \begin{math}v_1 = v_2 = 2\end{math} for \begin{math}\Delta=4k+2, k \in \mathbb{N}\end{math} form}
}
  \label{fig:a_plus_4k2}}
\end{figure}
\begin{figure}[htbp!] 
{
\captionsetup{justification=centering}
{
  \includegraphics[scale=1]{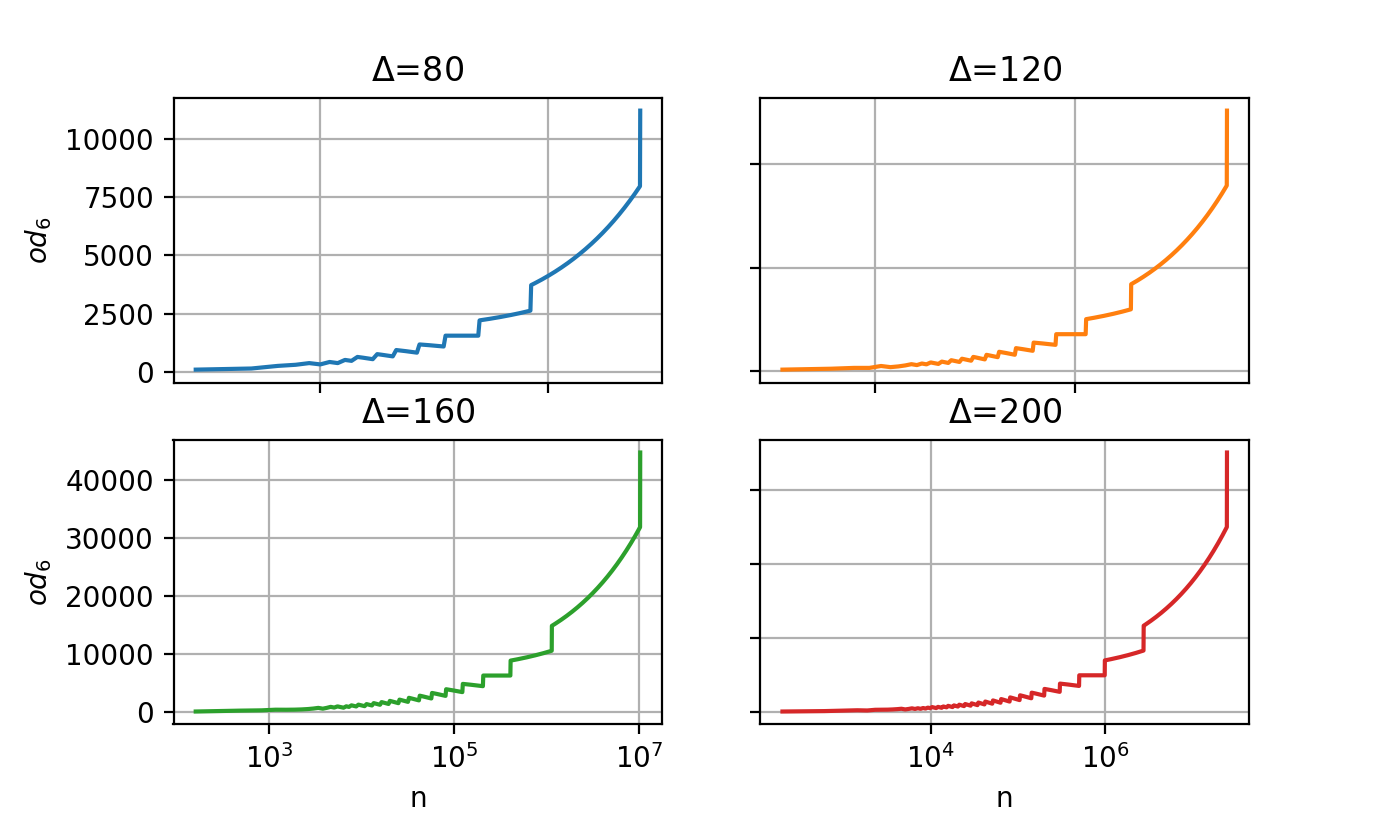}
  \caption{\begin{math}A_+\end{math} graph where \begin{math}v_1 = v_2 = 14\end{math} for \begin{math}\Delta=4k, k \in \mathbb{N}\end{math} form}
}
  \label{fig:roxy_v1_14_2}}
\end{figure}
\begin{figure}[htbp!] 
{
\captionsetup{justification=centering}
{
  \includegraphics[scale=1]{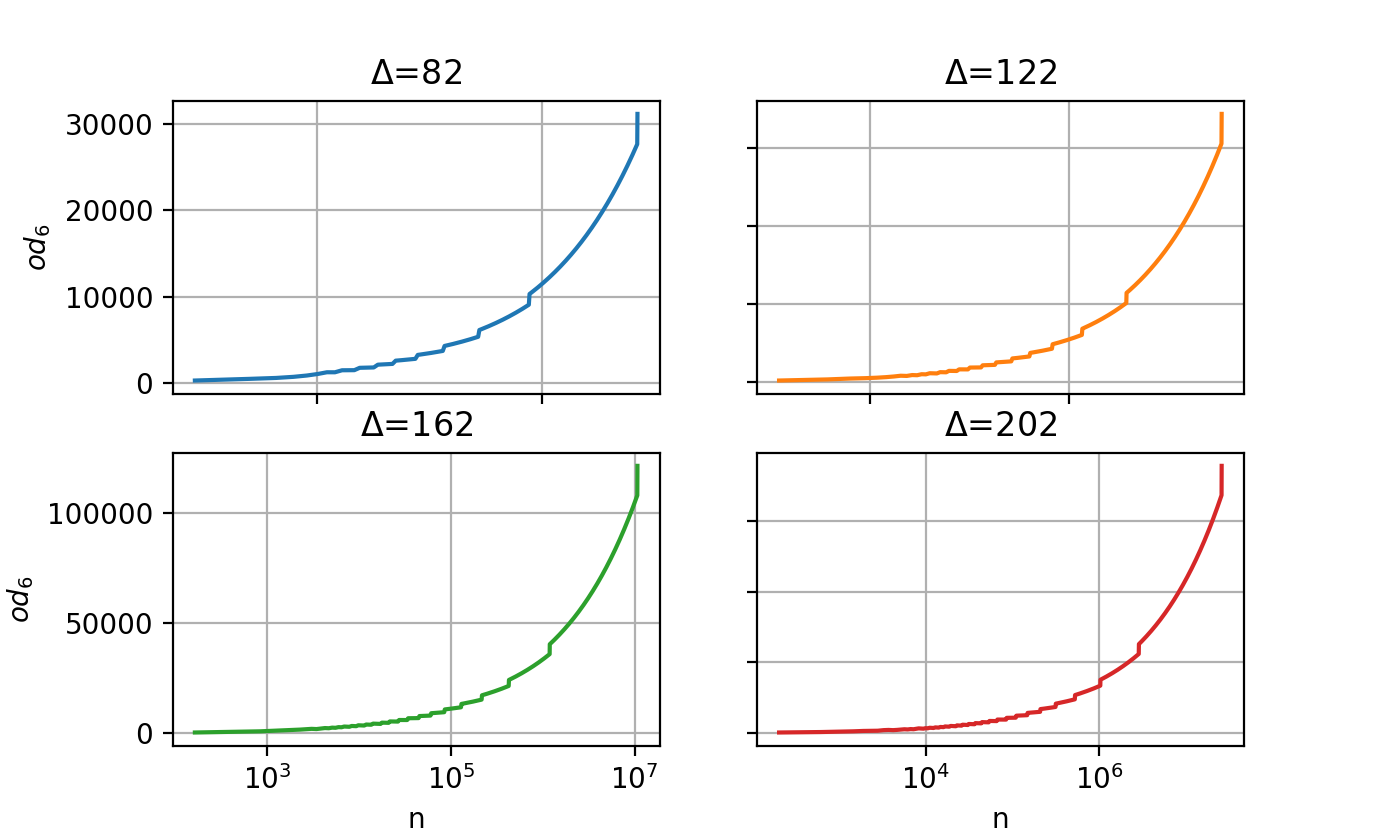}
  \caption{\begin{math}A_+\end{math} graph where \begin{math}v_1 = v_2 = 37\end{math} for \begin{math}\Delta=4k+2, k \in \mathbb{N}\end{math} form}
}
  \label{fig:roxy_v1_37_2}}
\end{figure}
\begin{figure}[htbp!] 
{
\captionsetup{justification=centering}
{
  \includegraphics[scale=1]{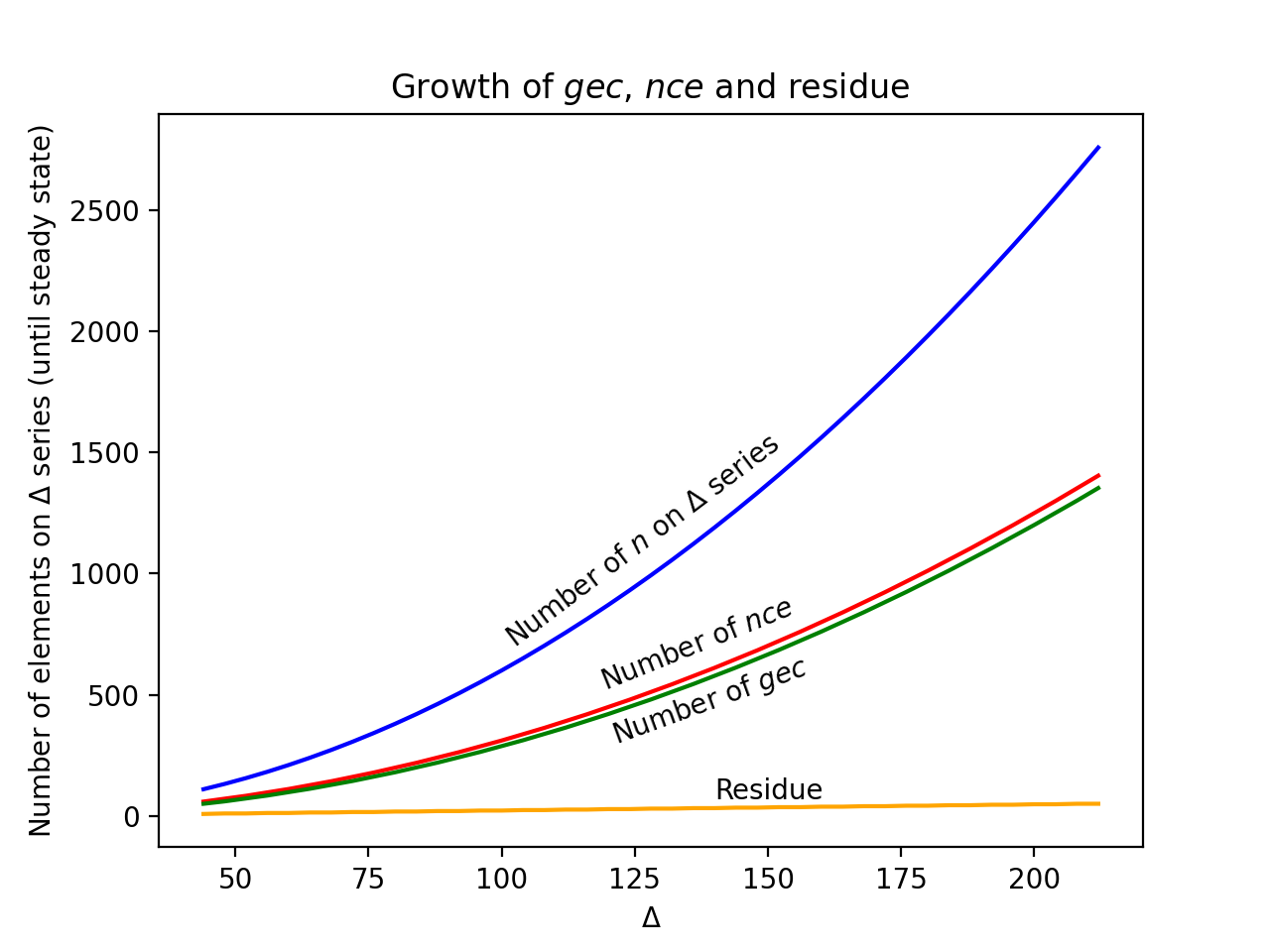}
  \caption{Growth of $gec, nce, residue$ when $\Delta=4k, k \in \mathbb{N}$ form} 
}
  \label{fig:growth_of_gec_2}
}
\end{figure}
\begin{figure}[htbp!] 
{
\captionsetup{justification=centering}
{
  \includegraphics[scale=1]{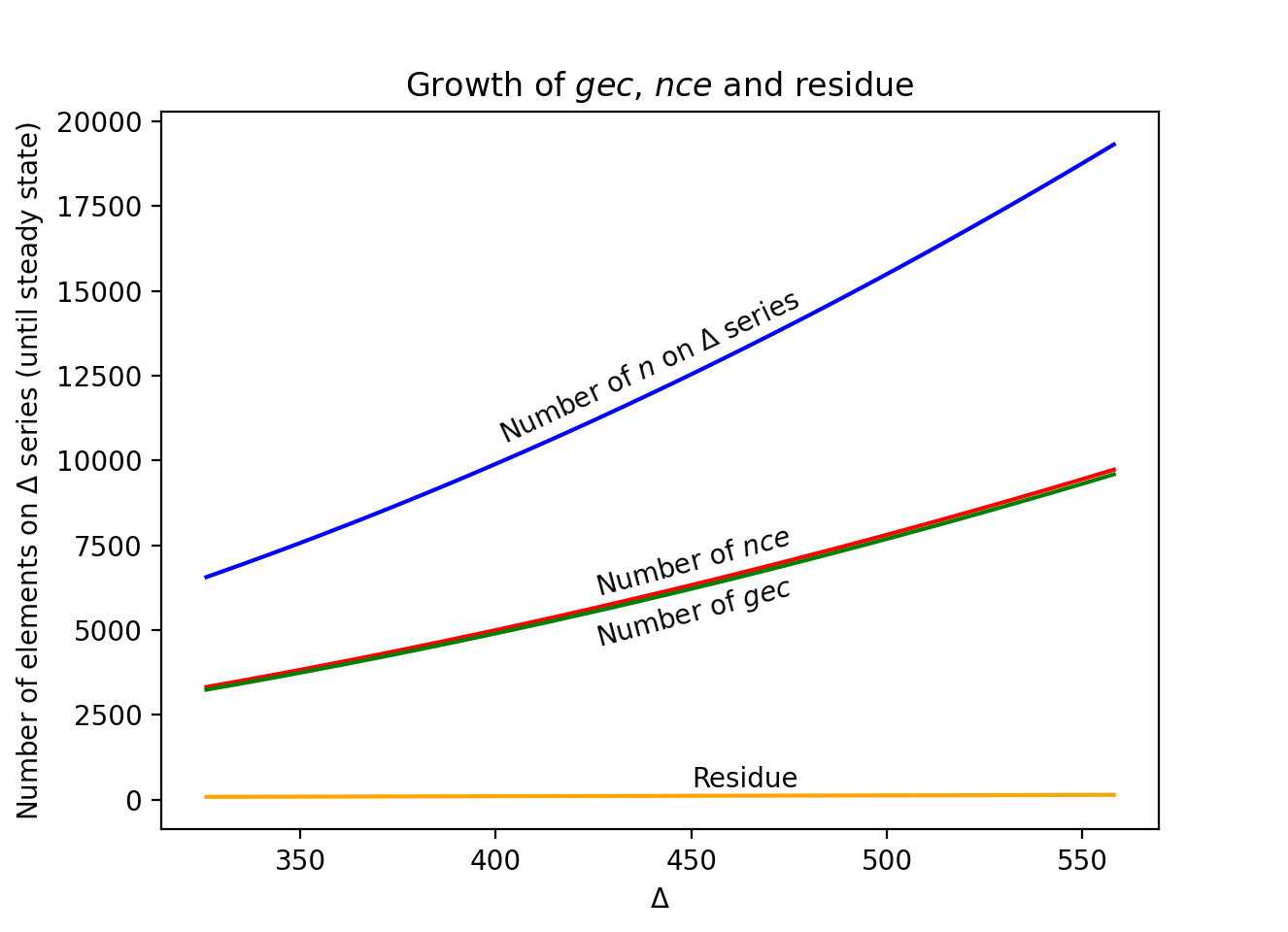}
  \caption{Growth of $gec, nce, residue$ when $\Delta=4k+2, k \in \mathbb{N}$ form} 
}
  \label{fig:growth_of_gec_4k2}
}
\end{figure}
\begin{figure}[htbp!] 
{
\captionsetup{justification=centering}
{
 \begin{center}
  \includegraphics[scale=1]{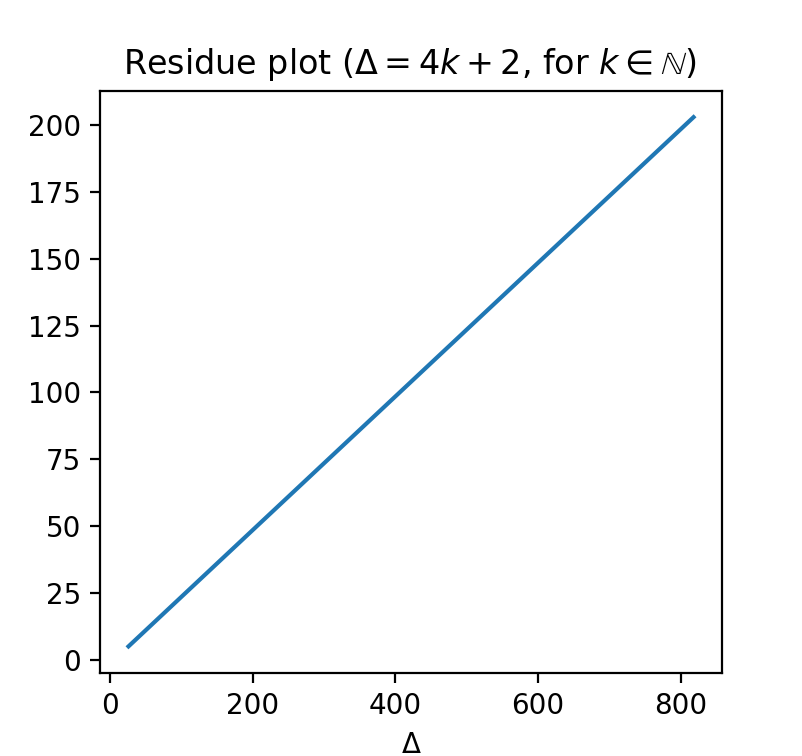}
 \end{center}
 \caption{The linear residue}
}
  \label{fig:residue_plot_4k2_2}
}
\end{figure}
\begin{figure}[htbp!] 
{
\captionsetup{justification=centering}
{
 \begin{center}
  \includegraphics[scale=1]{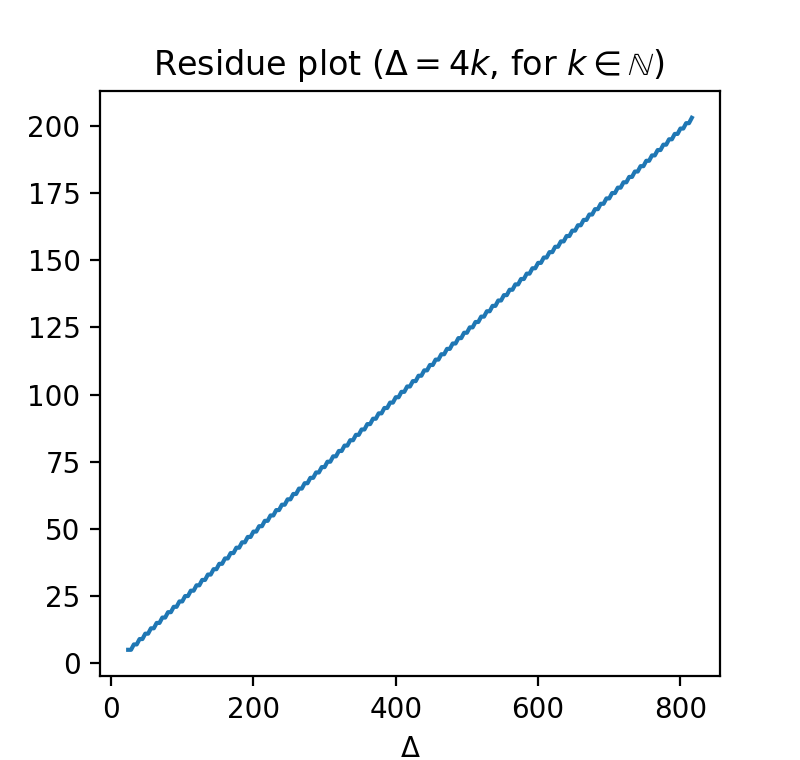}
 \end{center}
  \caption{The wood saw residue}
}
  \label{fig:residue_plot_4k_6}
}
\end{figure}
\section{Code to validate the hypotheses and observations} 
Github repo: https://github.com/vdeltasieve/arjun
\begin{enumerate}
\item \textbf{client\_ssv\_p\_check.py}: The code confirms the hypotheses in the following equations: (\ref{eqn:ssv1a}), (\ref{eqn:ssv1b}), (\ref{eqn:ssv2a}), (\ref{eqn:ssv2b}), (\ref{eqn:ssv2c}), (\ref{eqn:ssv2d}), (\ref{eqn:ssv2e}), (\ref{eqn:ssv4a}), (\ref{eqn:ssv4b}), (\ref{eqn:ssv4c}), (\ref{eqn:ssv4d}), (\ref{eqn:ssv4e})
\item \textbf{client\_od\_constant\_check.py}: The code confirms the hypothesis that once a $\Delta$ sieve zone is tuned for some known $\Delta$ in a particular observation deck, this $\Delta$ sieve zone will remain constant for all $\Delta$s. This code specifically confirms the $\Delta$ sieve zone results as per equations: (\ref{eqn:ssvgenod9a}), (\ref{eqn:ssvgenod9b}), (\ref{eqn:ssvgenod10}), (\ref{eqn:ssvgenod11})
\item \textbf{client\_od4\_dial\_rotation\_check.py}: The code confirms the hypothesis that $\Delta$ sieve zones will shift with changing $v_1$, $v_2$. This code specifically confirms this hypothesis for $od_4$
\item \textbf{client\_od5\_dial\_rotation\_check.py}: The code confirms the hypothesis that $\Delta$ sieve zones will shift with changing $v_1$, $v_2$. This code specifically confirms this hypothesis for $od_5$
\item \textbf{client\_roxy.py}: The code confirms observations in subsection 5.3 - "Reflection over \{X,Y\}"
\item \textbf{client\_delta\_sum\_equilibrium.py}: The code confirms the connection of $od_6$ on $\Delta$ series with $N$ on $\sum$ series, as described in section 8 - "The equilibrium of $\Delta_{|p-q|}$ and $\sum_{p+q}$"
\item \textbf{client\_dsc\_trapdoor.py}: The code validates the algorithm and confirms the example given in section 9 - "The Trapdoor"
\item \textbf{client\_odd\_delta\_ssv\_check.py}: The code confirms the hypothesis in the following equations: (\ref{eqn:odd_delta_1}), (\ref{eqn:odd_delta_2}), (\ref{eqn:odd_delta_3}), (\ref{eqn:odd_delta_4})
\end{enumerate} \
\clearpage\
\bibliography{references} 
\bibliographystyle{ieeetr}
\end{document}